\documentclass[a4paper,12pt]{amsart}
\usepackage{amsfonts}
\usepackage{amssymb}
\usepackage{ifthen}
\usepackage{graphicx}
\nonstopmode \numberwithin{equation}{section}
\usepackage[usenames]{color}
\newtheorem{thm}{Theorem}[section]
\newtheorem{lem}{Lemma}[section]
\newtheorem{op}{Open Problem}[section]
\newtheorem{cor}{Corollary}[section]
\newtheorem{prop}{Proposition}[section]

\newtheorem{cl}{Claim}[section]
\newtheorem{ca}[section]{Case}
\newtheorem{sca}[section]{Subcase}
\newtheorem{scl}[section]{Subclaim}
\newtheorem{conj}[equation]{Conjecture}

\theoremstyle{definition}
\newtheorem{defn}{Definition}[section]

\newtheorem{ques}[equation]{Question}
\newtheorem{rem}{Remark}[section]
\newtheorem{exam}{Example}[section]

\newcounter {own}
\def\theown {\thesection       .\arabic{own}}

\newenvironment{pf}[1][]{%
 \vskip 3mm
 \noindent
 \ifthenelse{\equal{#1}{}}%
  {{\slshape Proof. }}%
  {{\slshape #1.} }%
 }%
{\qed\bigskip}

\newcounter{alphabet}
\newcounter{tmp}
\newenvironment{Thm}[1][]{\refstepcounter{alphabet}%
\bigskip%
\noindent%
{\bf Theorem \Alph{alphabet}}%
\ifthenelse{\equal{#1}{}}{}{ (#1)}%
{\bf .} \itshape}{\vskip 8pt}

\makeatletter
\newcommand{\Ref}[1]{\@ifundefined{r@#1}{}{\setcounter{tmp}{\ref{#1}}\Alph{tmp}}}
\makeatother
\newcommand{\IR}{{\mathbb R}}

\newcommand{\ID}{{\mathbb D}}
\newcommand{\IB}{{\mathbb B}}

\newcommand{\diam}{{\operatorname{diam}}}

\newcommand{\dist}{{\operatorname{dist}}}

\def\be{\begin{equation}}
\def\ee{\end{equation}}

\newcommand{\bee}{\begin{enumerate}}
\newcommand{\eee}{\end{enumerate}}

\newcommand{\blem}{\begin{lem}}
\newcommand{\elem}{\end{lem}}
\newcommand{\bthm}{\begin{thm}}
\newcommand{\ethm}{\end{thm}}
\newcommand{\bcor}{\begin{cor}}
\newcommand{\ecor}{\end{cor}}
\newcommand{\beg}{\begin{exam}}
\newcommand{\eeg}{\end{exam}}
\newcommand{\begs}{\begin{examples}}
\newcommand{\eegs}{\end{examples}}
\newcommand{\bdefe}{\begin{defn}}
\newcommand{\edefe}{\end{defn}}
\newcommand{\bprob}{\begin{prob}}
\newcommand{\eprob}{\end{prob}}
\newcommand{\bques}{\begin{ques}}
\newcommand{\eques}{\end{ques}}
\newcommand{\bei}{\begin{itemize}}
\newcommand{\eei}{\end{itemize}}
\newcommand{\bcon}{\begin{conj}}
\newcommand{\econ}{\end{conj}}
\newcommand{\bop}{\begin{op}}
\newcommand{\eop}{\end{op}}

\newcommand{\bca}{\begin{ca}}
\newcommand{\eca}{\end{ca}}
\newcommand{\bsca}{\begin{sca}}
\newcommand{\esca}{\end{sca}}

\newcommand{\bcl}{\begin{cl}}
\newcommand{\ecl}{\end{cl}}

\newcommand{\bscl}{\begin{scl}}
\newcommand{\escl}{\end{scl}}

\newcommand{\bcons}{\begin{conjs}}
\newcommand{\econs}{\end{conjs}}
\newcommand{\bprop}{\begin{propop}}
\newcommand{\eprop}{\end{propop}}
\newcommand{\br}{\begin{rem}}
\newcommand{\er}{\end{rem}}
\newcommand{\brs}{\begin{rems}}
\newcommand{\ers}{\end{rems}}
\newcommand{\bo}{\begin{obser}}
\newcommand{\eo}{\end{obser}}
\newcommand{\bos}{\begin{obsers}}
\newcommand{\eos}{\end{obsers}}
\newcommand{\bpf}{\begin{pf}}
\newcommand{\epf}{\end{pf}}
\newcommand{\ba}{\begin{array}}
\newcommand{\ea}{\end{array}}
\newcommand{\beq}{\begin{eqnarray}}
\newcommand{\beqq}{\begin{eqnarray*}}
\newcommand{\eeq}{\end{eqnarray}}
\newcommand{\eeqq}{\end{eqnarray*}}

\newcommand{\ds}{\displaystyle}

\newcounter{minutes}
\divide\time by 60
\newcounter{hours}
\multiply\time by 60 \addtocounter{minutes}{-\time}
\begin{document}

\bibliographystyle{amsplain}
\title{The quasiconformal subinvariance property of John domains in $\protect \IR^n$ and its application}

%=========================================================================
\thanks{%$^\dagger$
File:~\jobname .tex,
          printed: \number\year-\number\month-\number\day,
          \thehours.\ifnum\theminutes<10{0}\fi\theminutes}
%=========================================================================

\author{M. Huang}
\address{M. Huang, Department of Mathematics,
Hunan Normal University, Changsha,  Hunan 410081, People's Republic
of China.
% and Department of Mathematics,
%Indian Institute of Technology Madras, Chennai-600 036, India
}
\email{manzihuang0208@gmail.com}

\author{Y. Li}
\address{Y. Li, Department of Mathematics,
Hunan Normal University, Changsha,  Hunan 410081, People's Republic
of China.} \email{yaxiangli@163.com}

\author{S. Ponnusamy$^\dagger $}
\address{S. Ponnusamy,
Indian Statistical Institute (ISI), Chennai Centre, SETS (Society
for Electronic Transactions and security), MGR Knowledge City, CIT
Campus, Taramani, Chennai 600 113, India. }
\email{samy@isichennai.res.in, samy@iitm.ac.in}
%\address{S. Ponnusamy, Department of Mathematics,
%Indian Institute of Technology Madras, Chennai-600 036, India.}
%\email{samy@iitm.ac.in}

\author{X. Wang $^\dagger {}^\dagger$
%${}^{~\mathbf{*}}$
}
\address{X. Wang, Department of Mathematics,
Hunan Normal University, Changsha,  Hunan 410081, People's Republic
of China.}
\email{xtwang@hunnu.edu.cn}

%\date{}
\keywords{Uniform domain, QED domain, broad domain, John domain, quasiconformal mapping, quasisymmetry.\\
$^\dagger$ {\tt This author is on leave from the
%currently at the
%Indian Statistical Institute (ISI), Chennai Centre, SETS (Society
%for Electronic Transactions and security), MGR Knowledge City, CIT
%Campus, Taramani, Chennai 600 113, India. {\it E-mail address}: samy@isichennai.res.in}\\
Department of Mathematics,
Indian Institute of Technology Madras, Chennai-600 036, India.\\
${}^\dagger {}^\dagger$ {\tt Corresponding author.}\\
}
}

\begin{abstract}
The main aim of this paper is to give a complete solution to one of
the open problems, raised by Heinonen from 1989, concerning the
subinvariance of John domains under quasiconformal mappings in
$\IR^n$. As application, the quasisymmetry of quasiconformal
mappings is discussed.
\end{abstract}

%\thanks{The research was partly supported by NSFs of China (No. 11071063 and No. 11101138).
%}

\maketitle \pagestyle{myheadings} \markboth{M. Huang, Y. Li, S.
Ponnusamy and X. Wang}{The quasiconformal subinvariance property of
John domains in $\protect \IR^n$ and its application}

\section{Introduction and main results}\label{sec-1}

In this paper, we study the subinvariance property of John domains
under quasiconformal mappings in $\IR^n$.
%This means that under
%certain conditions, $f$ preserves some property of subdomains of
%$D$.
%For example, if $D'$ is a $c$-QED domain, $f$ maps every
%$c$-QED domain $G\subset D$ onto a $c_1$-QED domain with $c_1=c_1(n,
%K, c)$ (see \cite{FHM}, p. 121). there is some constant $c\geq 1$
%depending only on the given coefficients such that for each John
%subdomain of $D$, its image under $f$ is a $c$-John subdomain in
%$D'$, see e.g. \cite{GNV}.
The motivation for this study stems from
one of the open problems raised by Heinonen in \cite{H}. The aim of
this paper is to give a complete solution to this open problem and
its application.
%Before the statement of the open problem and its
%background information, let us recall the concept of John domains.
John domains were introduced by John in \cite{Jo} in his study of
rigidity of local isometries. The term ``John domain" was coined by
Martio and Sarvas seventeen years later \cite{MS}. Among various
equivalent characterizations, we shall adopt the definition for John
domains in the sense of carrot property, which is called the carrot
definition for John domains in the following. We shall give the
precise definition together with other necessary notions and notations
in the next section.

Throughout the paper we assume that $D$ and $D'$ are subdomains in $\IR^n$ and
that $f: D\to D'$ is a $K$-quasiconformal mapping with
$K\geq 1$. See \cite{Vai, Vu-book-88} for definitions and properties
of $K$-quasiconformal mappings.

In \cite{FHM}, the authors showed the following.

\begin{Thm}\label{ThmB}$($\cite[p.~120-121]{FHM}$)$
Suppose $D'$ is $QED$. Then for each $QED$ subdomain $D_1$ in $D$,
$f(D_1)$ is also $QED$.
\end{Thm}

The reader is referred to \cite{Gem} for the definition of QED
domains. By Theorem \Ref{ThmB} and \cite[Theorem 5.6]{Vai2}, the
following is obvious.

\begin{Thm}\label{ThmA}
If $D'$ is $c$-uniform, then
for each $c_1$-uniform subdomain $D_1$ in $D$, $f(D_1)$ is still
$\rho$-uniform, where $\rho=\rho(n, K, c, c_1)$, which means that
the constant $\rho$ depends only on $n$, $K$, $c$ and $c_1$.
\end{Thm}

In the next result, the case $(1)$ is due to V\"ais\"al\"a from
\cite[Theorem 2.20]{Vai0}, whereas the case $(2)$ is obtained by
Heinonen \cite[Theorem 7.1]{H}.

\begin{Thm}\label{ThmC}
Suppose $D'$ is broad.
\begin{enumerate}
\item[{\rm (1)}]
Then for each John subdomain $D_1$ in $D$, $f(D_1)$ is a John domain;

\item[{\rm (2)}]
If both $D$ and $D'$ are bounded, then for each broad subdomain
$D_1$ in $D$, $f(D_1)$ is a broad domain.
\end{enumerate}
\end{Thm}

We refer to \cite{MS} for an early discussion on this topic.
However, a natural problem is that whether $f(D_1)$ is a John domain
for each John subdomain $D_1$ of $D$ when $D'$ is a John domain. This is an open problem raised by Heinonen \cite{H} (in fact, this open problem was put forward by V\"ais\"al\"a) in the
following form.

\bop\label{Con2}
Suppose that  $f$ is a quasiconformal mapping of
a domain $D$ in $\mathbb{R}^n$ onto a John domain $D'$ in $\mathbb{R}^n$. Is it then
true that every John subdomain of $D$ is mapped onto a John subdomain of $D'$ by $f$?
\eop

That is, if $D'$ is a $c$-John domain, then for every $c_1$-John
subdomain in $D$, is its image under $f$ a $c_2$-John domain in
$D'$, where $c_2=c_2(n, K, c, c_1)$ (see e.g. \cite{GNV} or
\cite{FHM})? If the answer is yes, in what follows, we shall say
that $f$ has the subinvariance property of John domains. Heinonen
himself discussed this problem and as consequence he obtained that for a quasiconformal mapping
from the unit ball $\mathbb{B}$ in $\IR^n$ onto a John domain $D'$, the image of each Stolz
cone $C_M(w)$ with vertex at $w\in \partial\mathbb{B}$, the boundary of $\mathbb{B}$, under $f$ is uniform, see \cite[Theorem 7.3]{H}, where
$C_M(w)$ is defined to be the interior of the closed
convex hull of $w$ and the hyperbolic ball centered at $0$ with
radius $M>0$. It is known that every ball in $\IR^n$ is a $\frac{\pi}{2}$-uniform domain.

In this paper, we consider Open Problem \ref{Con2} further. The
examples constructed in Section \ref{sec-3} show that the answer to Open problem \ref{Con2}
is negative when one of $D$ and $D'$ is unbounded.
This observation shows that it suffices for us to consider the case where
both $D$ and $D'$ are bounded. In this case, we get the following affirmative answer.

\begin{thm}\label{thm1.1}
Suppose that $D$ and $D'$ are bounded subdomains in $\IR^n$ and that $f: D\to D'$ is a
$K$-quasiconformal mapping. If
$D'$ is an $a$-John domain with center $y'_0$, then for every $c$-John subdomain $D_1$ in $D$ with center $z_0$, its image $f(D_1)$
 is a $\tau$-John subdomain in $D'$ with center $f(z_0)$, where $\tau=\tau\left(n, K, a, c, \frac{\diam(D)}{d_D(f^{-1}(y'_0))}\right)$.
\end{thm}

In \cite{Vai0}, V\"{a}is\"{a}l\"{a} considered the question: When is
$G\times \IR^1$ quasiconformally equivalent to the unit open ball
$\mathbb{B}^3$ in $\IR^3$, where $G$ is a domain in the plane $\IR^2$? He
proved that there is a quasiconformal mapping $f$ from $G\times \IR^1$ to $\mathbb{B}^3$ if
and only if $G$ satisfies the so-called internal chord-arc
condition, see \cite[Theorem 5.2]{Vai0}. In the proof of
\cite[Theorem 5.2]{Vai0}, the following result plays a key role.

\begin{Thm}\label{ThmE}$($\cite[Theorem 2.20]{Vai0}$)$
Suppose that $f: D\to D'$ is a $K$-quasiconformal mapping between
domains $D$ and $D'\subset \IR^n$, where $D$ is $\varphi$-broad.
Suppose also that $A\subset D$ is a pathwise connected set and that
$A'=f(A)$ has the $c_1$-carrot property in $D'$ with center $y'_0\in \overline{D'}$. If
$y'_0 \not=\infty$ and hence $y'_0\in D'$, we assume that
$\diam(A)\leq c_2d_{D}(f^{-1}(y'_0))$. If $y'_0=\infty$, we assume that $f$
extends to a homeomorphism $D\cup \{\infty\}\to D'\cup \{\infty\}$.
Then $f|_{A}$ is $\eta$-quasisymmetric in the metrics $\delta_D$ and
$\delta_{D'}$ with $\eta$ depending only on the data $\kappa_1=\kappa_1(n,
K, c_1,c_2,\varphi)$.
\end{Thm}

We remark that in Theorem \Ref{ThmE} and elsewhere later, the closure of a domain in $\IR^n$
is always taken from $\overline{\IR}^n$.

%the assumption ``$\diam(A)\leq c_2d_{D}(y_0)$" in
%Theorem \Ref{ThmE} means that ``$A$ is bounded" and ``$\diam(A)\leq
%c_2d_{D}(y_0)$".
V\"{a}is\"{a}l\"{a} pointed out in \cite{Vai0} that
this sufficient condition for
 a quasiconformal mapping to be quasisymmetric in the internal metric is of independent interest.
In order to prove the main result in \cite{H}, Heinonen discussed
the generalization of Theorem \Ref{ThmE}. He considered the case
where both $D$ and $D'$ are bounded. In this case, he proved

\begin{Thm}\label{Lem6-1'}$($\cite[Theorem 6.1]{H}$)$
Suppose that $D$ and $D'$ are bounded, that $f: D\to D'$ is
$K$-quasiconformal, and that $D$ is $\varphi$-broad. If $A\subset D$
is such that $f(A)$ is $b$-$LLC_2$ with respect to $\delta_{D'}$ in
$D'$, then $f|_A: A\to f(A)$ is weakly $H$-quasisymmetric in the
metrics $\delta_{D}$ and $\delta_{D'}$, where $H$ depends only on the data
$$\kappa_2=\kappa_2\left (n,K,b,\varphi, \frac{\delta_{D}(A)}{d_{D}(x_0)},
\frac{\delta_{D'}(f(A))}{d_{D'}(f(x_0))}\right ),
$$
$x_0$ is a fixed point in $A$ and $\delta_{D}(A)$ denotes the
$\delta_{D}$-diameter of $A$.
\end{Thm}

As an application of Theorem \ref{thm1.1}, we consider Theorem
\Ref{ThmE} further. Our result is as follows.

%Also under this
%replacement, our result overcomes the mentioned drawback in Theorem
%\Ref{Lem6-1'}.

%\begin{thm}\label{thm1.2} Suppose that $f: D\to D'$ is a $K$-QC map between  domains $D$,
%$D'\subset \IR^n$, where $D$ is $\varphi$-broad and $D'$ is bounded.
%Suppose also that $A\subset D$ is a pathwise connected set and that
%$A'=f(A)$ has the $c_1$-carrot property in $D'$ with center $y'_0\in
%\overline{D'}$. Then $f|_{A}$ is $\eta$-QUASISYMMETRIC in the metrics $\delta_D$
%and $\delta_{D'}$ with $\eta$ depending only on the data
%$\nu=(c_1,c_2,K,\varphi,n)$.
%\end{thm}

\begin{thm}\label{thm1.2}
Suppose that $f: D\to D'$ is a $K$-quasiconformal mapping between bounded
subdomains $D$ and $D'$ in $\IR^n$, where $D$ is a $b$-uniform
domain. Suppose also that $A\subset D$ is a pathwise connected set
and that $A'=f(A)$ has the $c_1$-carrot property in $D'$ with center
$y'_0$. Then the restriction $f|_{A}$ is
$\eta$-quasisymmetric in the metrics $\delta_D$ and $\delta_{D'}$
with $\eta$ depending only on the data $\kappa=\kappa\Big(n, K,
b,c_1,\frac{\diam(D')}{d_{D'}(f(x_0))}\Big)$, where $x_0\in D$ satisfies $d_D(x_0)=\sup\{d_D(x):\, x\in D\}$.
\end{thm}
%satisfies that $D$ has the $4b^2$-carrot property with center $x_0$.\end{thm}

%$(2)$ If  $D$ is unbounded and $f$ extends to a homeomorphism $D\cup
%\{\infty\}\to D'\cup \{\infty\}$, then the restriction $f|_{A}$ is
%$\eta$-quasisymmetric in the metrics $\delta_D$ and $\delta_{D'}$
%with $\eta$ depending only on the data $\kappa=(n, K, b,c_1)$.
%\end{thm}

\br
\noindent $(1)$ Although, in Theorem \ref{thm1.2}, we replace the condition ``$D$ being broad" in
 Theorem \Ref{ThmE} by the one ``$D$ being uniform", the dependence of the function $\eta$ on the center of $A'$
 in Theorem \Ref{ThmE} is removed.

\noindent $(2)$ Obviously, it follows from the proof of Theorem \ref{thm1.2}
that Theorem \ref{thm1.2} still holds when $D$ is assumed to be
broad and inner uniform.

\noindent $(3)$ By Lemma \ref{lem5-A-0'} in Section \ref{sec-7}, we see that $D$ has the $4b^2$-carrot property with center $x_0$.
\er

The arrangement of this paper is as follows. In Section \ref{sec-2},
we shall introduce necessary notions and notations, and record some
useful results. Section \ref{sec-3} contains three counterexamples to Heinonen's open problem and their proofs. These examples show that if one of the domains is unbounded then the answer to
Heinonen's open problem is negative. In view of these examples it suffices to consider the case where both
domains are bounded. Thus, to give a positive answer to Heinonen's open problem in this case, we need some auxiliary
results which we include in Sections \ref{sec-4} and \ref{sec-5}. Based on carrot arcs,
 in Section \ref{sec-4},  we obtain a method of construction of
uniform subdomains in a domain in $\IR^n$ and, in addition, we present four more basic
results that will be useful for our discussion in the sequel. Section \ref{sec-5} is devoted to a property for a
class of special arcs which satisfy the so-called ``$QH$-condition".
Prior to the proof of Theorem \ref{thm1.1}, in Section \ref{sec-6}, we give a way to construct
new uniform subdomains and new carrot arcs, and  obtain several other related properties.
Finally, we present a proof of the main result in this section, namely, Theorem \ref{thm6.1} from which we prove Theorem \ref{thm1.1}
easily.  In Section \ref{sec-7}, based on Theorem \ref{thm1.1} and a combination of some related results
obtained in Section \ref{sec-4}, we prove Theorem \ref{thm1.2}.

The authors thank the referee very much for their
careful reading of this paper and many valuable comments, which led to numerous improvements.
%%%%%%%%%%%%%%%%%%%%%%%%%%%%%%%%%%
%%%%%%%%%%%%%%%%%%%%%%%%%%%%%%%%%%
\section{Preliminaries}\label{sec-2}
%%%%%%%%%%%%%%%%%%%%%%%%%%%%%%%%%%
%%%%%%%%%%%%%%%%%%%%%%%%%%%%%%%%%%

%%%%%%%%%%%%%%%%%%%%%%%%%%%%%%%%%%
\subsection{Notation}
%%%%%%%%%%%%%%%%%%%%%%%%%%%%%%%%%%
Throughout the paper, we always use $\mathbb{B}(x_0,r)$ to denote
the open ball $\{x\in \IR^n:\,|x-x_0|<r\}$ centered at $x_0$ with
radius $r>0$. Similarly, for the closed balls and spheres, we use
the  notations $\overline{\mathbb{B}}(x_0,r)$ and $
\mathbb{S}(x_0,r)$, respectively. Obviously, $\IB=\mathbb{B}(0,1)$.

%%%%%%%%%%%%%%%%%%%%%%%%%%%%%%%%%%
\subsection{John domains, uniform domains and carrot
property}\label{sub2.2}
%%%%%%%%%%%%%%%%%%%%%%%%%%%%%%%%%%

\bdefe \label{def1} A domain $D$ in $\IR^n$  is said to be  {\it
$c$-uniform}  if there exists a constant $c$ with the property that
each pair of points $z_{1},z_{2}$ in $D$ can be joined by a
rectifiable arc $\gamma$ in $ D$ satisfying (cf. \cite{Martio-80})
$$\ell(\gamma)\leq c\,|z_{1}-z_{2}|
$$
and
\beq\label{sf-1} \ds\min_{j=1,2}\{\ell (\gamma [z_j, z])\}\leq c\,
d_D(z)
\eeq
for all $z\in \gamma$, where $\ell(\gamma)$ denotes the
arclength of $\gamma$, $\gamma[z_{j},z]$ the part of $\gamma$
between $z_{j}$ and $z$, and $d_D(z)$ the distance from $z$ to the
boundary $\partial D$ of $D$. Also we say that $\gamma$ is a {\it
double $c$-cone arc}. \edefe

 In order to introduce the definition of John domains, we need the following concept.

\bdefe\label{wed-1} A set $A\subset D$ in $\IR^n$  is said to have
the {\it $c$-carrot property} with center $x_0\in \overline{D}$ if
there exists a constant $c$ with the property that for each point
$z_{1}$ in $A$, $z_1$ and $x_0$ can be joined by a rectifiable arc
$\gamma$ in $ D$ satisfying (cf. \cite{RJ, Vai0})
$$\ell (\gamma [z_1, z])\leq c\, d_D(z)
$$
for all $z\in \gamma$. Also we say that $\gamma$ is a {\it
$c$-carrot arc} with center $x_0$.\edefe

Now, we are ready to introduce the carrot definition for John
domains.

 \bdefe\label{wed-1'} A domain  $D$ in $\IR^n$  is
said to be a {\it $c$-John} domain with center $x_0\in \overline{D}$
if $D$ has the $c$-carrot property with center $x_0\in
\overline{D}$. Especially, $x_0$ is taken to be $\infty$ when $D$ is
unbounded.\edefe

This is the carrot definition for John domains.
Now, there are plenty of alternative characterizations for uniform
and John domains, see \cite{Bro, FW,  Geo, Kil, Martio-80, Vai4,
Vai5, Vai6, Vai7, Vai8}.  The importance of uniform
and John domains   together with some special
domains throughout the function theory is well documented, see
\cite{FW, Kil, RJ, Vai2}. Moreover, John domains and uniform domains
in $\mathbb{R}^n$ enjoy numerous geometric and function theoretic
features in many areas of modern mathematical analysis, see
\cite{Bea, Bro, Geo, Has, Jo80, Jo81, Vai2, Vai7} (see also
\cite{Alv}).

Definitions \ref{def1} and \ref{wed-1'} (resp.
Definition \ref{wed-1}) are often referred to as the ``arclength"
definitions for uniform domains and John domains (resp. the carrot
property). When the arclength in Definitions \ref{def1} and
\ref{wed-1'} (resp. Definition \ref{wed-1}) is
replaced by the diameter, then it is called the ``diameter"
definition for uniform domains and John domains (resp. the carrot
property), and in the latter definition, the assumption on the
rectifiability of the arc is not necessary.
%
%. {\color{red} We remark that in the diameter definition of uniform and John domains, the
%hypothesis on the rectifiability of the arc is not necessary.  }

%Definitions \ref{def1}, \ref{wed-1'} and \ref{def-1} (resp.
%Definition \ref{wed-1}) are often referred to as the ``arclength"
%definitions for uniform domains and John domains (resp. the carrot
%property). When the arclength in Definitions \ref{def1},
%\ref{wed-1'} and \ref{def-1} (resp. Definition \ref{wed-1}) is
%replaced by the diameter, then it is called the ``diameter"
%definition for uniform domains and John domains (resp. the carrot
%property).

The following result reveals the close relationship between these two definitions.

\begin{Thm}\label{ThmF-1}$($\cite{MS, RJ, Vai'}$)$
The ``arclength" definition for uniform domains and John domains
$($resp. the carrot property$)$ is quantitatively equivalent to the
``diameter" one. In particular, for a John domain the center can be taken to be the same.
\end{Thm}

%%%%%%%%%%%%%%%%%%%%%%%%%%%%%%%%%%
\subsection{Quasihyperbolic metric and solid arcs}
%%%%%%%%%%%%%%%%%%%%%%%%%%%%%%%%%%

Let $\gamma$ be a rectifiable arc or path in $D$. Then the {\it
quasihyperbolic length} of $\gamma$ is defined to be the number
$\ell_{k_D}(\gamma)$ given by (cf. \cite{GP}):
$$\ell_{k_D}(\gamma)=\int_{\gamma}\frac{|dz|}{d_D(z)}.
$$
For $z_1$, $z_2$ in $D$, the {\it quasihyperbolic distance}
$k_D(z_1,z_2)$ between $z_1$ and $z_2$ is defined in the usual way:
$$k_D(z_1,z_2)=\inf\ell_{k_D}(\gamma),
$$
where the infimum is taken over all rectifiable arcs $\gamma$
joining $z_1$ to $z_2$ in $D$. An arc $\gamma$ from $z_1$ to $z_2$
is called a {\it quasihyperbolic geodesic} if
$\ell_{k_D}(\gamma)=k_D(z_1,z_2)$. Each subarc of a quasihyperbolic
geodesic is obviously a quasihyperbolic geodesic. It is known that a
quasihyperbolic geodesic between two points in $D$ always exists
(cf. \cite[Lemma 1]{Geo}). Moreover, for $z_1$, $z_2$ in $D$, we
have (cf. \cite{Vai4, Vu-book-88})
\beq
\nonumber k_{D}(z_1, z_2)&\geq &
\inf_{\gamma}\log\left (1+\frac{\ell(\gamma)}{\min\{d_D(z_1), d_D(z_2)\}}\right )
\geq
\log\left (1+\frac{|z_1-z_2|}{\min\{d_D(z_1), d_D(z_2)\}}\right )\\ \nonumber &\geq &
\Big|\log \frac{d_D(z_2)}{d_D(z_1)}\Big|,
\eeq
where $\gamma$ denote curves in $D$ connecting $z_1$ and $z_2$. If
$\gamma$ is a quasihyperbolic geodesic in $D$ joining $z_1$ to $z_2$, then
we know that for $x\in \gamma$,
\beq \label{h-w-01}
k_{D}(z_1, z_2)\geq
\log\left (1+\frac{\ell(\gamma)}{d_D(x)}\right ).
\eeq
If $|z_1-z_2|\le d_D(z_1)$, we have \cite{Vu-book-88}
\be\label{xxt} k_D(z_1,z_2)< \log\Big( 1+ \frac{
|z_1-z_2|}{d_D(z_1)-|z_1-z_2|}\Big).
\ee

The following characterization of uniform domains in terms of
quasihyperbolic metric is useful for our discussions.

\begin{Thm}\label{ThmF'} $($\cite[2.50 (2)]{Vu}$)$\quad
A domain $D\subset \IR^n$ is $c$-uniform if and only if there is a
constant $\mu_1$ such that for all $x$, $y\in D$,
$$k_{D}(x, y)\leq \mu_1\log \left(1+\frac{|x-y|}{\min\{d_{D}(x), d_{D}(y)\}}\right ),
$$
where $\mu_1=\mu_1(c)$.
\end{Thm}

This form  of the definition of uniform domains is due to
Gehring and Osgood \cite{Geo}. As a matter of fact, in \cite[Theorem
1]{Geo}, there was an additive constant in the inequality of Theorem
\Ref{ThmF'}, but it was shown by Vuorinen \cite[2.50 (2)]{Vu}
that the additive constant can be chosen to be zero.

Next, we recall a relationship between the quasihyperbolic distance of
points in a domain $D$ and the one of their images in $D'$ under
$f$.

\begin{Thm} \label{ThmF}$($\cite[Theorem 3]{Geo}$)$
Suppose that $D$ and $D'$ are domains in $\IR^n$, and that $f:\;
D\to D'$ is a $K$-quasiconformal mapping. Then for $z_1$, $z_2\in
D$,
$$k_{D'}(z'_1,z'_2)\leq \mu_2\max\{k_{D}(z_1,z_2), (k_{D}(z_1,z_2))^{\frac{1}{\mu_2}}\},
$$
where the constant $\mu_2=\mu_2(n, K)\geq 1$.
\end{Thm}

As a generalization of quasihyperbolic geodesics, V\"ais\"al\"a
introduced the concept of solids in \cite{Vai6}.

 Let $\alpha$ be an arc in $D$. The arc
may be closed, open or half open. Let $\overline{x}=(x_0,\ldots ,x_n)$,
$n\geq 1$, be a finite sequence of successive points of $\alpha$.
For $h\geq 0$, we say that $\overline{x}$ is {\it $h$-coarse quasihyperbolic} if
$k_{D}(x_{j-1}, x_j)\geq h$ for all $1\leq j\leq n$. Let $\Phi_{k_D}(\alpha,h)$
be the family of all $h$-coarse quasihyperbolic sequences of $\alpha$. Set
$$s_{k_D}(\overline{x})=\sum^{n}_{j=1}k_{D}(x_{j-1}, x_j)
$$
and
$$\ell_{k_D}(\alpha, h)=\sup \{s_{k_D}(\overline{x}):\, \overline{x}\in \Phi_{k_D}(\alpha,h)\}
$$
with the agreement that $\ell_{k_D}(\alpha, h)=0$ if
$\Phi_{k_D}(\alpha,h)=\emptyset$. Then the number $\ell_{k_D}(\alpha, h)$ is the
{\it $h$-coarse quasihyperbolic length} of $\alpha$.

\bdefe \label{def1.5} Let $D$ be a domain in $\mathbb{R}^n$. An arc
$\alpha\subset D$ is {\it $(\nu, h)$-solid} with $\nu\geq 1$ and
$h\geq 0$ if
$$\ell_{k_D}(\alpha[x,y], h)\leq \nu\;k_D(x,y)
$$
for all $x, y\in \alpha$.\edefe

%{\color{red} A {\it $(\nu,0)$-solid arc} is said to be a {\it
%$\nu$-neargeodesic}, i.e. an arc $\alpha\subset D$ is a
%$\nu$-neargeodesic if and only if $\ell_k(\alpha[x,y])\leq
%\nu\;k_D(x,y)$ for all $x, y\in \alpha$.}\edefe

%{\color{red} Obviously, a $\nu$-neargeodesic is a quasihyperbolic geodesic if and
%only if $\nu=1$.

%In \cite{Vai4}, V\"ais\"al\"a got the following property concerning
%the existence of neargeodesics in $E$.

%\begin{Thm}\label{LemA} $($\cite[Theorem 3.3]{Vai4}$)$
%Let $\{z_1,\, z_2\}\subset D$ and $\nu>1$. Then there is a
%$\nu$-neargeodesic in $D$ joining $z_1$ and $z_2$.
%\end{Thm}
%}

Obviously, each quasihyperbolic geodesic is $(1,0)$-solid. The
following result is due to V\"ais\"al\"a.

\begin{Thm}\label{Thm4-1} {\rm (\cite[Theorem 6.22]{Vai6})} Suppose that
$\gamma\subset D\not= \IR^n$ is a $(\nu,h)$-solid arc with endpoints
$a_0$, $a_1$ and that $D$ is a $c$-uniform domain. Then there is a
constant $c_2$ such that

\begin{enumerate}
\item  $\ds\min\Big\{\diam(\gamma [a_0, z]),\diam(\gamma [a_1, z])\Big\}\leq c_2\,d_{D}(z)
$ for all $z\in \gamma$, and

\item  $\diam(\gamma)\leq c_2\max\Big\{|a_0-a_1|,2(e^h-1)\min\{d_{D}(a_0),d_{D}(a_1)\}\Big\}$,
\end{enumerate}
where $c_2=c_2(\nu,h,c)\geq 1$.
\end{Thm}

\subsection{Internal metric, inner uniform domains and broad domains}\label{subsec2.5}
%%%%%%%%%%%%%%%%%%%%%%%%%%%%%%%%%%

For $x$, $y$ in $D$, the internal metric $\delta_{D}$ in $D$ is
defined by
$$\delta_D(x,y)=\inf \{\diam(\alpha):\; \alpha\subset D\;
\mbox{is a rectifiable arc joining}\; x\; \mbox{and}\; y \}.
$$

\bdefe \label{def1-0} A domain $D$ in $\IR^n$  is said to be  {\it
$c$-inner uniform}  if there exists a constant $c$ with the property
that each pair of points $z_{1},z_{2}$ in $D$ can be joined by a
rectifiable arc $\gamma$ in $ D$ satisfying \eqref{sf-1} and
$\ell(\gamma)\leq c\,\delta_D(z_1,z_2)$ (cf. \cite{Vai7}). \edefe

Obviously, ``uniformity" implies ``inner uniformity".

\bdefe \label{def3}
Let $\varphi: (0, \infty)\to  (0, \infty)$ be a
decreasing homeomorphism. We say that $D$ is {\it $\varphi$-broad}
if for each $t>0$ and each pair $(C_0, C_1)$ of continua in $D$ the
condition $\delta_{D}(C_0, C_1)\leq t\min\{\diam(C_0),\diam(C_1)\}$
implies
$${\rm Mod}\,(C_0, C_1; D)\geq \varphi(t),
$$
where $\delta_{D}(C_0, C_1)$ denotes the $\delta_{D}$-distance
between $C_0$ and $C_1$ and see, for example, \cite{Vai} for the definition of ${\rm Mod}\,(C_0, C_1; D)$, the modulus of a family
of the curves in $D$ connecting $C_0$ and $C_1$.
\edefe

Broad domains were introduced in \cite{Vai0} and it was later proved
that a simply connected planar domain is broad if and only if it is
inner uniform \cite{RJ, Vai7}. Further, the following is known.

\begin{Thm} \label{Lem4-1}$($\cite[Lemma 2.6]{Gem}$)$
If $D$ is a $c$-uniform domain, then $D$  is  $\varphi_c$-broad,
where $\varphi_c$ depends only on $c$.
\end{Thm}

%{\color{red}
It is important to recall that the notion of broad domains also goes under the
term L\"{o}ewner space. %especially in geometric group theory.
The notion of a Loewner space was introduced by Heinonen and Koskela \cite{HeK} in their
study of quasiconformal mappings of metric spaces; Heinonen's recent monograph \cite{Hei} renders 
an enlightening account of these ideas. See \cite{BaK, BoHK, Her, Ty} etc for more related discussions. 
%}

%\begin{Thm} \label{Thm-thu}$($\cite[??]{??}$)$
%If $D$ is a $c$-inner uniform domain, then $D$ is $??$-broad, where
%$??$ depends only on $c$.
%\end{Thm}

%%%%%%%%%%%%%%%%%%%%%%%%%%%%%%%%%%
\subsection{Quasisymmetric
mappings}
%%%%%%%%%%%%%%%%%%%%%%%%%%%%%%%%%%

\bdefe \label{def2'} Let $X_1$ and $X_2$ be two metric spaces with
distance written as $|x-y|$, and let $\eta: [0, \infty)\to [0,
\infty)$ be a homeomorphism. We say that an embedding $f: X_1\to X_2$ is {\it
$\eta$-quasisymmetric} if $|a-x|\leq t|a-y|$ implies
$$|f(a)-f(x)|\leq \eta(t)|f(a)-f(y)|
$$
for all $a$, $x$, $y\in X_1$, and if  there is a constant $\tau\geq
1$ such that $|a-x|\leq |a-y|$ implies
$$|f(a)-f(x)|\leq \tau|f(a)-f(y)|,
$$
then $f$ is said to be {\it weakly $\tau$-quasisymmetric}. \edefe
Obviously, ``quasisymmetry" implies ``weak quasisymmetry". 
%{\color{red} 
There are many references in literature in this line, see, for example, \cite{BaK, BoK, BoM, Vai0-1} etc. Among them, V\"ais\"al\"a proved the following
local quasisymmetry of quasiconformal mappings in $\IR^n$. 
%}

\begin{Thm} $($\cite[Theorem 2.4]{Vai0-1}$)$ \label{Vai-0}
Suppose that $D\subset\IR^n$, and that $f$ is a $K$-quasiconformal mapping of $D$ onto a
domain $D'\subset \IR^n$. Then for $x\in D$ and $0<\lambda<1$, the
restriction $f|_{\mathbb{B}(x,\lambda\,d_D(x))}$ is
$\eta$-quasisymmetric, where $\eta=\eta_{n, K, \lambda}$.
\end{Thm}

%%%%%%%%%%%%%%%%%%%%%%%%%%%%%%%%%%
\subsection{Linearly locally connected sets}
%%%%%%%%%%%%%%%%%%%%%%%%%%%%%%%%%%

\bdefe \label{def2} Suppose that $A\subset D$ and $b\geq 1$ is a
constant.  We say that $A$ is $b$-$LLC_2$ (resp. $b$-$LLC_2$ with
respect to $\delta_{D}$) in $D$ if for all $x\in A$ and $r>0$, any two
points in $A \backslash \overline{\mathbb{B}}(x, br)$ $($resp. $A
\backslash \overline{\mathbb{B}}_{\delta_{D}}(x, br))$ can be joined
in $D\backslash \overline{\mathbb{B}}(x, r)$ $($resp. $D\backslash
\overline{\mathbb{B}}_{\delta_{D}}(x, r))$, where
$$\mathbb{B}_{\delta_{D}}(x, r)=\{z\in \IR^n:\; \delta_{D}(z, x)< r\}.
$$
If $A=D$, then we say that $D$ is $b$-$LLC_2$ $($resp. $b$-$LLC_2$ with respect to
$\delta_{D})$.
\edefe

In \cite{H}, Heinonen proved

\begin{Thm}\label{Lem6-1}$($\cite[Lemma 7.2]{H}$)$
Let $D\subset \IR^n$, and let $A$ be a $\varphi$-broad subdomain of $D$. Then $A$ is $b$-$LLC_2$ with
$b=b( n, \varphi)$. In particular, $A$ is
 $\mu_3$-$LLC_2$ with respect to $\delta_{D}$ in $D$, where
$\mu_3=\mu_3(n,\varphi)$.
\end{Thm}

For convenience, in what follows in this paper, we always assume
that $x, y, z, \ldots$ denote points in a domain $D$ in $\IR^n$ and
$x', y', z', \ldots$ the images in $D'$ of $x, y, z, \ldots$ under
$f$, respectively. Also we assume that $\alpha, \beta, \gamma,
\ldots$ denote curves in $D$ and $\alpha', \beta', \gamma', \ldots$
the images in $D'$ of $\alpha, \beta, \gamma, \ldots$ under $f$,
respectively, and for a set $A$ in $D$, $A'$ always denotes its image in $D'$
under $f$.

%%%%%%%%%%%%%%%%%%%%%%%%%%%%%%%%%%
%%%%%%%%%%%%%%%%%%%%%%%%%%%%%%%%%%
\section{Counterexamples to Heinonen's open problem} \label{sec-3}
%%%%%%%%%%%%%%%%%%%%%%%%%%%%%%%%%%
%%%%%%%%%%%%%%%%%%%%%%%%%%%%%%%%%%

In this section, we shall construct three examples to show that the answer to  Open problem \ref{Con2}
is negative when one of $D$ and $D'$ is unbounded, which are as follows.

\begin{exam}\label{exa0}
Suppose that $\mathbb{H^+}$ denotes the upper half plane of
$\IR^2\cong \mathbb{C}$ and $\mathbb{D}=\{z:\, |z|<1\}$, the unit disk in $\mathbb{C}$, and
suppose that $f_1(z)=i\frac{z-i}{z+i}$, where $i^2=-1$. Then $f_1$ maps $\mathbb{H^+}$
onto $\mathbb{D}$ and does not have the
subinvariance property of John domains.
\end{exam}

\begin{exam}\label{exa1}
Suppose that $f_2(z)=i\frac{1+z}{1-z}$. Obviously, $f_2$ maps
$\mathbb{D}$ onto $\mathbb{H^+}$. Then  $f_2$ does not have the
subinvariance property of John domains.
\end{exam}

\begin{exam}\label{exa2}
Suppose that $f_3(z)=-\frac{1}{z}$. Obviously, $f_3$ maps
$\mathbb{H^+}$ onto $\mathbb{H^+}$. Then  $f_3$ does not have the
subinvariance property of John domains.
\end{exam}

\begin{figure}[!ht]
\begin{center}
\includegraphics%[width=10cm]
{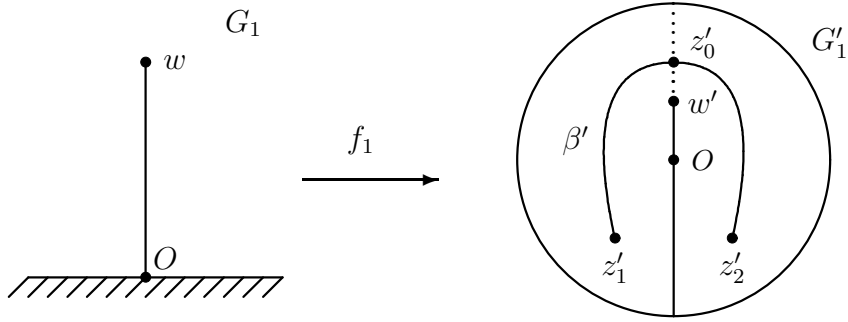}
\end{center}
\caption{ The domains $G_1$, $G'_1$, the points $z'_0$, $z'_1$,
 $z'_2$ and the arc $\beta'$ $(w=(8c-1)i$ and $w'=(1-\frac{1}{4c})i)$\label{fig3}}
\end{figure}

\noindent {\bf The proof of Example \ref{exa0}.} \quad
It is known that both $\mathbb{H}^+$ and $\mathbb{D}$ are $\frac{\pi}{2}$-John disks.
Suppose on the contrary that there is a constant $c\geq 1$ such that for every $4$-John
subdomain in $\mathbb{H}^+$, its image under the conformal mapping $f_1$ is
$c$-John in $\ID$, where $f_1(z)=i\frac{z-i}{z+i}$ and $i^2=-1$. We let
$$G_1=\mathbb{H}^+\backslash (0, (8c-1)i].
$$
Elementary calculations show that $G_1$ is a $4$-John domain and
$$G'_1=\mathbb{D}\backslash (-i, (1-\frac{1}{4c})i].
$$
We take
$$z'_1=-\frac{1}{2}-\frac{3}{4}i\;\;\mbox{and}\;\; z'_2=\frac{1}{2}-\frac{3}{4}i.
$$
Then there must exist a $c$-cone arc $\beta'$ in $G'_1$ joining $z'_1$ and $z'_2$. The fact that
$G'_1$ is a bounded $c$-John domain guarantees the existence of such
an arc $\beta'$. Let $z'_0$ be one of the intersection points of
$\beta'$ with the interval $[(1-\frac{1}{4c})i, i)$ in the imaginary
axis (see Figure \ref{fig3}). Obviously, $d_{G'_1}(z'_0)\leq
\frac{1}{8c},$ and thus
$$1<\min\{|z'_1-z'_0|, |z'_2-z'_0|\}\leq \min\{\ell(\beta'[z'_1,
z'_0]), \ell(\beta'[z'_0, z'_2])\}\leq cd_{D'_1}(z'_0)\leq
\frac{1}{8}.
$$
This obvious contradiction completes the proof of this example. \qed

%{\color{red} Figure 4: The domains $G_1$, $G'_1$, the points $z'_1$,
%$z'_2$ $z'_3$, and the arc $\beta'$}

% \begin{figure}[!ht]
%\includegraphics[width=0.75\textwidth]{figure04x} %,height=0.4\textwidth
%\caption{ The domains $G_1$, $G'_1$, the points $z'_0$, $z'_1$,
% $z'_2$ and the arc $\beta'$ $(w=(8c-1)i$ and $w'=(1-\frac{1}{4c})i)$\label{fig3}}
%\end{figure}

\begin{figure}[!ht]
\begin{center}
\includegraphics%[width=10cm]
{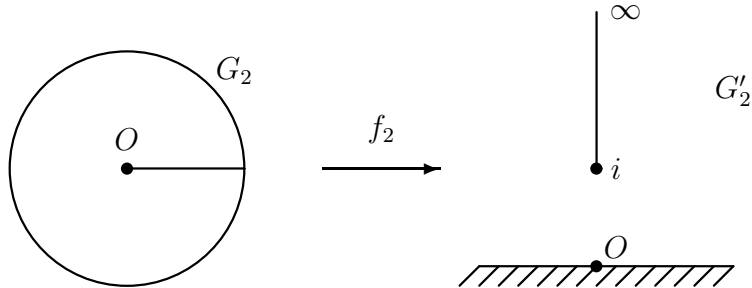}
\end{center}
\caption{The domains $G_2$ and $G'_2$\label{fig4}}
\end{figure}

\noindent {\bf The proof of Example \ref{exa1}.}  \quad Let $G_2=\mathbb{D}\backslash [0,
1)$. Since $f_2(z)=i\frac{1+z}{1-z}$, we see that the image of $G_2$
under $f_2$ is $G'_2=\mathbb{H}^+\backslash [i, +\infty)$ (see
Figure \ref{fig4}). It easily follows from \cite[Section 2.4]{RJ}
that $G_2$ is a John disk, but obviously, $G'_2$ is not a John domain. \qed

%{\color{red} Figure 5: The domains $G_2$ and $G'_2$}

% \begin{figure}[!ht]
%\includegraphics[width=0.75\textwidth]{figure05x} %,height=0.4\textwidth
%\caption{ The domains $G_2$ and $G'_2$\label{fig4}}
%\end{figure}

\begin{figure}[!ht]
\begin{center}
\includegraphics%[width=10cm]
{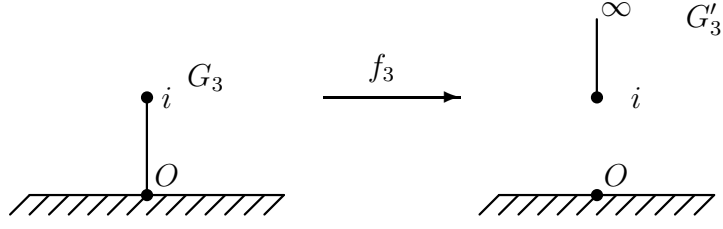}
\end{center}
\caption{ The domains $G_3$ and $G'_3$\label{fig3-1}}
\end{figure}

\noindent {\bf The proof of Example \ref{exa2}.}\quad Let $G_3=\mathbb{H}^+\backslash (0,
i]$. Then under the mapping $f_3(z)=-\frac{1}{z}$, the image $G'_3=\mathbb{H}^+\backslash [i, +\infty)$ (see
Figure \ref{fig3-1}). Elementary calculations show that $G_3$ is a $4$-John domain, but obviously, $G'_3$ is not a John domain. \qed

%%%%%%%%%%%%%%%%%%%%%%%%%%%%%%%%%%
%%%%%%%%%%%%%%%%%%%%%%%%%%%%%%%%%%
\section{The properties of carrot arcs}\label{sec-4}
%%%%%%%%%%%%%%%%%%%%%%%%%%%%%%%%%%
%%%%%%%%%%%%%%%%%%%%%%%%%%%%%%%%%%

 In this section, we shall prove some auxiliary results which are useful for the proofs of Theorems \ref{thm1.1} and \ref{thm1.2} given in the last two sections. Our first results concern the construction of uniform subdomains in a domain based on carrot arcs.

Suppose $G$ is a domain in $\IR^n$. For $z_1$ and $z_0$ in $G$, if
$\gamma$ is a $c$-carrot arc joining $z_1$ and $z_0$ in $G$ with
center $z_0$, then we have the following lemma.
%Antti's figure
\begin{figure}[!ht]
\includegraphics[width=10cm]{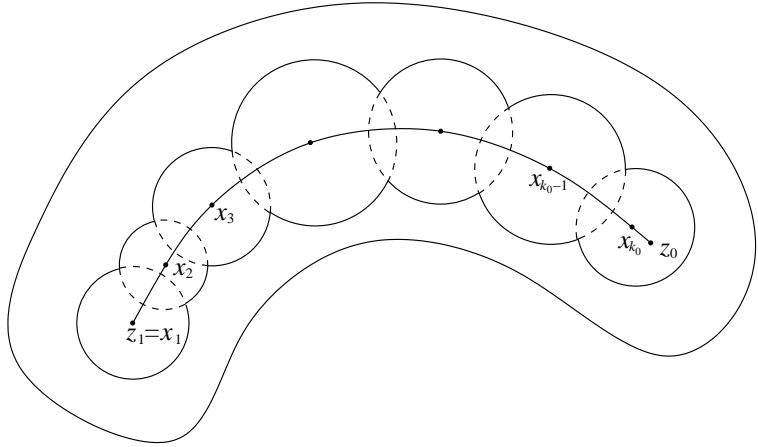}
\caption{The domain $G_{0,0}$\label{fig1}}
\end{figure}

\begin{lem}\label{lem-m-v-1}
There exists a simply connected domain $G_{1,0}=\bigcup
\limits_{i=1}^{k_0}B_i\subset G$ such that
\bee
\item\label{equ--1}
$z_{1}$, $z_{0}\in G_{1,0}$;
\item\label{equ--2}
for each $i\in \{1,\ldots, k_0\}$,
$$\frac{1}{3\mu_4}\,d_{G}(x_i)\leq r_i\leq \frac{1}{\mu_4}d_{G}(x_i);
$$
\item\label{equ--3}
if $k_0\geq 3$, then for all $i,j\in\{1,\ldots,k_0\}$ with
$|i-j|\geq 2$,
$$\dist(B_i, B_j)\geq \frac{1}{2^5\mu_5}\max\{r_i,r_j\};
$$
\item\label{equ--4}
if $k_0\geq 2$, then for each $i\in\{1,\ldots,k_0-1\}$
$$\ds r_i+r_{i+1}-|x_i-x_{i+1}|\ge \frac{1}{2^5\mu_5}\max\{r_i,r_{i+1}\};
$$
\item\label{equ--4-11}
$d_G(z_0)\leq 2^8\mu_4\mu_5d_{G_{1,0}}(z_0),$
\eee
where $B_i=\mathbb{B}(x_i, r_i)$, $x_i\in \gamma$, $x_i\not\in B_{i-1}$
for each $i\in \{2, \ldots, k_0\}$, $x_1=z_1$, $\mu_4\geq 2$ is a
constant and $\mu_5=2^{32+c^2\mu_4^2}$ $($see Figure \ref{fig1}$)$.
\end{lem}
\bpf
The proof is similar to that of \cite[Lemma 2.2]{Ml}: By replacing the radius
$\frac{1}{2}d_{D}(\cdot)$ in the proof of \cite[Lemma 2.2]{Ml} by
the one $\frac{1}{2\mu_4}d_{G}(\cdot)$, we shall get a finite sequence of balls in $D$; and then we multiply the radius of the last ball obtained by $1+\frac{1}{2^5\mu_5}$ times. The new finite sequence of balls is needed. \epf

%\begin{figure}
%
%\includegraphics[width=8cm]{figure01}
%
%\caption{the caption text}\label{fig1}
%\end{figure}

\br It easily follows from the construction in Lemma \ref{lem-m-v-1}
that the intersection $\partial G_{1,0}\cap
\partial G$ is an empty set.
\er

%Next, we  prove the uniformity of $G_{1,0}$. In order to prove the next lemma we need the following result which
%is from \cite{Yli}.%^

%\begin{Thm}\label{Lem5.2}{\rm (\cite[Theorem 1.2]{Yli})}
%Suppose that $E_1$ and $E_2$ are two convex domains in a Banach
%space, where $E_1$ is bounded and $E_2$ is $c$-uniform for some
%$c>1$, and suppose that there exist $z_0\in E_1\cap E_2$ and $r>0$
%such that $\IB(z_0,r) \subset E_1\cap E_2$. If there exist constants
%$R_1>0$ and $c_0>1$ such that $R_1 \leq c_0r$ and $E_1\subset
%\IB(z_0,R_1)$, then $E_1\cup E_2$ is a $\rho$-uniform domain with
%$\rho=(c+1)(2c_0+1)+c$.
%\end{Thm}

\begin{lem}\label{lem2.2-4-mv-3}
The domain $G_{1,0}$ constructed in Lemma {\rm \ref{lem-m-v-1}} is a
$2^{12}c^2\mu_4\mu_5$-uniform domain.
\end{lem} %\bpf
%By the definition of $G_{1,0}$, we only need to show that for $u_1$ and $u_2\in
%G_{1,0}$, there is a double $2^{12}c^2\mu_4\mu_5$-cone arc in
%$G_{1,0}$ connecting $u_1$ and $u_2$. The construction of such an
%arc is similar to that in the proof of
% \cite[Theorem $1.8$]{HW}. \epf
\bpf
By definition, we only need to show that for $u_1$ and $u_2\in
G_{1,0}$, there is a double $2^{12}c^2\mu_4\mu_5$-cone arc $\lambda$ in $G_{1,0}$
connecting $u_1$ and $u_2$.

We first consider the case where there is an $i\in
\{1,\ldots,k_0-1\}$ such that $u_1$, $u_2\in B_i\cup B_{i+1}$. Under
this assumption, clearly, the existence of $\lambda$ easily follows from
\cite[Theorem 1.2]{Yli}.

Next, we consider the case where there are $i, j\in \{1,\ldots,k_0\}$
such that $j-i\geq 2$, $u_1\in B_i$, $u_2\in B_j$ and $\{u_1, u_2\}$
is not contained in $B_t\cup B_{t+1}$ for all $t\in \{i, \ldots,
j-1\}$. It suffices to discuss the case: $u_1\notin [x_i, x_{i+1}]$
and $u_2\notin [x_{j-1}, x_j]$ since the discussions for other cases
are similar, where $[x_i, x_{i+1}]$ denotes the line segment with endpoints $x_i$ and
$x_{i+1}$.
In this case, we let
$$\lambda=[u_1,x_i]\cup[x_i, x_{i+1}]\cup\cdots\cup[x_{j-1}, x_j]\cup[x_j,u_2].
$$
For each $u\in \lambda$, if $u\in [u_1, x_i]$, then
$$\ell(\lambda[u_1, u])=|u_1-u|<d_{G_{1,0}}(u).
$$
Also if $u\in [u_2, x_j]$, then
$$\ell(\lambda[u_2, u])=|u_2-u|<d_{G_{1,0}}(u).
$$

In the following, we consider the case $u\in \lambda[x_i, x_j]$.
Clearly, there exists a $t\in \{i,\ldots, j\}$ such that $u\in B_t$.
If $u\in B_j$, then  by  Lemma \ref{lem-m-v-1},
$$\ell(\lambda[u_2, u])\leq 2r_j\leq 2^7\mu_5 d_{G_{1,0}}(u).
$$
Similarly, if $u\in B_i$, then
$$\ell(\lambda[u_1, u])\leq 2r_i\leq 2^7\mu_5 d_{G_{1,0}}(u).
$$
So it remains to consider the case: $u \in \lambda \backslash
\big(B_i\cup B_j\big)$. Then Lemma \ref{lem-m-v-1} implies
\begin{eqnarray*}
\ell(\lambda[u_1, u])&\leq& \ell(\lambda[u_1, x_t])+r_t\leq 3c\mu_4\;d_{G_{1,0}}(x_t)+r_t \leq 2^6(3c\mu_4+1)\mu_5\;d_{G_{1,0}}(u).
\end{eqnarray*}
Hence $\lambda$ is a $2^6(3c\mu_4+1)\mu_5$-cone arc. It follows
from Lemma \ref{lem-m-v-1} that
%\be\label{eq1.6.0-7}
$$\max\{r_i,r_j\}\leq 2^5\mu_5|u_1-u_2|,
$$
whence
\begin{eqnarray*}
\ell(\lambda[u_1, u_2])&\leq& r_i+\ell(\lambda[x_i, x_j])+r_j \leq 2^6\mu_5|u_1-u_2|+ cd_{G}(x_j) \leq 2^7c\mu_4\mu_5|u_1-u_2|.
\end{eqnarray*}
The proof of the existence of $\lambda$ is finished. \epf

Our next result is an estimate on the distance from the points in carrot arcs to the boundary of $G$, which plays a key role in the discussions in Section \ref{sec-5}. For the proof of the estimate, the following result is necessary.

\begin{lem}\label{lem-4-1''}
For $x,y\in G$, if $|x-y|\leq \mu_6d_G(y)$, then $d_G(y)\geq \frac{1}{2\mu_6}d_G(x)$, where
$\mu_6\geq 1$ is a constant.
%\end{enumerate}
\end{lem}

\bpf We divide the proof into two cases: $|x-y|<\frac{1}{2}d_G(x)$
and $|x-y|\geq \frac{1}{2}d_G(x)$. For the first case, we see that
$$d_G(y)\geq d_G(x)-|x-y|\> >\frac{1}{2}d_G(x).
$$
For the second case,
$$d_G(y)\geq \frac{1}{\mu_6}|x-y|\geq \frac{1}{2\mu_6}d_G(x).
$$
The proof is complete.\epf

Suppose $\chi\subset G$ is a  $\mu_7$-carrot arc from $x_1$ to $x_2$
with the center $x_2$ in $G$, where $\mu_7 \geq 1$ is a constant.
Then we have

\begin{lem}\label{lem-4-1}
For every $u\in \chi[x_1, x_2]$,
$$\ds d_G(u)\geq \frac{2\ell(\chi[x_1,u])+d_G(x_1)}{4\mu_7}.
$$
\end{lem}
\bpf For every $u\in \chi[x_1, x_2]$,
$$d_G(u)\geq \frac{\ell(\chi[x_1,u])}{\mu_7}.
$$
Then by Lemma \ref{lem-4-1''}, we have
$$d_G(u)\geq\frac{1}{2\mu_7}d_G(x_1).
$$
Hence
$$d_G(u)\geq \frac{2\ell(\chi[x_1,u])+d_G(x_1)}{4\mu_7}
$$
as required. \epf

Our last two results are on the estimate on the quasihyperbolic distance between two points in a ball in $G$ and the comparison of the distances from the images of the points in a ball in $G$ to the boundary of $G'$, which will be crucial for the discussions in Sections \ref{sec-5} and \ref{sec-6}.

\begin{lem}\label{mxll-1}
Suppose $G$ is a domain in $\mathbb{R}^n$ and $x\in G$.
Then for all $y\in \overline{\mathbb{B}}(x, \frac{1}{b}d_G(x))$ with $b>1$, we have
$$k_G(x, y)< \frac{1}{b-1}.
$$
\end{lem}\bpf
Since $d_G(z)\geq d_G(x)-|x-z|\geq (1-\frac{1}{b})d_G(x)$ for each $z\in[x, y]$, we get
$$k_{G}(x, y)\leq\int_{[x,y]}\frac{|dz|}{d_G(z)} \leq
\frac{|x-y|}{(1-\frac{1}{b})d_G(x)} < \frac{1}{b-1}.
$$
The proof is complete.
\epf

\begin{lem}\label{mxll-1''}
Suppose that $G$, $G'\subset\mathbb{R}^n$, that $f:G\to G'$ is a $K$-quasiconformal
mapping, and that $x\in G$. Then for each $y\in
\overline{\mathbb{B}}(x, \frac{24}{\rho_1}d_G(x))$, we have
$$|x'-y'|<\frac{1}{5}\min\{d_{G'}(x'),d_{G'}(y')\}\;\;\mbox{and}\;\;\frac{4}{5}d_{G'}(y')<
d_{G'}(x')< \frac{5}{4}d_{G'}(y'),
$$
where $\rho_1=(6\mu_2)^{4\mu_2}$ and
$\mu_2$ is from Theorem {\rm \Ref{ThmF}}.
\end{lem}\bpf
For $x\in G$ and $y\in \overline{\mathbb{B}}(x,
\frac{24}{\rho_1}d_G(x))$, by Lemma  \ref{mxll-1}, we have
$$k_G(x,y)\leq \frac{24}{\rho_1-24}<1,
$$
and so Theorem \Ref{ThmF} implies
\beqq
\log\left (1+\frac{|x'-y'|}{\min\{d_{G'}(x'),d_{G'}(y')\}}\right )
&\leq & k_{G'}(x', y')\\
&\leq & \mu_2\big(k_G(x, y)\big)^{\frac{1}{\mu_2}}\\
&<&\frac{1}{6}.
\eeqq
It follows that
$$|y'-x'|< \frac{1}{5}\min\{d_{G'}(x'),d_{G'}(y')\},
$$
whence
$$ d_{G'}(y')\geq d_{G'}(x')-|x'-y'|> \frac{4}{5}d_{G'}(x')
$$
and
$$ d_{G'}(x')\geq d_{G'}(y')-|x'-y'|> \frac{4}{5}d_{G'}(y'),
$$
from which  the proof follows.\epf

%%%%%%%%%%%%%%%%%%%%%%%%%%%%%%%%%%
%%%%%%%%%%%%%%%%%%%%%%%%%%%%%%%%%%
\section{ $QH$ condition }\label{sec-5}
%%%%%%%%%%%%%%%%%%%%%%%%%%%%%%%%%%
%%%%%%%%%%%%%%%%%%%%%%%%%%%%%%%%%%

In this section, we introduce a class of special arcs and prove a related
lemma. This lemma is useful for the proof of
Theorem \ref{thm6.1} given in the next section.

Suppose that $G$ is a domain in $\mathbb{R}^n$. For $w_1$, $w_2\in G$,
we suppose $\beta\subset G$ is an arc with endpoints $w_1$ and
$w_2$. Then we introduce the following concept.
\medskip

\begin{figure}[!ht]
\begin{center}
\includegraphics%[width=10cm]
{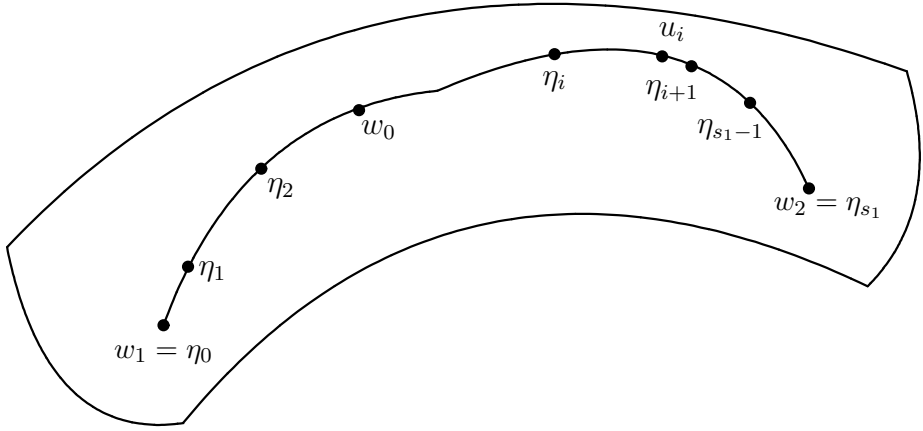}
\end{center}
\caption{The arc $\beta$ \label{fig5}}
\end{figure}

\noindent {\bf $QH$ condition.}\quad  We say that  $\beta$
satisfies {\it $QH\,(s_0,s_1, s_2)$-condition} in $G$ if there
exist positive constants $s_0$, $s_1$ and $s_2$, where $s_1$ is an
integer and $s_2\geq 1$, satisfying the following (see Figure \ref{fig5}):
\bee
\item[$(1)$] $k_{G}(w_1,w_2)\geq s_0$;

\item[$(2)$] There exist successive points $\eta_0 ~(=w_1),
\eta_1, \ldots, \eta_{s_1}~(=w_2)$ in $\beta$ such that for all $i,
j\in \{0,\ldots,s_1-1\}$,
$$\frac{1}{s_2}k_{G}(\eta_j,\eta_{j+1})\leq k_{G}(\eta_i,\eta_{i+1})\leq s_2
k_{G}(\eta_j,\eta_{j+1});
$$
\item[$(3)$] For each $i\in\{0,\ldots,s_1-2\}$ and $\eta\in
\beta[\eta_{i+1},w_2]$,
$$ k_{G}(\eta_i,\eta_{i+1})\leq s_2k_{G}(\eta_i,\eta).
$$
\eee
%\medskip

%{\color{red} The place for Figure 8:.}
%\begin{figure}[!ht]
%\includegraphics[width=0.85\textwidth]{figure08} %,height=0.4\textwidth
%\caption{The arc $\beta$\label{fig5}  {\color{red} Add $u_i$ between $\eta_i$ and $\eta_{i+1}$} }
%\end{figure}

Roughly speaking, the $QH$ condition says that $\beta$ can be partitioned into a number of pieces, such that
the $k_G$-distance between the endpoints of the pieces are comparable, and so that the $k_G$-distance between
$\eta_i$ and $\eta_{i+1}$ of the $i$th piece is comparable to the $k_G$-distance from the left endpoint $\eta_i$
to any point in the subcurve of $\beta$ from the right endpoint $\eta_{i+1}$ to the final endpoint $w_2$.
By the definition of $QH$ condition, the following is obvious.

\begin{cor}\label{scl-0-1}
Suppose $\beta$ is a quasihyperbolic geodesic in $G$ and
$0<s_0\leq \ell_{k_{G}}(\beta)<\infty$. Then it satisfies $QH\,(s_0,s_1,
s_2)$-condition in $G$, where $s_1=2$ and $s_2=1$.
\end{cor}

Let $w_0\in \beta$ be such that (see Figure \ref{fig5})
\beq \label{q-1"}
d_{G}(w_0)=\inf_{z\in \beta}d_{G}(z).
\eeq
For each $i\in \{0,\ldots,s_1-1\}$, we let $u_i\in \beta[\eta_i,\eta_{i+1}]$ (see Figure \ref{fig5})
satisfy
$$k_{G}(u_i,\eta_{i+1})=\frac{1}{4s_2} k_{G}(\eta_i,\eta_{i+1}).
$$

Suppose that $G$ is an $a$-John domain with center $x_0$, and that
$s_1\geq [(4a)^{2n}]$ is an integer, $s_2\geq 1$, $s_0\geq
8s_1s_2^2$. Then we have the following result.

\begin{thm}\label{scl-1-1}
Suppose $\beta$ denotes an arc in $G$ which satisfies
$QH\,(s_0,s_1, s_2)$-condition in $G$. Then
\bee
\item[$(1)$] for every $i\in \{0, 1, \ldots, s_1-1\}$,
$$\min\{d_{G}(u_i), d_{G}(\eta_{i+1})\}< \diam(\beta);
$$

\item[$(2)$] $k_{G}(w_1, w_2)\leq 48a^2s_1s_2^2\log \left
(1+\frac{\diam(\beta)}{d_{G}(w_0)}\right ). $
\eee
\end{thm}
\bpf We prove the conclusion $(1)$ in Theorem \ref{scl-1-1} by a method of
contradiction. Suppose on the contrary that there exists some $t\in
\{0, 1, \ldots, s_1-1\}$ such that
%\beq\label{h-l-0-1}
$$\min\{d_{G}(u_t), d_{G}(\eta_{t+1})\}\geq \diam(\beta).
$$
%\eeq
Then
$$[u_t, \eta_{t+1}]\subset \overline{\mathbb{B}}\big (u_t, \frac{1}{2}d_{G}(u_t)\big )
\cup \overline{\mathbb{B}} \big (\eta_{t+1},
\frac{1}{2}d_{G}(\eta_{t+1})\big ),
$$
which implies
$$k_{G}(u_t, \eta_{t+1})< 2.
$$
Then by the definition of $QH$ condition, we see that for each $i\in
\{0, 1, \ldots, s_1-1\}$,
$$k_{G}(\eta_i,\eta_{i+1})\leq s_2k_{G}(\eta_t,\eta_{t+1})= 4s_2^2k_{G}(u_t, \eta_{t+1}),
$$
and so
$$k_{G}(w_1,w_2)\leq \sum_{i=0}^{s_1-1}k_{G}(\eta_i,\eta_{i+1})
\leq 4s_2^2\sum_{i=0}^{s_1-1}k_{G}(u_t,\eta_{t+1}) < 8s_1s_2^2,
$$
which contradicts the assumption that ``$k_{G}(w_1,w_2)\geq s_0$", and
thus  the proof of $(1)$ in Theorem \ref{scl-1-1} is complete.
\smallskip

\begin{figure}[!ht]
\begin{center}
\includegraphics%[width=10cm]
{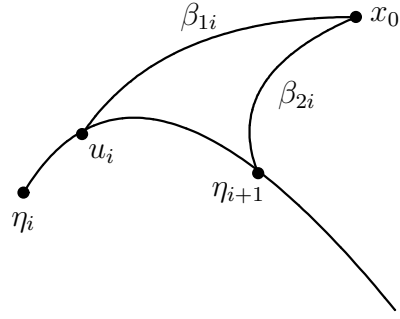}
\end{center}
\caption{The arc $\beta_i=\beta_{1i}\cup\beta_{2i}$ \label{fig6}}
\end{figure}

To prove $(2)$ of this lemma, we let $\beta_{1i}$ be an $a$-carrot
arc joining $u_i$ and $x_0$ in $G$, and let $\beta_{2i}$ be an
$a$-carrot arc joining $\eta_{i+1}$ and $x_0$ in $G$ for each $i\in
\{0,\ldots,s_1-1\}$ (see Figure \ref{fig6}). This can be done
because $G$ is an $a$-John domain with center $x_0$. Let
$$\beta_i=\beta_{1i}\cup\beta_{2i}.$$

%{\color{red} The place for Figure 9:.}

%\begin{figure}[!ht]
%\includegraphics[width=0.35\textwidth]{figure09} %,height=0.4\textwidth
%\caption{The arc $\beta_i$\label{fig6}  {\color{red} Add $x_0$ and change $\beta_{i}$ to $\beta_{1i}$ and $\beta_{2i}$} }
%\end{figure}

If there exists an $s\in \{0,\ldots,s_1-1\}$ such
that
$$ \ell(\beta_s)\leq \diam(\beta),
$$
then we see that $x_0\in G$. Thus from Lemma \ref{lem-4-1}, we deduce that
\begin{eqnarray*}
\nonumber k_{G}(u_s,\eta_{s+1})&\leq&
\int_{\beta_s}\frac{|dx|}{d_{G}(x)}\leq
\int_{\beta_{1s}}\frac{|dx|}{d_{G}(x)}+
\int_{\beta_{2s}}\frac{|dx|}{d_{G}(x)}
 \\
 \nonumber&<&8a\log\left (1+\frac{\diam(\beta)}{d_{G}(w_0)}\right ),
\end{eqnarray*}
whence
\beq\label{h-l-22}
k_{G}(w_1,w_2)&\leq& \sum_{i=0}^{s_1-1}k_{G}(\eta_i,\eta_{i+1})\leq 32as_1s_2^2\log\left(1+\frac{\diam(\beta)}{d_{G}(w_0)}\right ).
\eeq
In the remaining case, that is, for each $i\in \{0,\ldots,s_1-1\}$,
$$ \ell(\beta_i)> \diam(\beta).
$$
It is possible that $x_0=\infty$. To apply Lemma \ref{lem-4-1} to continue the proof, we choose a point from each $\beta_i$
 $(i\in
\{0,\ldots,s_1-1\})$ to replace $x_0$ in the following way: It follows from $(1)$ in this lemma that for each $i\in \{0, 1,
\ldots, s_1-1\}$, there exists $v_i\in \beta_i$ satisfying
$$\frac{1}{2a}\diam(\beta)\leq d_{G}(v_i)<\diam(\beta).
$$
Then we easily obtain that for all $i\in \{0,\ldots,s_1-1\}$,
$$\mathbb{B}\big (v_i,\frac{1}{2}d_{G}(v_i)\big )\subset
\mathbb{B}\big(\eta_0,(a+\frac{3}{2})\diam(\beta)\big ),
$$
since $|v_i-\eta_0|< (1+a)\diam(\beta)$.
We see that there exist $p\not=q \in \{0,\ldots,s_1-1\}$ such that
\beq\label{eq(h-3-5')}
\mathbb{B}\big (v_p,\frac{1}{2}d_{G}(v_p)\big )\cap \mathbb{B}\big (v_q, \frac{1}{2}d_{G}(v_q)\big
)\not=\emptyset,
\eeq
because otherwise,
\begin{align*}
\Big(a+\frac{3}{2}\Big)^n {\rm Vol\,}\big (\mathbb{B}(\eta_0,
\diam(\beta))\big )
=& {\rm Vol\,}\big (\mathbb{B}(\eta_0,(a+\frac{3}{2})\diam(\beta))\big)\\
> & \Big(\frac{1}{4a}\Big)^ns_1 {\rm Vol\,}\big
(\mathbb{B}(\eta_0,\diam(\beta) )\big ),
\end{align*}
where ``${\rm Vol\,}$" denotes the volume. This is clearly a
contradiction since $s_1\geq [(4a)^{2n}]$.

\begin{figure}[!ht]
\begin{center}
\includegraphics%[width=10cm]
{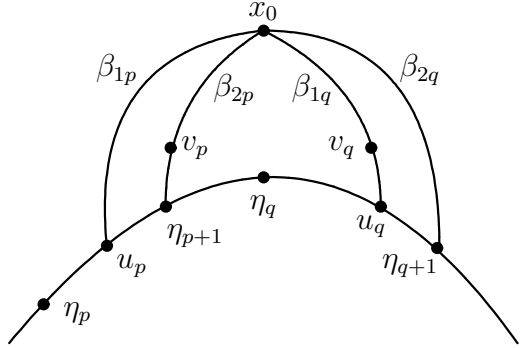}
\end{center}
\caption{The points $v_p$ and $v_q$ in the case $(1)$\label{fig7} }
\end{figure}

We divide the rest of the arguments into four cases:
\begin{enumerate}
\item\label{1-1}
$v_p\in \beta_{2p}\,\;\mbox{and}\;\, v_q\in \beta_{1q}\;\;(\mbox{see
Figure \ref{fig7}});$
\item\label{1-2}
$v_p\in \beta_{2p}\,\;\mbox{and}\;\, v_q\in \beta_{2q}$;
\item\label{1-3}
$v_p\in \beta_{1p}\,\;\mbox{and}\;\, v_q\in \beta_{1q}$;
\item\label{1-4}
$v_p\in \beta_{1p}\,\;\mbox{and}\;\, v_q\in \beta_{2q}$,
\end{enumerate}
where $p<q$.

%\begin{figure}[!ht]
%\includegraphics[width=0.40\textwidth]{figure10} %,height=0.4\textwidth
%\caption{The points $v_p$ and $v_q$ in the case $(1)$\label{fig7}
%{\color{red} Add $x_0$, $\beta_{1p}$, $\beta_{2p}$, $\beta_{1q}$ and $\beta_{2q}$}}
%\end{figure}

It suffices to discuss the first case since the discussions for the
remaining three
cases %\textbf{(b)}, \textbf{(c)} and \textbf{(d)}
%\eqref{1-2}, \eqref{1-3} and \eqref{1-4}
are similar.
First, we exploit Lemma \ref{lem-4-1} to estimate $k_{G}(\eta_{p+1}, u_q)$ in terms of $\frac{\diam(\beta)}{d_{G}(w_0)}$.
It follows from Lemma \ref{lem-4-1} that
\begin{align*}
k_{G}(v_p, \eta_{p+1})\leq & \ell_{k_{G}}(\beta_p[v_p, \eta_{p+1}])
\leq
4a^2\log\left (1+\frac{\diam(\beta)}{d_{G}(w_0)}\right ),
\end{align*}
$$k_{G}(u_q, v_q) \leq 4a^2\log\Big(1+\frac{\diam(\beta)}{d_{G}(w_0)}\Big)
$$
and by \eqref{eq(h-3-5')}, obviously,
$$k_{G}(v_p, v_q)<2.
$$
Thus
\beq\label{11-0-1}
k_{G}(\eta_{p+1}, u_q) &\leq & k_{G}(\eta_{p+1}, v_p)+k_{G}(v_p, v_q)+k_{G}(v_q, u_q)\\
\nonumber&\leq & 8a^2\log\left
(1+\frac{\diam(\beta)}{d_{G}(w_0)}\right)+2.
\eeq

Since $k_{G}(w_1, w_2)>s_1^2$, we see that
$\diam(\beta)>\frac{2}{3}d_{G}(w_0)$. For otherwise,
$$|w_1-w_0|\leq \frac{2}{3}d_{G}(w_0) ~\mbox{ and }~|w_2-w_0|\leq \frac{2}{3}d_{G}(w_0),
$$
which implies that $k_{G}(w_1, w_2)< 3$. It is impossible.
Therefore, an elementary computation shows that
$$2\leq 4a^2\log\left (1+\frac{\diam(\beta)}{d_{G}(w_0)}\right),
$$
whence \eqref{11-0-1} yields
\beq\label{11-0}
k_{G}(\eta_{p+1}, u_q)\leq  12a^2\log\left (1+\frac{\diam(\beta)}{d_{G}(w_0)}\right ).
\eeq

Next, we establish a relationship between $k_{G}(\eta_i,\eta_{i+1})$ and $k_{G}(\eta_{p+1},u_q)$ for  $i\in \{0,\ldots,s_1-1\}$.
We know from the $QH$ condition and the choice of $u_q$
that for each $i\in \{0,\ldots,s_1-1\}$, if $q>p+1$, then
$$k_{G}(\eta_i,\eta_{i+1})\leq s_2k_{G}(\eta_{p+1},\eta_{p+2})\leq s_2^2k_{G}(\eta_{p+1},u_q),
$$
and if $q=p+1$, then
$$k_{G}(\eta_i,\eta_{i+1})\leq s_2k_{G}(\eta_{p+1},\eta_{p+2})\leq
\frac{4}{4s_2-1}s_2^2k_{G}(\eta_{p+1},u_q),
$$
since $k_{G}(\eta_{p+1},u_q)\geq k_{G}(\eta_{p+1},\eta_{p+2})-k_{G}(\eta_{p+2},u_q)$.
Hence we have proved that for each $i\in \{0,\ldots,s_1-1\}$,
\beq \label{ny-1}
k_{G}(\eta_i,\eta_{i+1})\leq
s_2k_{G}(\eta_{p+1},\eta_{p+2})\leq 4s_2^2k_{G}(\eta_{p+1},u_q).
\eeq

Now, we are ready to finish the proof for this case. By (\ref{11-0}) and \eqref{ny-1}, we have
\beq\label{h-l-2}
k_{G}(w_1,w_2)&\leq& \sum_{i=0}^{s_1-1}k_{G}(\eta_i,\eta_{i+1})
\leq  48a^2s_1s_2^2\log\left(1+\frac{\diam(\beta)}{d_{G}(w_0)}\right ).
\eeq
The combination of \eqref{h-l-22} and \eqref{h-l-2} completes
the proof of $(2)$ in Theorem \ref{scl-1-1}.
\epf

%%%%%%%%%%%%%%%%%%%%%%%%%%%%%%%%%%
%%%%%%%%%%%%%%%%%%%%%%%%%%%%%%%%%%
\section{Quasiconformal subinvariance property of John domains}\label{sec-6}
%%%%%%%%%%%%%%%%%%%%%%%%%%%%%%%%%%
%%%%%%%%%%%%%%%%%%%%%%%%%%%%%%%%%%

In this section, we assume that $D$ and $D'$
are bounded subdomains in $\IR^n$, that $D'$ is an $a$-John domain with center $y'_0\in
D'$, that $f: D\to D'$ is a $K$-quasiconformal
mapping, and that $D_1\subset D$ is a $c$-John domain with center $z_0$. For any $z_1\in D_1$, we shall
construct a carrot arc in $D'_1$ to join  $z'_1$ and $z'_0$ in the
sense of ``diameter", which is stated as Theorem \ref{thm6.1}. Clearly, this approach shows that $D'_1$ is a
John domain, from which our main result in this paper, Theorem \ref{thm1.1},  will be easily proved. In order to get such a carrot arc in $D'$, we first construct a carrot arc in $D$, which is included in Lemma \ref{lem4-h-1} or Lemma \ref{lem4-1-0}, and then we shall prove that the image of the obtained carrot arc under $f$ is our desired carrot arc in $D'$. This will be reached through a series of lemmas which are divided into two groups. To this end, we
 let $\gamma_1\subset D_1$ denote a rectifiable arc with
endpoints $z_{1}$ and $z_0$ satisfying
$$\ell (\gamma_1 [z_1, z])\leq c\, d_{D_1}(z)
$$
for all $z\in \gamma_1 $, i.e. $\gamma_1$ is a $c$-carrot arc with center $z_0$.

%%%%%%%%%%%%%%%%%%%%%%%%%%%%%%%%%%
\subsection{The constructions of a uniform domain and an arc in $D_1$}%\label{sec-6-1}
%%%%%%%%%%%%%%%%%%%%%%%%%%%%%%%%%%
It follows from Lemma
\ref{lem-m-v-1} that the following result is obvious.

\begin{lem}\label{lem2.1}
There exists a simply connected domain $D_{1,0}=\bigcup
\limits_{i=1}^{k_{1}}B_{1,i}\subset D_1$ such that
\bee
\item%\label{equ--1}
$z_{1}$, $z_0\in D_{1,0}$;
\item%\label{equ--2}
for each $i\in \{1,\ldots, k_{1}\}$,
$$\frac{1}{3\rho_1}\,d_{D_1}(x_{1,i})\leq r_{1,i}\leq \frac{1}{\rho_1}d_{D_1}(x_{1,i});
$$
\item%\label{equ--3}
if $k_{1}\geq 3$, then for all $i,j\in\{1,\ldots,k_{1}\}$ with
$|i-j|\geq 2$,
$$\dist(B_{1,i}, B_{1,j})\geq \frac{1}{2^5\rho_2}\max\{r_{1,i},r_{1,j}\};
$$
\item%\label{equ--4}
if $k_{1}\geq 2$, then for each $i\in\{1,\ldots,k_{1}-1\}$,
$$\ds r_{1,i}+r_{1,i+1}-|x_{1,i}-x_{1,i+1}| \ge \frac{1}{2^5\rho_2}\max\{r_{1,i},r_{1,i+1}\},
$$
\eee
where $B_{1,i}=\mathbb{B}(x_{1,i}, r_{1,i})$, $x_{1,i}\in \gamma_{1}$,
$x_{1,i}\not\in B_{1,i-1}$ for each $i\in \{2, \ldots, k_{1}\}$, $\rho_1$ is from Lemma \ref{mxll-1''} and
$\rho_2=2^{32+c^2\rho_1^2}$.
\end{lem}

Moreover, Lemma \ref{lem2.2-4-mv-3} gives

\begin{lem}\label{lem2.2-4}
The domain $D_{1,0}$ constructed in Lemma {\rm \ref{lem2.1}} is a
$\rho_3$-uniform domain, where $\rho_3=2^{12}c^2\rho_1\rho_2$.
\end{lem}

\begin{figure}[!ht]
\includegraphics[width=0.85\textwidth]{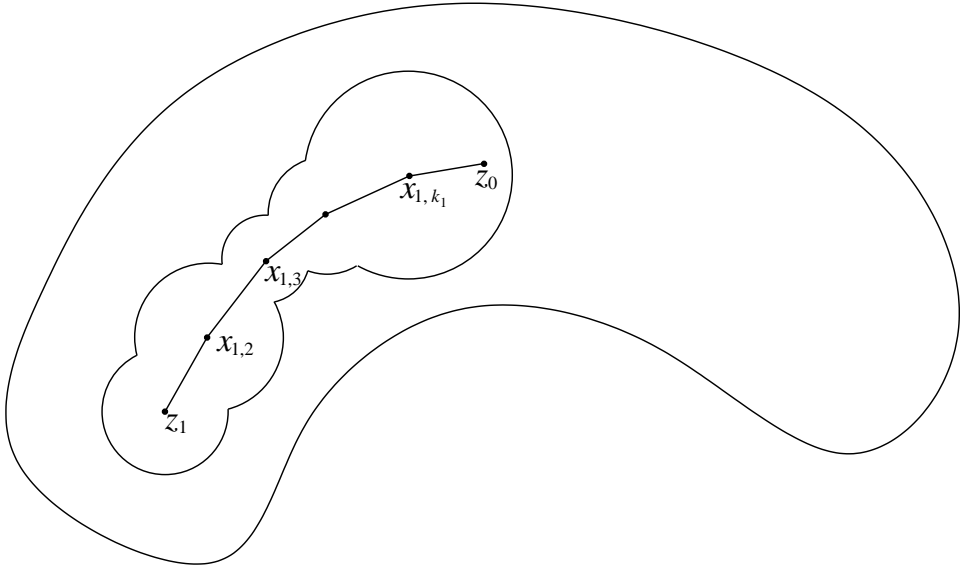} %,height=0.4\textwidth
\caption{The arc $\gamma_{1,0}$ joining $z_0$ and $z_1$
\label{fig2}
}
\end{figure}

Now we set
$$\gamma_{1,0}=[z_{1},x_{1,2}]\cup\cdots\cup[x_{1,k_1-1},x_{1,k_1}]\cup[x_{1,k_1},z_0]\;\,
 (\mbox{see Figure \ref{fig2}}).
 $$
Obviously, $\gamma_{1,0}\subset D_{1,0}$. In the rest of this section,
our main aim is to prove that the image $\gamma'_{1,0}$ of $\gamma_{1,0}$ under $f$ is the
desired carrot arc in $D'$. First, we need properties on $\gamma_{1,0}$, $D_{1,0}$ and $D'_{1,0}$,
which are included in the next two subsections.

%%%%%%%%%%%%%%%%%%%%%%%%%%%%%%%%%%
\subsection{The properties of $\gamma_{1,0}$ }
%%%%%%%%%%%%%%%%%%%%%%%%%%%%%%%%%%

In this subsection, we show several lemmas on $\gamma_{1,0}$ which will be used later on.

The following lemma shows that $\gamma_{1,0}$ is a carrot arc in $D_1$ with
center $z_0$.

\begin{lem}\label{lem4-h-1}
For each $z\in \gamma_{1,0}$, we have $\ell(\gamma_{1,0}[z_1,
z])\leq \rho_5 d_{D_1}(z)$, where $\rho_5=\frac{8}{7}c$.
\end{lem}\bpf
Clearly, for every $z\in \gamma_{1,0}$, there exists an $i\in\{1,\ldots,k_1\}$ such that
$z\in B_{1,i}= \mathbb{B}(x_{1,i}, r_{1,i})$. It follows from Lemma \ref{lem2.1} that
$$\ell(\gamma_{1,0}[z_1, z])\leq \ell(\gamma_1[z_1, x_{1,i}])+r_{1,i}
\leq c d_{D_1}(x_{1,i})+r_{1,i}\leq \rho_5 d_{D_1}(z),
$$
since $d_{D_1}(x_{1,i})\leq |x_{1,i}-z|+d_{D_1}(z)$ and $r_{1,i}=\frac{1}{\rho_1}d_{D_1}(x_{1,i})$.
Hence the lemma holds.
\epf

 Further, from Lemma \ref{lem2.1} and the similar reasoning as in the proof of
\cite[Theorem $1.8$]{HW}, we have the following.

\begin{lem}\label{lem4-1-0}
Suppose that $x_1\in \gamma_{1,0}$. Then for $x_2\in \gamma_{1,0}[x_1, z_0]$,
the part $\gamma_{1,0}[x_1,x_2]$ of $\gamma_{1,0}$ is a double
$(2^3c\rho_1)^2(2^6\rho_2+1)$-cone arc in $D_{1,0}$.
\end{lem}

The next two results show that every part of $\gamma_{1,0}$ is
solid in $D$ and in $D_{1}$, respectively.

\begin{lem}\label{lem4-1-6}
For $x_1$, $x_2\in \gamma_{1,0}$,
$\ell_{k_{D}}(\gamma_{1,0}[x_1,x_2])\leq
2^{16}c^3\rho_1\rho_2k_{D}(x_1,x_2)$, where $\rho_1$ and $\rho_2$ are from Lemma \ref{lem2.1}.
\end{lem}\bpf
In view of Lemmas \ref{lem-4-1}, \ref{lem4-1-0} and \ref{lem4-h-1}, we have
\begin{eqnarray*}
\ell_{k_{D}}(\gamma_{1,0}[x_1, x_2])&=& \int_{\gamma_{1,0}[x_1, x_2]}\frac{|dx|}{d_{D}(x)}
\leq 2\rho_5 \log\left(1+\frac{2\ell(\gamma_{1,0}[x_1, x_2])}{d_{D}(x_1)}\right ) \\
\nonumber&<&
2^{16}c^3\rho_1\rho_2k_{D}(x_1,x_2)
\end{eqnarray*}
as required.
\epf

Similarly, we have

\begin{lem}\label{lem4-1-6-00'}
For $u_1$, $u_2\in \gamma_{1,0}$, $\ell_{k_{D_1}}(\gamma_{1,0}[u_1,u_2])
\leq 2^{16}c^3\rho_1\rho_2k_{D_1}(u_1,u_2)$.
\end{lem}

We easily infer from Lemmas  \ref{lem-4-1} and \ref{lem4-h-1} the following which will be used later on.

\begin{lem}\label{lem-4-1'}
Suppose $x_1$ and $x_2\in \gamma_{1,0}$. Then for every $u\in
\gamma_{1,0}[x_1, x_2]$, we have

%\noindent $(1)$
\begin{enumerate}
\item[{\rm (1)}]
$\ds d_{D_1}(u)\geq \frac{2\ell(\gamma_{1,0}[x_1, u])+d_{D_1}(x_1)}{4\rho_5}$, and

%\noindent $(2)$
\item[{\rm (2)}]
$\ds d_{D}(u)\geq \frac{2\ell(\gamma_{1,0}[x_1, u])+d_{D}(x_1)}{4\rho_5}.$
\end{enumerate}
\end{lem}

%%%%%%%%%%%%%%%%%%%%%%%%%%%%%%%%%%
\subsection{The properties of the domain $D_{1,0}$ and its image $D'_{1,0}$}
%%%%%%%%%%%%%%%%%%%%%%%%%%%%%%%%%%

In the proof of the main result in this section, Theorem \ref{thm6.1}, the following two results on $D_{1,0}$ and $D'_{1,0}$ are useful.

\begin{lem}\label{lem4-h-0'}
For each $z\in D_{1,0}$, we have
$$d_{D_{1}}(z)\geq (\rho_1-1)d_{D_{1,0}}(z),
$$
and for each $z\in \gamma_{1,0}[z_1,x_{1,k_1}]$,

%\noindent $(1)$
\begin{enumerate}
\item[{\rm (1)}]
$d_{D_{1}}(z)\leq 2^6(3\rho_1+1)\rho_2d_{D_{1,0}}(z)$;

%\noindent $(2)$
\item[{\rm (2)}]
$d_{D'_{1}}(z')\leq \big (2^7(3\rho_1+1)\rho_2\mu_2\big )^{\mu_2}d_{D'_{1,0}}(z')$,
\end{enumerate}
where $\mu_2$ is from Theorem {\rm \Ref{ThmF}}.
\end{lem} \bpf
For $z\in D_{1,0}$, by Lemma \ref{lem2.1}, there exists a $j\in
\{1,\ldots, k_1\}$ such that $z\in B_{1,j}$. Further, we assume that $B_{1,i}$ is the last ball from
$B_{1,1}$ to $B_{1,k_1}$ such that $z\in B_{1,i}$. Then
$$d_{D_1}(z)\geq d_{D_1}(x_{1,i})-|z-x_{1,i}|\geq (\rho_1-1)d_{D_{1,0}}(x_{1,i})
\geq (\rho_1-1)d_{D_{1,0}}(z),
$$
from which the first inequality follows.

For the inequality $(1)$, let $z\in \gamma_{1,0}[z_1,x_{1,k_1}]$. We still assume that
$B_{1,i}$ is the last ball from $B_{1,1}$ to $B_{1,k_1}$ such that
$z\in B_{1,i}$. Again, Lemma \ref{lem2.1} implies
$$ d_{D_{1,0}}(z)\geq \frac{1}{2^6\rho_2}r_{1,i}
$$
and
$$d_{D_1}(z)\leq d_{D_1}(x_{1,i})+r_{1,i}\leq (3\rho_1+1)r_{1,i},
$$
whence
$$d_{D_1}(z)\leq 2^6(3\rho_1+1)\rho_2d_{D_{1,0}}(z),
$$
which shows that $(1)$ holds.

For the remaining inequality, let $u\in\mathbb{S}\big(z, d_{D_{1,0}}(z)\big)$. Then we have
$$\log\left (1+\frac{1}{2^6(3\rho_1+1)\rho_2}\right )\leq \log\left (1+\frac{|u-z|}{d_{D_1}(z)}\right)
\leq  k_{D_1}(z,u)\leq  \int_{[z,u]}\frac{|dw|}{d_{D_1}(w)}\leq
\frac{1}{\rho_1-2},
$$
since for $w\in [z, u]$, $d_{D_1}(w)\geq d_{D_1}(z)-|w-z|$ and $d_{D_1}(z)\geq (\rho_1-1)d_{D_{1,0}}(z)$,
whence we see from Theorem \Ref{ThmF} that
$$k_{D'_1}(z',u')\geq \left (\frac{1}{\mu_2}\log\Big(1+\frac{1}{2^6(3\rho_1+1)\rho_2}\Big)\right )^{\mu_2}>
\frac{1}{\big((2^7(3\rho_1+1)\rho_2)\mu_2\big)^{\mu_2}-1}.
$$
Necessarily, we have
$$u'\in \IR^n \backslash \overline{\mathbb{B}}\Big(z', \frac{1}{(2^7(3\rho_1+1)\rho_2\mu_2)^{\mu_2}}d_{D'_{1}}(z')\Big),
$$
then the proof of $(2)$ easily follows from the fact that
$f\big(\overline{\mathbb{B}}(z, d_{D_{1,0}}(z))\big)\subset
\overline{D'_{1,0}}$. \epf

%Formulation of our next result needs the following result.

%\begin{Thm} \label{lem-j-j} $($\cite[Lemma 2.1]{HLVW}$)$
%Suppose that $G$ is a $c$-uniform domain in $\IR^n$, and for $x$ and $y$ in $G$, let $\beta$ denote
%a $(\nu,h)$-solid arc in $G$ with the endpoints $x$ and $y$. If $d_{G}(x_0)=\sup
%\{d_{G}(p):\, p\in \beta\}$, then for all  $z\in \beta[x, x_0]$,
%$$|x-z|\leq \rho_4\;d_{G}(z),
%$$
%and for all  $z\in \beta[y, x_0]$,
%$$|y-z|\leq \rho_4\;d_{G}(z),
%$$
%where $\rho_4=\rho_4(\nu, h, c)$.
%\end{Thm}

It follows from \cite[Theorem 3]{Geo} and \cite[Theorem 4.15]{Vai4} that solid arcs have the quasiconformal invariance property, i.e., the image of each solid arc under a conformal mapping is still a solid arc.
Since each quasihyperbolic geodesic is a $(0,1)$-solid arc,
we easily see that the following corollary is a consequence of Lemma
\ref{lem2.2-4}, \cite[Lemma 2.1]{HLVW}, Theorems \Ref{ThmF'} and \Ref{Thm4-1}.

\begin{cor}\label{Cor-1}
There exists a constant $\mu_8$ such that for $w_1$ and $w_2\in D_{1, 0}$, the following hold:
\bee
\item[$(1)$] $k_{D_{1,0}}(w_1, w_2)\leq \mu_8 \log \left
(1+\frac{|w_1-w_2|}{\min\{d_{D_{1,0}}(w_1), d_{D_{1,0}}(w_2)\}}\right ). $

Suppose that $\zeta'$ is a quasihyperbolic geodesic
joining $w'_1$ and $w'_2$ in $D'_{1,0}$. Then

\item[$(2)$] $\zeta$ is a $(\nu, h_1)$-solid arc and
$$\diam(\zeta)\leq \mu_8\max\{|w_1-w_2|,\min\{d_{D_{1,0}}(w_1),d_{D_{1,0}}(w_2)\}\},
$$
where $(\nu, h_1)$ depends only on $(n, K)$, and

\item[$(3)$] for all $w\in \zeta$, $\min\Big\{\diam(\zeta [w_1,
w]),\diam(\zeta [w, w_2])\Big\}\leq \mu_8\,d_{D_{1,0}}(w);$

\item[$(4)$] Let $v_{0,0}\in \zeta$ be such that
$$d_{D_{1,0}}(v_{0,0})=\sup\limits_{p\in \zeta}d_{D_{1,0}}(p).
$$
\eee
Then for $w\in \zeta[w_j, v_{0,0}]$ $(j=1,2)$,
$$|w_j-w|\leq \mu_8\;d_{D_{1,0}}(w),
$$
where $\mu_8=\mu_8(n,K,\rho_3, \nu, h_1)=\mu_8(n,K,\rho_3)$ and $\rho_3$ is from Lemma \ref{lem2.2-4}.
\end{cor}

Now, we are in a position to state and prove our main result in this section.

%%%%%%%%%%%%%%%%%%%%%%%%%%%%%%%%%%
\subsection{The statement of Theorem \ref{thm6.1}} %\label{sec-6-1}
%%%%%%%%%%%%%%%%%%%%%%%%%%%%%%%%%%

The following is the main result in this section which shows that $\gamma'_{1,0}$ is a carrot arc with center $z'_0$.

\begin{thm}\label{thm6.1}
There is a constant $\rho$ such that for any $z'\in \gamma'_{1,0}$, $$\diam(\gamma'_{1,0}[z'_1,z'])<
\rho d_{D'_{1}}(z'),$$ where $\rho=\rho\left(n, K, a, c, \frac{\diam(D)}{d_D(f^{-1}(y'_0))}\right)$.
\end{thm}

%%%%%%%%%%%%%%%%%%%%%%%%%%%%%%%%%%
\subsection{The proof of Theorem \ref{thm6.1}} %\label{sec-6-1}
%%%%%%%%%%%%%%%%%%%%%%%%%%%%%%%%%%

In this subsection, by a method of contradiction, we shall prove that the constant $\rho$ in Theorem \ref{thm6.1} can be taken to be $2\rho_{10}$.
Here
\bee
\item[$(1)$]  $\ds \rho_{10}=\rho_8^{\rho_9}$,

\item[$(2)$]  $\ds \rho_9=\rho_8^{\Gamma(\Lambda_{T_1+3})}$,

\item[$(3)$]   $\ds \rho_8=48a^2\mu_2^6\rho_7^8$,

\item[$(4)$]   $\ds \rho_7=2^{12}a^{2n}c^2\big (\mu_8\rho_2\rho_6\big )^{2\mu_2\mu_8}$,

\item[$(5)$]   $\ds \rho_6=\frac{2^{12}\mu_8^2(3\rho_1+1)\rho_2\rho_5}{\eta_1^{-1}(\frac{1}{\rho_5})}$,
\eee
where we need to remember the following:
\bee
\item[(i)]  $\mu_2$ (resp. $\mu_8$) is from Theorem \Ref{ThmF} (resp. Corollary \ref{Cor-1}); $\rho_1$ and $\rho_2$ are from Lemma \ref{lem2.1} and $\rho_3$ (resp. $\rho_5$)  is from Lemma
\ref{lem2.2-4} (resp. Lemma
\ref{lem4-h-1});

\item[(ii)]   in $\rho_6$, $\eta_1=\eta_{n, K, \lambda_1}$,
$\lambda_1=\frac{1}{5}$ and $\eta$ is from Theorem
\Ref{Vai-0};

\item[(iii)] in $\rho_9$, $\Lambda_0=\rho_8^2$ and for each $i\in\{1,\ldots,T_1+3\}$ ($T_1=[16a]^{2n}$),
$\Lambda_i= \Gamma(\Lambda_{i-1})$, here
$$\Gamma(t)=\Big((1+216\rho_7t)t\rho(t)H(t)\Big){^{2\mu_2^2\mu_8\mu(t)}};
$$

\item[(iv)]  in $\Gamma(t)$, $\rho(t)=\rho(n,
K, \rho_3, 2^{47+(3a\rho_7t)^2}a^2\rho_7t^2)$ and $\rho$ is from theorem \Ref{ThmA};

\item[(v)]  in $\Gamma(t)$, $\mu(t)=\mu_1(c_1(t))$, $c_1(t)=2^{47+(3a\rho_7t)^2}a^2\rho_7t^2$
 and $\mu_1$ is from Theorem \Ref{ThmF'};

\item[(vi)]  in $\Gamma(t)$, $H(t)$ is defined as follows:
$$H(t)=\kappa_2\left (n,K,\mu_{0t},\varphi(t),6a(\rho_7+1),
\rho_7^{2\mu_2}\frac{\diam(D)}{d_D(y_0)}\right );
$$

\item[(vii)] in $H(t)$, $\varphi(t)=\varphi_{c_2(t)}$, $c_2(t)=2^{47+(3a\rho_7t)^2}a^2\rho_7t^2$ and $\varphi$ is
from Theorem \Ref{Lem4-1};

\item[(viii)] in $H(t)$, $\mu_{0t}=\mu_3(n,\varphi_{\rho(t)})$ and $\mu_3$ is from Theorem \Ref{Lem6-1};

\item[(ix)] in $H(t)$, $\kappa_2$ is from Theorem \Ref{Lem6-1'}.
\eee

Suppose on the
contrary that there exists a point $w'\in \gamma'_{1,0}$ satisfying
\beq\label{eq(h-3-2)}
\diam(\gamma'_{1,0}[z'_1,w'])\geq 2\rho_{10}
d_{D'_{1}}(w').
\eeq

Based on this assumption, we shall prove a series of lemmas. Through these lemmas, a contradiction will be reached, from which the proof of Theorem \ref{thm6.1} is complete. In the following, we divide these lemmas into two groups.
\medskip

%%%%%%%%%%%%%%%%%%%%%%%%%%%%%%%%%%%
%%%%%%%%%%%%%%%%%%%%%%%%%%%%%%%%%%%
\subsubsection{{\bf Lemmas: Part I}}\ \ \
%%%%%%%%%%%%%%%%%%%%%%%%%%%%%%%%%%%
%%%%%%%%%%%%%%%%%%%%%%%%%%%%%%%%%%%
\medskip

{\bf The sketch for this part}:\quad
In this part, under the assumption of \eqref{eq(h-3-2)}, we
will pick up two special points $z'_{1,0}$ and
 $z'_{0,0}$ from the image $\gamma'_{1,0}$ of $\gamma_{1,0}$
 constructed in Section \ref{sec-6}.
Then these two points determine a quasihyperbolic geodesic
$\gamma'_{2,0}$ in $D'_{1,0}$. We shall apply the obtained results,
especially  Theorem \ref{scl-1-1}, to analyze $\gamma_{2,0}$ and
$\gamma'_{2,0}$ together with some related points in these two arcs, and
seven Lemmas will be proved.\medskip

\begin{figure}[!ht]
\includegraphics%[width=10cm]
{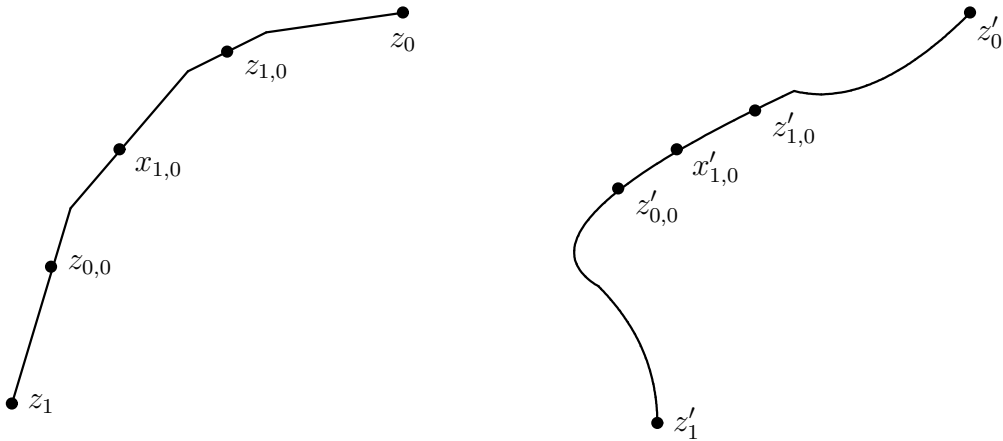}
\caption{The points $z_{0,0}$, $z_{1,0}$, $x_{1,0}$ in
$\gamma_{1,0}$ and their images under $f$ in
$\gamma'_{1,0}$\label{fig8}
%{\color{red} The position of $z'_{0,0}$ should be lifted to $x'_{1,0}$   }
}
\end{figure}

We let $z'_{1,0}$ be the first point in $\gamma'_{1,0}$ along the
direction from $z'_1$ to $z'_0$ such that
\beq\label{eq(h-3-3)}
\diam(\gamma'_{1,0}[z'_1,z'_{1,0}])= \rho_{10} d_{D'_{1}}(z'_{1,0}).
\eeq
The existence of $z'_{1,0}$ follows from \eqref{eq(h-3-2)}. And let
$$z'_{0,0}\in \gamma'_{1,0}[z'_1, z'_{1,0}]\cup \mathbb{S}\big(z'_{1,0},
\frac{1}{4}\diam(\gamma'_{1,0}[z'_1,z'_{1,0}])\big)
$$
be such that $\gamma'_{1,0}[z'_{0,0},z'_{1,0}]\subset\overline{\mathbb{B}}(z'_{1,0},
\frac{1}{4}\diam(\gamma'_{1,0}[z'_1,z'_{1,0}]))$. Obviously,
$z'_{0,0}$ satisfies the following  (see Figure \ref{fig8}):
\beq \label{24-1}
2\diam(\gamma'_{1,0}[z'_{0,0},z'_{1,0}])\leq
\diam(\gamma'_{1,0}[z'_1,z'_{1,0}])\leq
2\diam(\gamma'_{1,0}[z'_1,z'_{0,0}])
\eeq
and
%{\color{red} The place for Figure 11:.}
\beq \label{24-1-mei}
|z'_{1,0}-z'_{0,0}|=\frac{1}{4}\diam(\gamma'_{1,0}[z'_1,z'_{1,0}]).
%\;(\mbox{see Figure \ref{fig8}}).
\eeq

Next, we determine the positions of points $z_{0,0}$ and
$z_{1,0}$ in $\gamma_{1,0}$.
\begin{lem}\label{cl8-h-1} The points $ z_{0,0}$
and $z_{1,0}$ are not contained in the part
$\gamma_{1,0}[x_{1,k_1},z_0]$ of $\gamma_{1,0}$, where the point
$x_{1,k_1}$ is from Lemma \ref{lem2.1} $($see Figure \ref{fig2}$)$.
\end{lem}

Now, we  prove this Lemma. For each
$x\in\gamma_{1,0}[x_{1,k_1},z_0]$, by Lemma \ref{lem2.1}, we see
that $x\in B_{1, k_1}=\mathbb{B}(x_{1, k_1},
\frac{1}{\rho_1}d_{D_1}(x_{1, k_1}))$, and thus, Lemma
\ref{mxll-1''} implies
\be\label{00}
|x'-x'_{1,k_1}|< \frac{1}{5}d_{D_1'}(x'_{1,k_1})\;\;\mbox{and}\;\;d_{D'_1}(x'_{1,k_1})
\leq \frac{5}{4}d_{D'_1}(x').
\ee

Suppose $z_{0,0}\in \gamma_{1,0}[x_{1,k_1},z_0]$ or $z_{1,0}\in
\gamma_{1,0}[x_{1,k_1},z_0]$. Then for each
$w'\in\gamma'_{1,0}[z'_{1,0},z'_0]$,  we deduce from the choice of
$z'_{1,0}$ and (\ref{00}) that
\begin{eqnarray*}
\diam(\gamma'_{1,0}[z'_1,w'])&\leq& \diam(\gamma'_{1,0}[z'_1,x'_{1,k_1}])
+\diam(\gamma'_{1,0}[x'_{1,k_1},z'_0])\\ \nonumber&\leq&
\rho_{10}d_{D'_1}(x'_{1,k_1})+\frac{2}{5}d_{D'_1}(x'_{1,k_1}) \leq \frac{2+5\rho_{10}}{4}d_{D'_1}(w')\\
\nonumber&<& 2\rho_{10} d_{D'_{1}}(w'),
\end{eqnarray*}
which, together with the choice of $z'_{1,0}$, shows that for all $w'\in
\gamma'_{1,0}$, $\diam(\gamma'_{1,0}[z'_1,w'])< 2\rho_{10}
d_{D'_{1}}(w').$  It follows from \eqref{eq(h-3-2)} that this is impossible. Hence Lemma \ref{cl8-h-1} is true.
\medskip

Further, we have

\begin{lem}\label{co1} $(1)$ \label{eq-h-1}
For $z\in\gamma_{1,0}[z_{0,0},z_{1,0}]$,
$$d_{D}(z)< \rho_1\ell(\gamma_{1,0}[z_1,z])\;\;\mbox{and} \;\;d_{D}(z)<
\rho_1\rho_5 d_{D_1}(z);
$$

$(2)$ \label{cl5-h-1} For $x_1, x_2\in \gamma_{1,0}[z_{0,0},
z_{1,0}]$,
$$k_{D_1}(x_1,x_2)\leq 2^6\rho_2 k_{D}(x_1,x_2).
$$
\end{lem}

First, we prove $(1)$. Suppose on the contrary that there exists
some point $z\in\gamma_{1,0}[z_{0,0},z_{1,0}]$ such that
%\be\label{eq3-12-0}
$$d_{D}(z)\geq \rho_1\ell(\gamma_{1,0}[z_1,z]).
$$
We shall use the local quasisymmetry of $f$ to get a contradiction. We first do some preparation.

Let $x_{1,0}$ be the first point of $\gamma_{1,0}[z_{0,0},z_{1,0}]$
from $z_{0,0}$ to $z_{1,0}$  satisfying (see Figure \ref{fig8})
%\be\label{eq3-12}
$$d_{D}(x_{1,0})\geq \rho_1\ell(\gamma_{1,0}[z_1,x_{1,0}]).
%\;\;(\mbox{see Figure \ref{fig8}}).
$$
It is possible that $x_{1,0}=z_{0,0}$ or $z_{1,0}$. For
every $y\in \gamma_{1,0}[z_1, x_{1,0}]$, since $y\in
\overline{\mathbb{B}}(x_{1,0}, \frac{1}{\rho_1}d_{D}(x_{1,0}))$,  by
Lemma \ref{mxll-1''}, we get
$$|y'-x'_{1,0}| < \frac{1}{5} d_{D'}(x'_{1,0}),
$$
which yields
$$\gamma'_{1,0}[z'_1, x'_{1,0}]\subset \IB \big(x'_{1,0}, \frac{1}{5}d_{D'}(x'_{1,0})\big),
$$
and so by (\ref{24-1}),
\be\label{sum-1}\diam(\gamma'_{1,0}[x'_{1,0}, z'_{1,0}])\leq \diam(\gamma'_{1,0}[z'_1, x'_{1,0}])\leq\frac{2}{5}d_{D'}(x'_{1,0}).
\ee

\begin{figure}[!ht]
\includegraphics[width=0.45\textwidth]{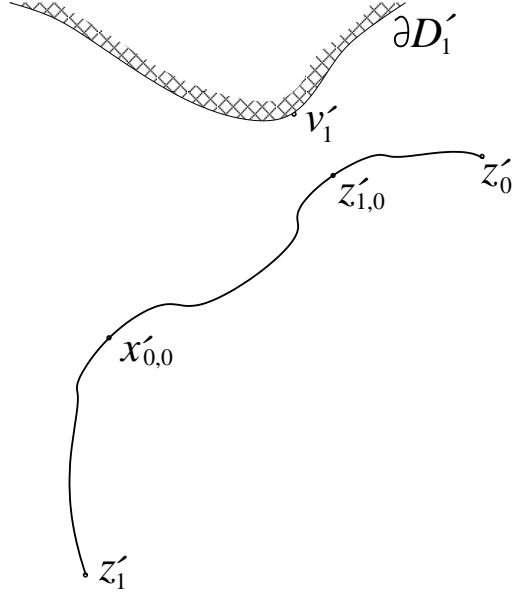} %,height=0.4\textwidth
\caption{The points $z'_{1,0}$ in
$\gamma'_{1,0}[z'_1,z'_{1,0}]$ and $v'_1$ in $\partial D'_1$\label{fig9}}
\end{figure}

Let
$$B'_1=\mathbb{B}\Big(x'_{1,0},\frac{3}{5}d_{D'}(x'_{1,0})\Big),
$$
 and let $v'_1$ be a point
in the boundary of $D'_1$ such that
\be\label{sun-2}
|z'_{1,0}-v'_1|=d_{D'_1}(z'_{1,0})\;\;(\mbox{see Figure
\ref{fig9}}).
\ee

It follows from (\ref{eq(h-3-3)}) and (\ref{sum-1}) that
$$d_{D'_1}(z'_{1,0})=\frac{1}{\rho_{10}}\diam(\gamma'_{1,0}[z'_1, z'_{1,0}])\leq \frac{3}{5\rho_{10}}d_{D'}(x'_{1,0}).$$

Then we easily know that $v'_1, z'_{0,0}, z'_{1,0}\in
B'_1.$
We are ready to get a contradiction by using the local quasisymmetry of $f^{-1}$.
By Theorem \Ref{Vai-0}, the restriction $f^{-1}|_{B'_1}$ is
$\eta_1$-quasisymmetric, where $\eta_1=\eta_{n, K, \lambda_1}$ and
$\lambda_1=\frac{3}{5}$. It follows from Lemma
\ref{lem4-h-1}, (\ref{eq(h-3-3)}), \eqref{24-1-mei} and
\eqref{sun-2} that
$$\frac{1}{\rho_5}\leq \frac{|v_1-z_{1,0}|}{|z_{0,0}-z_{1,0}|}
\leq\eta_1\Big(\frac{|v'_1-z'_{1,0}|}{|z'_{0,0}-z'_{1,0}|}\Big)\leq\eta_1(\frac{4}{\rho_{10}}).
$$
This is the  desired contradiction. Hence the first inequality in
Lemma \ref{eq-h-1}$(1)$ is true. The second inequality in Lemma
\ref{eq-h-1}$(1)$ easily follows from the first one and Lemma
\ref{lem4-h-1}.

In order to give a proof of $(2)$ in Lemma \ref{co1}, we divide
the discussions into two cases: $|x_2-x_1|\leq
\frac{1}{2}d_{D_1}(x_1)$ and $|x_2-x_1|> \frac{1}{2}d_{D_1}(x_1)$.

For the first case, Lemma \ref{eq-h-1}$(1)$ implies that
\beq\label{eq3-12-1'}
k_{D_1}(x_1,x_2)&\leq& \int_{[x_1,x_2]}\frac{|dx|}{d_{D_1}(x)} \leq
\frac{2|x_2-x_1|}{d_{D_1}(x_1)} <   3\rho_1\rho_5\,k_{D}(x_1,x_2).
\eeq
For the second case, that's, $|x_2-x_1|>\frac{1}{2}d_{D_1}(x_1)$,  it
follows from Lemmas \ref{lem4-1-0},
\ref{lem-4-1'} and \ref{eq-h-1}$(1)$ that
\beq \label{eq3-12-1}
k_{D_1}(x_1,x_2)&\leq&
\int_{\gamma_{1,0}[x_1,x_2]}\frac{|dx|}{d_{D_1}(x)}
\leq 2\rho_5\log\left(1+\frac{2\ell(\gamma_{1,0}[x_1,x_2])}{d_{D_1}(x_1)}\right) \\
\nonumber &\leq& 2\rho_5\  \log\left (1+
2^{14}c^2\rho_1^2\rho_2\rho_5\frac{|x_2-x_1|}{ d_{D}(x_1)}\right)\\
\nonumber
  &<& 2^6\rho_2\,k_{D}(x_1,x_2),
\eeq
where, in the third inequality, the inequalities
$$\frac{|x_2-x_1|}{d_D(x_1)}\geq\frac{|x_2-x_1|}{\rho_1\rho_5d_{D_1}(x_1)}>\frac{1}{2\rho_1\rho_5}
$$
are used. The combination of \eqref{eq3-12-1'} and \eqref{eq3-12-1}
completes the proof of the inequality $(2)$ in Lemma \ref{eq-h-1}.
Hence the proof of Lemma \ref{eq-h-1} is complete.

\medskip

In order to apply Theorem \ref{scl-1-1} to continue the proof, the
following lower bound for the quasihyperbolic distance in $D'_1$
between $z'_{0,0}$ and $z'_{1,0}$ is useful.

\begin{lem}\label{cl6-h-3''} $ k_{D'_1}(z'_{0,0}, z'_{1,0})\geq \log
\frac{\rho_{10}}{4}\;$ and $\;\min\{|z_{0,0}-z_{1,0}|,
d_{D_{1,0}}(z_{1,0})\}\geq \rho_9d_{D_{1,0}}(z_{0,0})$. \end{lem}

The proof of the first inequality is obvious, since it follows from
(\ref{eq(h-3-3)}) and (\ref{24-1-mei}) that
$$k_{D'_1}(z'_{0,0}, z'_{1,0})\geq \log
\frac{|z'_{0,0}-z'_{1,0}|}{d_{D'_1}(z'_{1,0})}=\log
\frac{\rho_{10}}{4}.
$$
For the second inequality, we infer from Theorem \Ref{ThmF} and the proved inequality that $k_{D_1}(z_{0,0}, z_{1,0})>1$, and thus
$$k_{D_1}(z_{0,0}, z_{1,0})\geq \frac{1}{\mu_2}k_{D'_1}(z'_{0,0},
z'_{1,0}) \geq\frac{1}{\mu_2}\log \frac{\rho_{10}}{4}.
$$
Consequently, Corollary \ref{Cor-1}, Lemmas \ref{lem4-h-1},
\ref{lem4-h-0'} and \ref{cl8-h-1} yield
\begin{eqnarray*}
\frac{1}{\mu_2}\log
\frac{\rho_{10}}{4}&\leq & k_{D_{1,0}}(z_{0,0},z_{1,0})\leq
\mu_8\log\Big(1+\frac{|z_{0,0}-z_{1,0}|}{\min\{d_{D_{1,0}}(z_{0,0}),d_{D_{1,0}}(z_{1,0})\}}\Big)\\
\nonumber&\leq&
\mu_8\log\Big(1+\frac{2^6(3\rho_1+1)\rho_2\rho_5d_{D_{1,0}}(z_{1,0})}{\min\{d_{D_{1,0}}(z_{0,0}),d_{D_{1,0}}(z_{1,0})\}}\Big),
\end{eqnarray*}
and necessarily, we see
$$\min\{d_{D_{1,0}}(z_{0,0}),d_{D_{1,0}}(z_{1,0})\}=d_{D_{1,0}}(z_{0,0}),
$$
whence $$\min\{|z_{0,0}-z_{1,0}|, d_{D_{1,0}}(z_{1,0})\}\geq
\rho_9d_{D_{1,0}}(z_{0,0}).
$$
The proof of Lemma \ref{cl6-h-3''} is complete.
\medskip

\begin{figure}[!ht]
\begin{center}
\includegraphics%[width=10cm]
{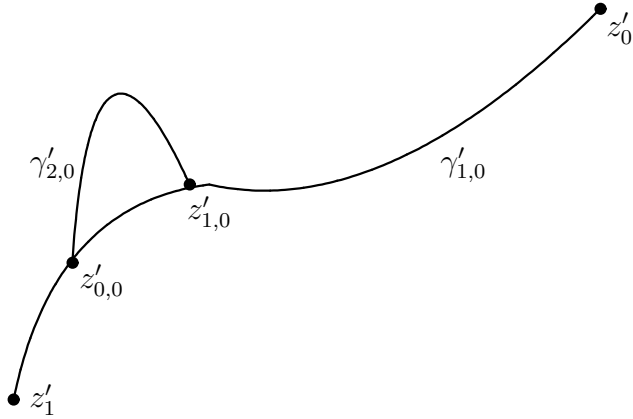}
\end{center}
\caption{The arc $\gamma'_{2,0}$ in $D'_{1,0}$\label{fig10}}
\end{figure}

Let $\gamma'_{2,0}$ be a quasihyperbolic geodesic joining $z'_{0,0}$
and $z'_{1,0}$ in $D'_{1,0}$ (see Figure \ref{fig10}). By
Proposition \ref{scl-0-1}, we see that $\gamma'_{2,0}$ satisfies the
$QH\, (\ell_{k_{D'_{1,0}}}(\gamma'_{2,0}), 2, 1)$-condition. Further,
we shall prove that if the quasihyperbolic distance between two
points in $\gamma'_{2,0}$ is big enough, then the part of
$\gamma'_{2,0}$ with these two points as its endpoints satisfies
some much ``stronger" $QH$ condition. This is included in Lemma
\ref{cl8-0''}. Before the statement and the proof of Lemma
\ref{cl8-0''}, we need some inequalities, which are stated in Lemmas
\ref{lem5-h-0'} and \ref{cl5-h-w-1}. We recall that $D_{1,0}$ is a
$\rho_3$-uniform domain (see Lemma \ref{lem2.2-4}).

%\begin{figure}[!ht]
%\includegraphics[width=0.45\textwidth]{figure13} %,height=0.4\textwidth
%\caption{The arc $\gamma'_{2,0}$ in $D'_{1,0}$\label{fig10}}
%\end{figure}

\begin{lem}\label{lem5-h-0'} For $z\in \gamma_{2,0}$, we have

\begin{enumerate}
\item\label{q-1-0}
$(\rho_1-1)d_{D_{1,0}}(z)\leq d_{D_{1}}(z)\leq
2^5\rho_2\rho_6d_{D_{1,0}}(z)$;
\item\label{q-1-1} $d_{D}(z)\leq
2^5\rho_1\rho_2\rho_6d_{D_{1,0}}(z)$;
\item\label{q-0}
$d_{D_1}(z_{0,0})< 2^9\mu_8\rho_2d_{D_{1}}(z)$;
\item\label{q-2-0}
$d_{D'_1}(z')\leq (2^6\rho_2\rho_6\mu_2)^{\mu_2}d_{D'_{1,0}}(z')$;
\item\label{q-2-1} $d_{D'}(z')\leq
(2^6\rho_1\rho_2\rho_6\mu_2)^{\mu_2}d_{D'_{1,0}}(z')$.
\end{enumerate}
\end{lem}

For a proof of this Lemma, it follows from the fact
``$\gamma_{2,0}\subset D_{1,0}$" and  Lemma \ref{lem4-h-0'} that for
$z\in \gamma_{2,0}$,
\be\label{13-1}
d_{D_{1}}(z)\geq(\rho_1-1)d_{D_{1,0}}(z),
\ee
and we obtain from Corollary \ref{Cor-1} that for $z\in \gamma_{2,0}$,
\be\label{mlm-3-3}
\min\{|z_{0,0}-z|, |z_{1,0}-z|\}\leq \mu_8\;d_{D_{1,0}}(z),
\ee
where $\mu_8=\mu_8(n,K,\rho_3)$ is from Corollary \ref{Cor-1}.

Let $z$ be a point in $\gamma_{2,0}$. To prove that $z$ satisfies the inequalities
\eqref{q-1-0} and \eqref{q-1-1} in Lemma \ref{lem5-h-0'}, we only need to
consider the case $\min\{|z_{0,0}-z|, |z_{1,0}-z|\}=|z_{0,0}-z|$
since the proof for the other case $\min\{|z_{0,0}-z|,
|z_{1,0}-z|\}=|z_{1,0}-z|$ is similar.

Further, we distinguish two possibilities: $|z_{0,0}-z|\leq \frac{1}{2}d_{D_{1,0}}(z_{0,0})$ and $|z_{0,0}-z|> \frac{1}{2}d_{D_{1,0}}(z_{0,0})$. For the first case, we have
$$d_{D_1}(z)\leq d_{D_1}(z_{0,0})+|z_{0,0}-z|\leq
\frac{3}{2}d_{D_1}(z_{0,0}),
$$
$$d_{D_{1,0}}(z)\geq d_{D_{1,0}}(z_{0,0})-|z_{0,0}-z|\geq
\frac{1}{2}d_{D_{1,0}}(z_{0,0})
$$
and
$$d_{D}(z)\leq d_{D}(z_{0,0})+|z_{0,0}-z|\leq \frac{3}{2}d_D(z_{0,0}).
$$
Since Lemma \ref{cl8-h-1} implies that $z_{0,0}\in \gamma_{1,0}[z_1, x_{1, k_1}]$, we see from Lemma
\ref{lem4-h-0'} that
\be\label{1}d_{D_{1,0}}(z)\geq
\frac{1}{2}d_{D_{1,0}}(z_{0,0})\geq
\frac{1}{2^7(3\rho_1+1)\rho_2}d_{D_1}(z_{0,0})\geq
\frac{1}{2^6(9\rho_1+3)\rho_2}d_{D_1}(z),
\ee
and Lemma \ref{co1} leads to
\beq\label{11-11}
d_{D_{1,0}}(z)&\geq& \frac{1}{2^7(3\rho_1+1)\rho_2}d_{D_1}(z_{0,0})\\
\nonumber &> & \frac{1}{2^7(3\rho_1+1)\rho_1\rho_2\rho_5}d_D(z_{0,0})\\
\nonumber&\geq& \frac{1}{2^6\rho_1(9\rho_1+3)\rho_2\rho_5}d_D(z).
\eeq

For the second case, that's, $|z_{0,0}-z|>
\frac{1}{2}d_{D_{1,0}}(z_{0,0})$, again, by Lemmas \ref{lem4-h-0'}, \ref{cl8-h-1} and \ref{co1}, we see that
\begin{eqnarray*}
d_{D_1}(z)\leq d_{D_1}(z_{0,0})+|z_{0,0}-z|&\leq& 2^6(3\rho_1+1)\rho_2d_{D_{1,0}}(z_{0,0})+|z_{0,0}-z|
\\
\nonumber&\leq& \big(2^7(3\rho_1+1)\rho_2+1\big)|z_{0,0}-z|
\end{eqnarray*}
and
\begin{eqnarray*}
d_D(z)\leq d_D(z_{0,0})+|z_{0,0}-z|&\leq&
2^6\rho_1(3\rho_1+1)\rho_2\rho_5d_{D_{1,0}}(z_{0,0})+|z_{0,0}-z|
\\
\nonumber&\leq& \big(2^7\rho_1(3\rho_1+1)\rho_2\rho_5+1\big)|z_{0,0}-z|,
\end{eqnarray*}
which, together with (\ref{mlm-3-3}), implies
\beq\label{2}
d_{D_{1,0}}(z)\geq \frac{1}{\mu_8}|z_{0,0}-z|\geq
\frac{1}{\mu_8\big(2^7(3\rho_1+1)\rho_2+1\big)}d_{D_1}(z)
\eeq
and
\be\label{2-h-1}
d_{D_{1,0}}(z)\geq \frac{1}{\mu_8}|z_{0,0}-z|\geq
\frac{1}{\mu_8\big(2^7\rho_1(3\rho_1+1)\rho_2\rho_5+1\big)}d_D(z).
\ee

We conclude from \eqref{13-1}, \eqref{1}, (\ref{11-11}), \eqref{2} and (\ref{2-h-1}) that the
inequalities \eqref{q-1-0} and \eqref{q-1-1} in Lemma \ref{lem5-h-0'} are true.

It follows directly  from  (\ref{mlm-3-3}), Lemmas \ref{lem-4-1''} and \ref{cl6-h-3''} that
$$d_{D_{1,0}}(z)\geq \frac{1}{2\mu_8}\min\{d_{D_{1,0}}(z_{0,0}), d_{D_{1,0}}(z_{1,0})\}
\geq \frac{1}{2\mu_8}\;d_{D_{1,0}}(z_{0,0}).
$$
Hence, by Lemmas \ref{lem4-h-0'} and \ref{cl8-h-1}, we have
$$d_{D_1}(z)\geq (\rho_1-1)d_{D_{1,0}}(z)
\geq \frac{\rho_1-1}{2\mu_8}d_{D_{1,0}}(z_{0,0})>\frac{1}{2^9\mu_8\rho_2}d_{D_1}(z_{0,0}),
$$
which shows that Lemma \ref{lem5-h-0'}\eqref{q-0} holds.

Based on the inequalities \eqref{q-1-0} and \eqref{q-1-1} in Lemma \ref{lem5-h-0'},
the proofs of \eqref{q-2-0} and \eqref{q-2-1} in Lemma  \ref{lem5-h-0'}
easily follow from a similar argument as in that of Lemma \ref{lem4-h-0'}(2).
Hence the proof of Lemma \ref{lem5-h-0'} is complete.

\begin{lem}\label{cl5-h-w-1} For $w'_1,w'_2\in\gamma'_{2,0}$, if
$k_{D'_{1,0}}(w'_1, w'_2)\geq \mu_2^2\rho_7$, then
$$k_{D'}(w'_1, w'_2)\geq \frac{1}{\mu_2^2\rho_7}k_{D'_{1,0}}(w'_1, w'_2).
$$
\end{lem}

We shall apply Theorem \Ref{ThmF} to obtain a proof of Lemma \ref{cl5-h-w-1}. Indeed,
since the assumption $k_{D'_{1,0}}(w'_1, w'_2)\geq \mu_2^2\rho_7$ implies $k_{D_{1,0}}(w_1, w_2)\geq 1$, by
Theorem \Ref{ThmF}, we deduce that
\be\label{eq1.2w-1}
k_{D_{1,0}}(w_1, w_2)\geq \frac{1}{\mu_2}k_{D'_{1,0}}(w'_1, w'_2)\geq \mu_2\rho_7.
\ee
Further, we obtain from Corollary \ref{Cor-1} and Lemma \ref{lem5-h-0'} that
\beq\label{d-1}
k_{D_{1,0}}(w_1, w_2)&\leq& \mu_8\log \left (1+\frac{|w_1-w_2|}{\min\{d_{D_{1,0}}(w_1), d_{D_{1,0}}(w_2)\}}\right )\\ \nonumber
&\leq&\mu_8\log \left (1+2^5\rho_1\rho_2\rho_6\frac{|w_1-w_2|}{\min\{d_{D}(w_1),
d_{D}(w_2)\}}\right)\\
\nonumber  &<& \rho_7\log \left (1+\frac{|w_1-w_2|}{\min\{d_{D}(w_1), d_{D}(w_2)\}}\right )\\
\nonumber  &\leq& \rho_7k_{D}(w_1,w_2),
\eeq
whence the combination with (\ref{eq1.2w-1}) shows $k_{D}(w_1,w_2)\geq\mu_2$, which implies
$k_{D'}(w'_1, w'_2)\geq 1$. Therefore, (\ref{eq1.2w-1}), \eqref{d-1}
and Theorem \Ref{ThmF} yield
$$k_{D'}(w'_1, w'_2)\geq \frac{1}{\mu_2}k_{D}(w_1, w_2)
\geq \frac{1}{\mu_2^2\rho_7} k_{D'_{1,0}}(w'_1, w'_2)
$$
as required.
\medskip

Our next lemma is as follows.

\begin{lem}\label{cl8-0''} Suppose $s_1= [\rho_7] $ and $s_2=\mu_2^2\rho_7$.
For $w'_1$, $w'_2\in\gamma'_{2,0}$,  if $k_{D'_{1,0}}(w'_1,w'_2)\geq
\mu_2^2\rho_7^2$ and $k_{D'}(w'_1,w'_2)\geq s_0>0$, then
\bee
\item[$(1)$] $\gamma'_{2,0}[w'_1,w'_2]$ satisfies $QH\,(s_0,s_1,s_2)$-condition in $D'$;

\item[$(2)$] $k_{D'_{1,0}}(w'_1,w'_2)\leq
48a^2\mu_2^6\rho_7^4\log\Big(1+\frac{\ell(\gamma'_{2,0}[w'_1,w'_2])}{d_{D'}(y'_0)}\Big)$
provided $s_0\geq 8s_1s_2^2$, where $d_{D'}(y'_0)=\inf_{z'\in\gamma'_{2,0}[w'_1,w'_2]}d_{D'}(z')$.
\eee
\end{lem}
\medskip

To prove Lemma \ref{cl8-0''}(1), since $w'_1\not=w'_2$, we only need to check Conditions
$(2)$ and $(3)$ in the definition of $QH$ condition.

First, we partition the part $\gamma'_{2,0}[w'_1,w'_2]$ of $\gamma'_{2,0}$. Let $x'_0=w'_1$, and let $x'_1, \ldots, x'_{[\rho_7]}\in
\gamma'_{2,0}[w'_1,w'_2]$ be successive points such that for each
$i\in \{1,\ldots,[\rho_7]\}$,
$$k_{D'_{1,0}}(x'_{i-1},x'_{i})=\frac{1}{[\rho_7]}\, k_{D'_{1,0}}(w'_1,w'_2),
$$
and thus $k_{D'_{1,0}}(x'_{i-1},x'_{i})
\geq\mu_2^2\rho_7.$ Then we see from Lemma \ref{cl5-h-w-1} that for all $i\not=j\in\{1,\ldots,[\rho_7]\}$,
$$k_{D'}(x'_{i-1},x'_{i})\leq k_{D'_{1,0}}(x'_{i-1},x'_{i})=
k_{D'_{1,0}}(x'_{j-1},x'_{j})\leq \mu_2^2\rho_7k_{D'}(x'_{j-1},x'_{j}).
$$
Similarly, for each $x'\in\gamma'_{2,0}[x'_{i},w'_2]$, again, we get from Lemma \ref{cl5-h-w-1}
that
\beq \nonumber k_{D'}(x'_{i-1},x'_{i})&\leq&
\ell_{k_{D'_{1,0}}}(\gamma'_{2,0}[x'_{i-1},x'_{i}])\leq
\ell_{k_{D'_{1,0}}}(\gamma'_{2,0}[x'_{i-1},x']) =
k_{D'_{1,0}}(x'_{i-1},x')\\ \nonumber &\leq&
\mu_2^2\rho_7k_{D'}(x'_{i-1},x').
\eeq
We complete the proof of $(1)$ in Lemma \ref{cl8-0''}.

Since $s_0\geq 8s_1s_2^2$, again, it follows from Lemma \ref{cl5-h-w-1} together with Theorem \ref{scl-1-1} and $(1)$
in Lemma \ref{cl8-0''} that
\begin{eqnarray*} \hspace{1cm}
k_{D'_{1,0}}(w'_1,w'_2) &\leq& \mu_2^2\rho_7k_{D'}(w'_1,w'_2)\leq
48a^2\mu_2^6\rho_7^4\log\left(1+
\frac{\diam(\gamma'_{2,0}[w'_1,w'_2])}{d_{D'}(y'_0)}\right ).
\end{eqnarray*}
The proof of Lemma \ref{cl8-0''} is complete.
\medskip

\begin{lem}\label{14-2}
\begin{enumerate}

\item\label{hq-2} For each $z'\in\gamma'_{2,0}$,
$$\diam(\gamma'_{1,0}[z'_{0,0},z'_{1,0}])\leq\frac{5}{4}\rho_{10}e^{2^7\mu_2\mu_8(c\rho_1\rho_5)^2}d_{D'_1}(z');
$$

\item\label{hq-1} $\diam(\gamma'_{2,0})\leq
\frac{8}{5}\Big(6+\frac{1}{\rho_{10}}\Big)\diam(\gamma'_{1,0}[z'_{0,0},z'_{1,0}])$;

\item\label{hq-3}
$k_{D'_{1,0}}(z'_{1,0}, z'_{0,0})\leq \rho_8\rho_9.$
\end{enumerate}
\end{lem}

For a proof of Lemma \ref{14-2}, we first need some preparation.
%the similar
%discussions as in the proofs of (\ref{1}) and \eqref{2}, together
%with Lemma \ref{lem4-h-0'}, show that \be\label{mei-1}d_{D_1}(z)\geq
%(\rho_1-1)d_{D_{1,0}}(z)\geq
%\frac{1}{2^{11}\mu_8\rho_2}d_{D_1}(z_{0,0}).\ee

From Lemma \ref{cl6-h-3''}, we know that
\be\label{mei-2}
\max\{|z_{0,0}-z_{1,0}|,\min\{d_{D_{1,0}}(z_{0,0}),d_{D_{1,0}}(z_{1,0})\}\}=|z_{0,0}-z_{1,0}|.
\ee
Then Corollary \ref{Cor-1} yields
$$\diam(\gamma_{2,0})\leq \mu_8|z_{0,0}-z_{1,0}|,
$$
and so by Lemma \ref{lem4-h-1}, we see that for each $z\in\gamma_{2,0}$,
\be\label{mei-3}
d_{D_1}(z)\leq \diam(\gamma_{2,0})+d_{D_1}(z_{1,0})
\leq(\mu_8\rho_5+1)d_{D_1}(z_{1,0}).
\ee
Since $z\in \gamma_{2,0}$, it follows from Lemma \ref{lem2.1} that there must exist some
$t_0\in\{1,\ldots,k_1\}$ such that $z\in B_{{1,t_0}}=\mathbb{B}(x_{1,t_0}, r_{1,t_0})$.
Then Lemma \ref{mxll-1''} shows that
\be\label{mei-5}
|z'-x'_{1,t_0}|< \frac{1}{5}d_{D'_1}(x'_{1,t_0})\;\;\mbox{and}\;\;d_{D'_1}(x'_{1,t_0})<
\frac{5}{4}d_{D'_1}(z').
\ee

Now, it is possible to prove Lemma \ref{14-2}. For a proof of \eqref{hq-2} and \eqref{hq-1}
in the lemma, we need upper bounds for the quotients
$$\frac{|z'_{0,0}-z'|}{\diam(\gamma'_{1,0}[z'_{0,0},z'_{1,0}])} ~\mbox{ and }~
\frac{\diam(\gamma'_{1,0}[z'_{0,0},z'_{1,0}])}{d_{D'_1}(z')}.
$$
To get these upper bounds, we divide the proof into three cases based on the
positions of $x_{1, t_0}$ in $\gamma_{1,0}$.

\begin{figure}[!ht]
\begin{center}
\includegraphics%[width=10cm]
{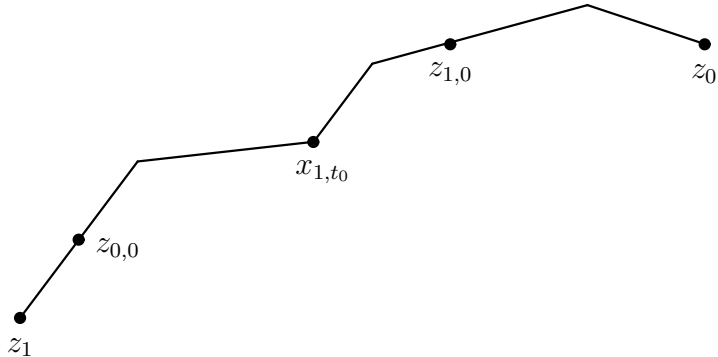}
\end{center}
\caption{The case $x_{1,t_0}\in \gamma_{1,0}[z_{0,0}, z_{1,0}]$\label{fig11}}
\end{figure}

For the first case when $x_{1,t_0}\in\gamma_{1,0}[z_{0,0},z_{1,0}]$ (see Figure \ref{fig11}),
we know from (\ref{eq(h-3-3)}), (\ref{24-1-mei}) and (\ref{mei-5}) that
\beq\label{meir-1}
|z'_{0,0}-z'| &\leq & |z'-x'_{1,t_0}|+|z'_{0,0}-x'_{1,t_0}|
\\ \nonumber &<&
\frac{1}{5}d_{D'_1}(x'_{1,t_0})
+|z'_{0,0}-x'_{1,t_0}|
\\ \nonumber &\leq & \frac{1}{5}\big(|x'_{1,t_0}-z'_{1,0}|+d_{D'_1}(z'_{1,0})\big)+\diam(\gamma'_{1,0}[z'_{0,0},
z'_{1,0}])
\\ \nonumber &\leq & \frac{6}{5}\diam(\gamma'_{1,0}[z'_{0,0},
z'_{1,0}])+\frac{1}{5\rho_{10}}\diam(\gamma'_{1,0}[z'_{1},
z'_{1,0}])
\\ \nonumber &= &
\frac{6}{5}\diam(\gamma'_{1,0}[z'_{0,0},
z'_{1,0}])+\frac{4}{5\rho_{10}}|z'_{1,0}- z'_{0,0}|
\\ \nonumber
&\leq & \frac{2}{5}\Big(3+\frac{2}{\rho_{10}}\Big)
\diam(\gamma'_{1,0}[z'_{0,0},z'_{1,0}]),
\eeq
and the choice of $z'_{1,0}$, (\ref{24-1}) and (\ref{mei-5}) lead to
\beq\label{meir-2}
d_{D'_1}(z')&> &
\frac{4}{5}d_{D'_1}(x'_{1,t_0})\geq
\frac{4}{5\rho_{10}}\diam(\gamma'_{1,0}[z'_1, x'_{1,t_0}])\\
\nonumber &\geq & \frac{4}{5\rho_{10}}\diam(\gamma'_{1,0}[z'_1,
z'_{0,0}])
\geq \frac{2}{5\rho_{10}}\diam(\gamma'_{1,0}[z'_1, z'_{1,0}])\\
\nonumber &\geq&
\frac{4}{5\rho_{10}}\diam(\gamma'_{1,0}[z'_{0,0},z'_{1,0}]).
\eeq

%\begin{figure}[!ht]
%\includegraphics[width=0.55\textwidth]{figure14} %,height=0.4\textwidth
%\caption{The case $x_{1,t_0}\in \gamma_{1,0}[z_{0,0}, z_{1,0}]$\label{fig11}}
%\end{figure}

For the second case when $x_{1,t_0}\in\gamma_{1,0}[z_1,z_{0,0}]$, it
follows from Lemmas \ref{lem4-h-1}, \ref{lem-4-1'} and
\ref{lem5-h-0'} that
\begin{eqnarray*}k_{D_1}(z_{0,0},x_{1,t_0})&\leq&
\int_{\gamma_{1,0}[x_{1,t_0},z_{0,0}]}\frac{|dx|}{d_{D_1}(x)} \leq
2\rho_5\log\Big(1+\frac{2\ell(\gamma_{1,0}[x_{1,t_0},z_{0,0}])}{d_{D_1}(x_{1,t_0})}\Big)\\
\nonumber&\leq&
2\rho_5\log\Big(1+\frac{2\rho_5d_{D_1}(z_{0,0})}{d_{D_1}(x_{1,t_0})}\Big)\leq
2\rho_5\log\Big(1+2^{10}\mu_8\rho_2\rho_5\frac{d_{D_1}(z)}{d_{D_1}(x_{1,t_0})}\Big)
\\ \nonumber&\leq &
2\rho_5\log\Big(1+\frac{2^{10}\mu_8\rho_2\rho_5(\rho_1+1)}{\rho_1}\Big)
\\ \nonumber&<&
2^7\mu_8(c\rho_1\rho_5)^2,
\end{eqnarray*}
since $d_{D_1}(z)\leq (1+\frac{1}{\rho_1})d_{D_1}(x_{1,t_0})$. Hence we see from Theorem
\Ref{ThmF} that
\begin{eqnarray*}
\max\bigg\{\log\Big(1+\frac{|z'_{0,0}-x'_{1,t_0}|}{d_{D'_1}(z'_{0,0})}\Big),
\Big|\log\frac{d_{D'_1}(x'_{1,t_0})}{d_{D'_1}(z'_{0,0})}\Big|\bigg\}
&\leq& k_{D'_1}(z'_{0,0},x'_{1,t_0})<
2^7\mu_2\mu_8(c\rho_1\rho_5)^2,
\end{eqnarray*}
which implies
\be\label{m-1}
\max\{|z'_{0,0}-x'_{1,t_0}|,d_{D'_1}(x'_{1,t_0})\}\leq
e^{2^7\mu_2\mu_8(c\rho_1\rho_5)^2}d_{D'_1}(z'_{0,0})
\ee
and
\be\label{m-1-1}
d_{D'_1}(z'_{0,0})\leq e^{2^7\mu_2\mu_8(c\rho_1\rho_5)^2}d_{D'_1}(x'_{1,t_0}).
\ee
Then, by (\ref{eq(h-3-3)}), (\ref{24-1-mei}) and (\ref{mei-5}), we have
\beq\label{meir-3}
|z'_{0,0}-z'|&\leq& |z'-x'_{1,t_0}|+|z'_{0,0}-x'_{1,t_0}|
\\ \nonumber&<&
\frac{1}{5}d_{D'_1}(x'_{1,t_0})+\diam(\gamma'_{1,0}[z'_1,z'_{1,0}])\\
\nonumber&\leq&
\frac{1}{5}(|x'_{1,t_0}-z'_{1,0}|+d_{D'_1}(z'_{1,0}))+\diam(\gamma'_{1,0}[z'_1,z'_{1,0}])\\
\nonumber&\leq &
\frac{1}{5}\Big(6+\frac{1}{\rho_{10}}\Big)\diam(\gamma'_{1,0}[z'_{1},z'_{1,0}])
\\
\nonumber&\leq&
\frac{4}{5}\Big(6+\frac{1}{\rho_{10}}\Big)\diam(\gamma'_{1,0}[z'_{0,0},z'_{1,0}]),
\eeq
and (\ref{24-1}), (\ref{mei-5}) and the choice of $z'_{1,0}$ lead to
\beq\label{meir-4}
d_{D'_1}(z') &\geq& \frac{4}{5}d_{D'_1}(x'_{1,t_0})\geq
\frac{4}{5e^{2^7\mu_2\mu_8(c\rho_1\rho_5)^2}}d_{D'_1}(z'_{0,0})
\\ \nonumber &\geq &
\frac{4}{5\rho_{10}e^{2^7\mu_2\mu_8(c\rho_1\rho_5)^2}}\diam(\gamma'_{1,0}[z'_{1},z'_{0,0}])
\\ \nonumber &\geq &
\frac{2}{5\rho_{10}e^{2^7\mu_2\mu_8(c\rho_1\rho_5)^2}}\diam(\gamma'_{1,0}[z'_{1},z'_{1,0}])
\\ \nonumber &\geq &
\frac{4}{5\rho_{10}e^{2^7\mu_2\mu_8(c\rho_1\rho_5)^2}}\diam(\gamma'_{1,0}[z'_{0,0},z'_{1,0}]).
\eeq

For the last case, that's, $x_{1,t_0}\in\gamma_{1,0}[z_{1,0},z_0]$, we
see from Lemmas \ref{lem4-h-1} and \ref{lem-4-1'}, together with
(\ref{mei-3}), that
\begin{eqnarray*}
k_{D_1}(z_{1,0},x_{1,t_0})&\leq& \int_{\gamma_{1,0}[z_{1,0},x_{1,t_0}]}\frac{|dx|}{d_{D_1}(x)}
\leq
2\rho_5\log\Big(1+\frac{2\ell(\gamma_{1,0}[z_{1,0},x_{1,t_0}])}{d_{D_1}(z_{1,0})}\Big)\\
\nonumber&\leq&
2\rho_5\log\Big(1+\frac{2\rho_5d_{D_1}(x_{1,t_0})}{d_{D_1}(z_{1,0})}\Big)
\leq
2\rho_5\log\Big(1+\frac{2\rho_1\rho_5}{\rho_1-1}\cdot\frac{d_{D_1}(z)}{d_{D_1}(z_{1,0})}\Big)
\\ \nonumber&\leq &
2\rho_5\log\Big(1+\frac{2\rho_1\rho_5(\mu_8\rho_5+1)}{\rho_1-1}\Big)
\\ \nonumber&<&
5c\mu_8\rho_5,
\end{eqnarray*}
since $d_{D_1}(x_{1,t_0})\leq \frac{\rho_1}{\rho_1-1}d_{D_1}(z)$ for $z\in B_{1,t_0}$. By
replacing $2^7\mu_8(c\rho_1\rho_5)^2$ with $5c\mu_8\rho_5$, a
similar reasoning as in the proofs of (\ref{m-1}) and (\ref{m-1-1})
shows that
$$\max\{|z'_{1,0}-x'_{1,t_0}|,d_{D'_1}(x'_{1,t_0})\}\leq
e^{5c\mu_2\mu_8\rho_5}d_{D'_1}(z'_{1,0})
$$
and
$$ d_{D'_1}(z'_{1,0})\leq e^{5c\mu_2\mu_8\rho_5}
d_{D'_1}(x'_{1,t_0}).
$$
Hence, by (\ref{eq(h-3-3)}), (\ref{24-1-mei}) and (\ref{mei-5}), we have
\beq\label{meir-6}
|z'_{0,0}-z'|&\leq&
(|z'_{0,0}-z'_{1,0}|+|z'_{1,0}-x'_{1,t_0}|+|x'_{1,t_0}-z'|)\\
\nonumber&\leq&
\diam(\gamma'_{1,0}[z'_{0,0},z'_{1,0}])+\frac{6}{5}e^{5c\mu_2\mu_8\rho_5}d_{D'_1}(z'_{1,0})\\
\nonumber&\leq&
\Big(1+\frac{24}{5\rho_{10}}e^{5c\mu_2\mu_8\rho_5}\Big)\diam(\gamma'_{1,0}[z'_{0,0},z'_{1,0}]),
\eeq
and (\ref{eq(h-3-3)}), (\ref{24-1}) and (\ref{mei-5}) lead to
\beq\label{meir-7}
d_{D'_1}(z')&>& \frac{4}{5}d_{D'_1}(x'_{1,t_0})\geq
\frac{4}{5e^{5c\mu_2\mu_8\rho_5}}d_{D'_1}(z'_{1,0})\\
\nonumber&\geq&
\frac{8}{5\rho_{10}e^{5c\mu_2\mu_8\rho_5}}\diam(\gamma'_{1,0}[z'_{0,0},z'_{1,0}]).
\eeq
The inequalities \eqref{meir-1}, (\ref{meir-2}), \eqref{meir-3}, (\ref{meir-4}),  \eqref{meir-6} and (\ref{meir-7})
are our requirements. Also, (\ref{meir-2}), (\ref{meir-4}) and (\ref{meir-7}) show
that Lemma \ref{14-2}(\ref{hq-2}) is true.

Further, we get from (\ref{meir-1}),
(\ref{meir-3}) and (\ref{meir-6}) that for $z'\in\gamma'_{2,0}$,
\be\label{meiren-1}
|z'_{0,0}-z'|\leq \frac{4}{5}\Big(6+\frac{1}{\rho_{10}}\Big)\diam(\gamma'_{1,0}[z'_{0,0},z'_{1,0}]).
\ee
Let $x'_1\in\gamma'_{2,0}$ be such that
$$|z'_{0,0}-x'_1|= \frac{1}{2}\diam(\gamma'_{2,0}).
$$
Then we know from (\ref{meiren-1}) that
$$\diam(\gamma'_{2,0})= 2|z'_{0,0}-x'_1|\leq \frac{8}{5} \Big(6+\frac{1}{\rho_{10}}\Big)\diam(\gamma'_{1,0}[z'_{0,0},z'_{1,0}]),
$$
whence Lemma \ref{14-2}(\ref{hq-1}) is also true.

For a proof of $(\ref{hq-3})$ in Lemma \ref{14-2}, we let
$s_1=[\rho_7]$, $s_2=\mu_2^2\rho_7$ and $s_0=8s_1s_2^2$. Since  Lemma \ref{cl6-h-3''} implies
$$k_{D'_{1,0}}(z'_{0,0}, z'_{1,0})\geq k_{D'_1}(z'_{0,0}, z'_{1,0})\geq \log\frac{\rho_{10}}{4}>\mu_2^2\rho_7^2,
$$
we see from Lemma \ref{cl5-h-w-1} that
$$k_{D'}(z'_{1,0}, z'_{0,0})\geq \frac{1}{\mu_2^2\rho_7}k_{D'_{1,0}}(z'_{1,0}, z'_{0,0})
\geq \frac{1}{\mu_2^2\rho_7}\log\frac{\rho_{10}}{4}>s_0.
$$
Hence Lemma \ref{cl8-0''} shows that $\gamma'_{2,0}$ satisfies $(s_0,s_1,
s_2)$-$QH$ condition. It follows from Lemma \ref{cl8-0''} together with $(\ref{hq-2})$ and
$(\ref{hq-1})$ in this Lemma that
\begin{eqnarray*} \hspace{1cm}
k_{D'_{1,0}}(z'_{1,0}, z'_{0,0}) &\leq&
48a^2\mu_2^6\rho_7^4\log\left(1+
\frac{\diam(\gamma'_{2,0})}{d_{D'}(w'_0)}\right )\\
\nonumber &\leq&48a^2\mu_2^6\rho_7^4\log\left(1+
\frac{\frac{8}{5}(6+\frac{1}{\rho_{10}})\diam(\gamma'_{1,0}[z'_{0,0},z'_{1,0}])}{d_{D'}(w'_0)}\right
)\\
\nonumber &<& 48a^2\mu_2^6\rho_7^4\log\left(1+
\frac{13\rho_{10}e^{2^7\mu_2\mu_8(c\rho_1\rho_5)^2}d_{D'_1}(w'_0)}{d_{D'}(w'_0)}\right
)\\ \nonumber &<&\rho_8\rho_9,
\end{eqnarray*}
where $w'_0\in \gamma'_{2,0}$ satisfies
$$d_{D'}(w'_0)=\inf\{d_{D'}(z'):\; z'\in \gamma'_{2,0}\}.
$$
Hence the inequality Lemma \ref{14-2}(\ref{hq-3}) holds, and thus
the proof of Lemma \ref{14-2} is complete.

%%%%%%%%%%%%%%%%%%%%%%%%%%%%%%%%%%%
%%%%%%%%%%%%%%%%%%%%%%%%%%%%%%%%%%%
\subsubsection{{\bf Lemmas: Part II}}\ \ \
%%%%%%%%%%%%%%%%%%%%%%%%%%%%%%%%%%%
%%%%%%%%%%%%%%%%%%%%%%%%%%%%%%%%%%%
\medskip

{\bf The sketch for this part}:\quad In this part, first, we choose a proper point
$\omega'_0$ in $\gamma'_{2,0}$ and partition
the part $\gamma'_{2,0}[\omega'_0, z'_{1,0}]$ of $\gamma'_{2,0}$ with the aid
of some special points $\{y'_i\}_{i=1}^{m+1}$. After that, we shall pick
up a needed point $\{y'_{i_1}\}$ from  $\{y'_i\}_{i=1}^{m+1}$
and another special point $v'_3$ from $\gamma'_{2,0}$, and then
partition the part $\gamma'_{2,0}[y'_{i_1}, v'_3]$ of $\gamma'_{2,0}$ by using
the points $\{x'_i\}_{i=1}^{4[\Gamma(\Lambda_{T_1+3})]+1}$. Based on the obtained points, we shall construct the corresponding
carrot arcs, balls etc. With the aid of the related points, carrot arcs,
balls etc, seven Lemmas will be proved.
\medskip

We begin to let $\omega'_0$ to be the first point in $\gamma'_{2,0}$ along the
direction from $z'_{1,0}$ to $z'_{0,0}$ such that
\beq \label{m-w-l-s-1}
d_{D'_{1,0}}(\omega'_0)=\sup\limits_{p'\in \gamma'_{2,0}}d_{D'_{1,0}}(p').
\eeq
It is possible that $\omega'_0=z'_{1,0}$ or $z'_{0,0}$. Clearly,
there exists a nonnegative integer $m$ such that
\be\label{india-mw-1}
2^{m}\, d_{D'_{1,0}}(z'_{1,0}) \leq d_{D'_{1,0}}(\omega'_0)< 2^{m+1}\, d_{D'_{1,0}}(z'_{1,0}).
\ee

\begin{figure}[!ht]
\begin{center}
\includegraphics%[width=10cm]
{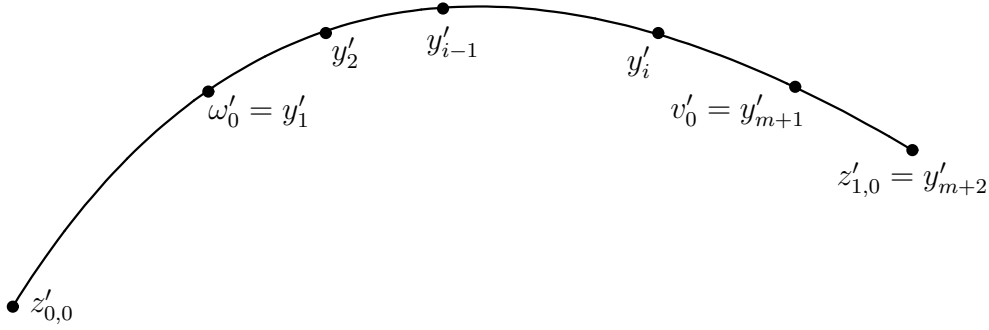}
\end{center}
\caption{The arc $\gamma'_{2,0}$ and the related point $\omega_0'$
\label{fig12}}
\end{figure}

Let $v'_0$ be the first point in $\gamma'_{2,0}[z'_{1,0},\omega'_0]$
from $\omega'_0$ to $z'_{1,0}$ satisfying (see Figure \ref{fig12})
$$d_{D'_{1,0}}(\omega'_0)=2^{m}\, d_{D'_{1,0}}(v'_0).
%\;\;\mbox{(see Figure \ref{fig12})}.
$$

Now, we give a partition to $\gamma'_{2,0}[\omega'_0,v'_0]$.
Let $y_1'=\omega'_0$. If $v'_0= y'_1$, we let $y'_2=v'_0$. If $v'_0\not= y'_1$, then we let
$y'_2,\ldots ,y'_{m+1}\in \gamma'_{2,0}[\omega'_0,v'_0]$ be the
points such that for each $i\in \{2,\ldots,m+1\}$, $y'_i$ denotes
the first point from $\omega'_0$ to $v'_0$ with
%\be\label{hws-eq(4.3)}
\beq \label{m-w-l-s-2}
d_{D'_{1,0}}(y'_i)=\frac{1}{2^{i-1}}\, d_{D'_{1,0}}(y'_1).
\eeq
Obviously, $y'_{m+1}=v'_0$. If $v'_0\not= z'_{1,0}$, then we use
$y'_{m+2}$ to denote $z'_{1,0}$ (see Figure \ref{fig12}).

%{\color{red}}

%\begin{figure}[!ht]
%\includegraphics[width=0.75\textwidth]{figure16} %,height=0.4\textwidth
%\caption{The arc $\gamma'_{2,0}$ and the related points\label{fig12}}
%\end{figure}

\begin{lem}\label{India-samy}
Suppose that $w'_1\in \gamma'_{2,0}[y'_{m+1},\omega'_0]$ and
$w'_2\in\gamma'_{2,0}[w'_1,\omega'_0]$ which satisfy
$$d_{D'_{1,0}}(w'_2)\geq \Gamma(t)\;d_{D'_{1,0}}(w'_1)
$$
for some $t>(\mu_2\rho_7)^2$. Then there exists some
$w'\in\gamma'_{2,0}[w'_1,w'_2]$ such that
$$\ell(\gamma'_{2,0}[w'_1,w'])\geq t\;d_{D'_{1,0}}(w'),
$$
where $\Gamma(t)$ is defined as in Section \ref{sec-6}.\end{lem}

We prove Lemma \ref{India-samy} by a method of contradiction. Suppose that  for each
$w'\in\gamma'_{2,0}[w'_1,w'_2]$, we have
\beq \label{india-s-1}
\ell(\gamma'_{2,0}[w'_1,w'])< t\;d_{D'_{1,0}}(w').
\eeq

\begin{figure}[!ht]
\includegraphics[width=0.85\textwidth]{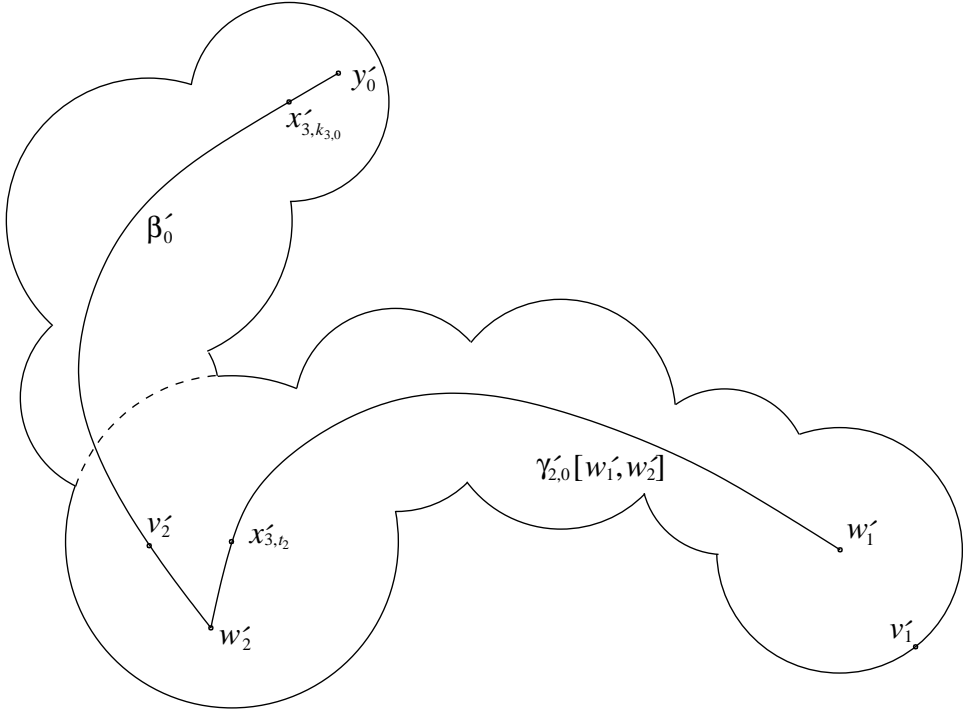} %,height=0.4\textwidth
\caption{The arc $\beta'_{2,0}$ and the domain $D'_{3,0}$ together
with its subdomain $D'_{3,1}$
%{\color{red}$\gamma_{2,0}[w'_1, w'_2]$ should be $\gamma'_{2,0}[w'_1, w'_2]$}
\label{fig13}}
\end{figure}

In the following, based on a new carrot arc $\beta'_{2,0}$, we shall construct a uniform domain $D'_{3,0}$ in $D'$
together with its uniform subdomain  $D'_{3,1}$ and then show that the restriction $f^{-1}|_{D'_{3,1}}$ is weakly
$H(t)$-quasisymmetric in $\delta_{D'_{3,0}}$ and $\delta_{D_{3,0}}$
(see Subsection \ref{subsec2.5} for the definition of internal metrics).
This weak quasisymmetry of $f^{-1}$ will help us to get a contradiction.

At first, let us construct a new arc in $D'$. Since $D'$ is a bounded $a$-John domain with center $y'_0$,
there must exist an $a$-carrot arc $\beta'_0$ in $D'$ joining
$w'_2$ and $y'_0$ with center $y'_0$.
We set (see Figure \ref{fig13})
$$\beta'_{2,0}=\gamma'_{2,0}[w'_1,w'_2]\cup\beta'_0.
%\; \mbox{ (see Figure \ref{fig13})}.
$$
If $y'_0=w'_2$, we assume that $\beta'_0=\{w'_2\}$.

By using \eqref{india-s-1}, we shall prove that $\beta'_{2,0}$ is a carrot arc.
For $z'\in \beta'_{2,0}$, if $z'\in\gamma'_{2,0}[w'_1,w'_2]$, then \eqref{india-s-1}
implies that
$$\ell(\gamma'_{2,0}[w'_1,z'])< t\;d_{D'_{1,0}}(z');
$$
On the other hand, if $z'\in\beta'_0$, since $|w'_2-z'|\leq
a\,d_{D'}(z')$,  Lemma \ref{lem-4-1''} yields
$$d_{D'}(w'_2)\leq 2a d_{D'}(z'),
$$
and so \eqref{india-s-1} leads to
\begin{eqnarray*}
\ell(\beta'_{2,0}[w'_1,z'])&=&\ell(\gamma'_{2,0}[w'_1,w'_2])+\ell(\beta'_0[w'_2,z'])\leq
t\;d_{D'_{1,0}}(w'_2)+a\,d_{D'}(z')\\
\nonumber &\leq&  2atd_{D'}(z')+a\,d_{D'_{1,0}}(z')\\ \nonumber &<&
3at\,d_{D'_{1,0}}(z').
\end{eqnarray*}
Hence we have proved that $\beta'_{2,0}$ is a $3at$-carrot arc from
$w'_1$ to $y'_0$ with center $y'_0$, and thus Lemmas \ref{lem-m-v-1} and
\ref{lem2.2-4-mv-3} guarantee the following.

\begin{prop}\label{00-11}
There exists a simply connected $2^{47+(3a\rho_7t)^2}a^2\rho_7t^2$-uniform
domain $D'_{3,0}=\bigcup \limits_{i=1}^{k_{3,0}}B_{3,i}\subset D'$ such that
\bee
\item\label{equ--m-1}
$w'_1$, $y'_{0}\in D'_{3,0}$;
\item\label{equ--m-2}
For each $i\in \{1,\ldots, k_{3,0}\}$,
$$\frac{1}{3\rho_7}\,d_{D'}(x'_{3,i})\leq r_{3,i}\leq
\frac{1}{\rho_7}d_{D'}(x'_{3,i}),
$$
\eee
where $B_{3,i}=\mathbb{B}(x'_{3,i}, r_{3,i})$, $x'_{3,i}\in
\beta'_{2,0}$, $x'_{3,i}\not\in B_{3,i-1}$ for each $i\in \{2,
\ldots, k_{3,0}\}$ and $x'_{3,1}=w'_1$  $($see Figure \ref{fig13}$)$.
\end{prop}

Next, we construct a subdomain of $D'_{3,1}$. Let
$$t_2=\max\{t: x'_{3,i}\in \gamma'_{2,0}\;\;\mbox{for each }\;\;i\in\{1,\ldots,t\}\}.
$$
Then we take (see Figure \ref{fig13})
$$D'_{3,1}=\bigcup \limits_{i=1}^{t_2}B_{3,i}.
%\;\; \mbox{ (see Figure \ref{fig13})}.
$$

We easily know that $D'_{3,1}$
is a $2^{47+(3a\rho_7t)^2}a^2\rho_7t^2$-uniform domain.
Since for each $i\in\{1,\ldots,t_2\}$, Lemma \ref{lem5-h-0'} and
Proposition \ref{00-11} imply
$$r_{3,i}\leq \frac{d_{D'}(x'_{3,i})}{\rho_7}\leq
\frac{(2^6\rho_1\rho_2\rho_6\mu_2)^{\mu_2}d_{D'_{1,0}}(x'_{3,i})}{\rho_7}<d_{D'_{1,0}}(x'_{3,i}),
$$
we see that
$$D'_{3,1}\subset D'_{1,0},
$$
whence it follows from
Theorem \Ref{ThmA} and the fact ``$D_{1,0}$ being $\rho_3$-uniform" that
$D_{3,1}$ is a $\rho(t)$-uniform domain, where $\rho(t)=\rho(n,
K, \rho_3, 2^{47+(3a\rho_7t)^2}a^2\rho_7t^2)$ and $\rho$ is from
Theorem \Ref{ThmA}.  Hence, Theorem  \Ref{Lem4-1} implies that
$D_{3,1}$ is $\varphi_{\rho(t)}$-broad, and then Theorem \Ref{Lem6-1}
shows that $D_{3,1}$ is $\mu_{0t}$-$LLC_2$ with respect to $\delta_{D_{3,0}}$ in $D_{3,0}$,
where $\mu_{0t}=\mu_3(n,\varphi_{\rho(t)})$ and $\mu_3$ is from Theorem \Ref{Lem6-1}.
Since Proposition \ref{00-11} and Theorem \Ref{Lem4-1} guarantee that
$D'_{3,0}$ is a $\varphi(t)$-broad domain, where $\varphi(t)=\varphi_{c_2(t)}$, $c_2(t)=2^{47+(3a\rho_7t)^2}a^2\rho_7t^2$ and $\varphi$ is
from Theorem \Ref{Lem4-1},
we see from Theorem
\Ref{Lem6-1'} that the restriction $f^{-1}|_{D'_{3,1}}$ is weakly
$H(t)$-quasisymmetric in the metrics $\delta_{D'_{3,0}}$ and $\delta_{D_{3,0}}$, where
$$H(t)=\kappa_2\left (n,K,\mu_{0t},\varphi(t), \frac{\delta_{D'_{3,0}}(D'_{3,1})}{d_{D'_{3,0}}(x'_{3,k_{3,0}})},
\frac{\delta_{D_{3,0}}(D_{3,1})}{d_{D_{3,0}}(x_{3,k_{3,0}})}\right ).
$$

Now, we need to estimate the ratios
$\frac{\delta_{D'_{3,0}}(D'_{3,1})}{d_{D'_{3,0}}(x'_{3,k_{3,0}})}$ and
$ \frac{\delta_{D_{3,0}}(D_{3,1})}{d_{D_{3,0}}(x_{3,k_{3,0}})}$
to
improve the constant $H(t)$. We first get an estimate for the second quantity.
%$ \frac{\delta_{D_{3,0}}(D_{3,1})}{d_{D_{3,0}}(x_{3,k_{3,0}})}$.
\medskip

For each $z'\in \mathbb{S}(x'_{3,k_{3,0}}, r_{3,k_{3,0}})$, we see from
Lemma \ref{mxll-1''} that
$$|x_{3,k_{3,0}}-z|< \frac{1}{5}d_D(x_{3,k_{3,0}}),
$$
and so (\ref{xxt}) implies
\beq\label{india-s-2}
k_D(x_{3,k_{3,0}}, z)\leq \log\Big(1+\frac{|x_{3,k_{3,0}}-z|}{d_D(x_{3,k_{3,0}})-|x_{3,k_{3,0}}-z|}\Big)<1.
\eeq
Moreover, by Lemma \ref{mxll-1} and Proposition \ref{00-11}, we have
\beq
\nonumber \frac{1}{\rho_7-1}\geq
k_{D'}(x'_{3,k_{3,0}},z')\geq
\log\Big(1+\frac{|x'_{3,k_{3,0}}-z'|}{d_{D'}(x'_{3,k_{3,0}})}\Big)\geq
\log\big(1+\frac{1}{3\rho_7}\big),
\eeq
whence it follows from \eqref{india-s-2} and Theorem \Ref{ThmF} that
\beq
\nonumber \Big(\frac{1}{6\mu_2\rho_7}\Big)^{\mu_2}&<&\Big(\frac{1}{\mu_2}\log(1+\frac{1}{3\rho_7})\Big)^{\mu_2}
\leq \big(\frac{1}{\mu_2}k_{D'}(x'_{3,k_{3,0}},z')\big)^{\mu_2}\leq
k_D(x_{3,k_{3,0}},z)\\ \nonumber  &\leq&
\log\Big(1+\frac{|x_{3,k_{3,0}}-z|}{d_D(x_{3,k_{3,0}})-|x_{3,k_{3,0}}-z|}\Big),
\eeq
and so
$$|x_{3,k_{3,0}}-z|\geq \Big(1-e^{-(\frac{1}{6\mu_2\rho_7})^{\mu_2}}\Big)d_D(x_{3,k_{3,0}}),
$$
which shows that
$$d_{f^{-1}(B_{3,k_{3,0}})}(x_{3,k_{3,0}})
>\Big(1-e^{-(\frac{1}{6\mu_2\rho_7})^{\mu_2}}\Big)d_D(x_{3,k_{3,0}}).
$$
Since $y'_0\in B_{3,k_{3,0}}$, it follows from
Lemma \ref{mxll-1''} that
\be\label{meiyouc}
\nonumber d_{D_{3,0}}(x_{3,k_{3,0}})=d_{f^{-1}(B_{3,k_{3,0}})}(x_{3,k_{3,0}})
>\frac{4}{5}\Big(1-e^{-(\frac{1}{6\mu_2\rho_7})^{\mu_2}}\Big)d_D(y_{0}),
\ee
whence we have
\be\label{xxmmm}
\frac{\delta_{D_{3,0}}(D_{3,1})}{d_{D_{3,0}}(x_{3,k_{3,0}})}\leq  \frac{\diam(D)}{d_{D_{3,0}}(x_{3,k_{3,0}})}
\leq \frac{5e^{(\frac{1}{6\mu_2\rho_7})^{\mu_2}}}{4(e^{(\frac{1}{6\mu_2\rho_7})^{\mu_2}}-1)}\cdot
\frac{\diam(D)}{d_D(y_{0})}< \rho_7^{2\mu_2}\cdot\frac{\diam(D)}{d_D(y_{0})}.
\ee

In order to get an estimate for $\frac{\delta_{D'_{3,0}}(D'_{3,1})}{d_{D'_{3,0}}(x'_{3,k_{3,0}})}$,
we let $v'_1\in \partial D'_{3,1}$ be such that
$$|v'_1-y'_0|\geq \frac{1}{2}\delta_{D'_{3,0}}(D'_{3,1})\;\;\mbox{ (see Figure \ref{fig13})}.
$$
%The existence of such a point $v'_2$ in $\partial D'_{3,1}$ is obvious.
Since $y'_0\in B_{3,k_{3,0}}$, we know from Proposition \ref{00-11} that
$$d_{D'}(y'_0)\leq d_{D'}(x'_{3,k_{3,0}})+|x'_{3,k_{3,0}}-y'_0|\leq (1+\frac{1}{\rho_7})d_{D'}(x'_{3,k_{3,0}}),
$$
and, since $D'$ has the $a$-carrot property with center $y'_0$, we see that
there is an $a$-carrot arc $\tau'$ in $D'$ connecting $v'_1$ and $y'_0$ with center $y'_0$.
Hence,
\beq\label{z-2}\nonumber
\delta_{D'_{3,0}}(D'_{3,1})&\leq& 2|v'_1-y'_0|\leq 2\ell(\tau')\leq 2ad_{D'}(y'_0)
\\ \nonumber&\leq&
2a(1+\frac{1}{\rho_7})d_{D'}(x'_{3,k_{3,0}})\leq 6a(1+\rho_7)r_{3,k_{3,0}}
\\
\nonumber&=& 6a(1+\rho_7)d_{D'_{3,0}}(x'_{3,k_{3,0}})
\eeq
so that
$$\frac{\delta_{D'_{3,0}}(D'_{3,1})}{d_{D'_{3,0}}(x'_{3,k_{3,0}})}\leq 6a(1+\rho_7).
$$
Hence we infer from (\ref{xxmmm}) that we can take $H(t)$ to be as follows:
$$H(t)= \kappa_2\Big(n,K,\mu_{0t},\varphi(t),6a(\rho_7+1), \rho_7^{2\mu_2}\frac{\diam(D)}{d_D(y_0)}\Big).
$$
%If $D$ and $D'$
%are unbounded, it follows from \eqref{sf-2} and \eqref{meiyouc} that
%\beq\label{xxmmm-1}\frac{\diam(D_{3,1})}{d_{D_{3,0}}(x_{3,k_{3,0}})}\leq
%\frac{5e^{(\frac{1}{6\mu_2\rho_7})^{\mu_2}}}{4(e^{(\frac{1}{6\mu_2\rho_7})^{\mu_2}}-1)}\cdot
%\frac{\diam(D_{1,0})}{d_D(y_{0})}\leq \rho_7^{2\mu_2}.\eeq

For each $i\in\{1,\ldots, k_{3,0}-1\}$, we get from $(2)$ in Proposition \ref{00-11} that
$$\frac{\rho_7-1}{\rho_7+1}d_{D'}(x'_{3,{i}})
\leq d_{D'}(x'_{3,{i+1}})\leq \frac{\rho_7+1}{\rho_7-1}d_{D'}(x'_{3,{i}}),
$$
from which the following easily follows:
$$\frac{\rho_7-1}{3(\rho_7+1)}r_{3,i}\leq r_{3,i+1}\leq \frac{3(\rho_7+1)}{\rho_7-1}r_{3,i},
$$
which guarantees that there is $t_3\in \{1,\ldots,t_2\}$ such that $x'_{3,t_3}\in
\gamma'_{2,0}[w'_1,x'_{3,t_2}]$ and
\be\label{mmeir-5-5}
\nonumber  \frac{1}{48}\delta_{D'_{3,0}}(x'_{3,t_2},w'_1)\leq
\delta_{D'_{3,0}}(x'_{3,t_3},w'_1)\leq
\frac{1}{2}\delta_{D'_{3,0}}(x'_{3,t_2},w'_1).
\ee
Obviously,
\be\label{mmeir-5-6}
\nonumber  \frac{1}{2}\delta_{D'_{3,0}}(x'_{3,t_2},w'_1)\leq\delta_{D'_{3,0}}(x'_{3,t_2},x'_{3,t_3})
\leq \frac{3}{2}\delta_{D'_{3,0}}(x'_{3,t_2},w'_1),
\ee
whence
\beq\label{d-man-1}
\frac{\delta_{D'_{3,0}}(x'_{3,t_3},w'_1)}{\delta_{D'_{3,0}}(x'_{3,t_2},x'_{3,t_3})}\leq 1.
\eeq

In order to apply the weak quasisymmetry of $f^{-1}$ and then to arrive at a contradiction,
we still need to estimate $\frac{\delta_{D_{3,0}}(x_{3,t_3},w_1)}{\delta_{D_{3,0}}(x_{3,t_2},x_{3,t_3})}$.

Let $v'_2$ be an intersection point of $\beta'_0$ with the ball $B_{3,t_2}$ (see Figure \ref{fig13}).
Since
$$|w'_2-v'_2|\leq \ell(\beta'_0[w'_2, v'_2])\leq ad_{D'}(v'_2),
$$
Lemma \ref{lem-4-1''} implies
$$d_{D'}(v'_2)\geq\frac{1}{2a}d_{D'}(w'_2),
$$
and so by Lemma \ref{lem5-h-0'} and Proposition \ref{00-11},
\begin{eqnarray*}
d_{D'_{1,0}}(x'_{3,t_2})&\geq& \frac{1}{(2^6\rho_1\rho_2\rho_6\mu_2)^{\mu_2}}d_{D'}(x'_{3,t_2})\\
&\geq & \frac{\rho_7}{(\rho_7+1)(2^6\rho_1\rho_2\rho_6\mu_2)^{\mu_2}}d_{D'}(v'_2)\\
\nonumber&>&
\frac{1}{3a(2^6\rho_1\rho_2\rho_6\mu_2)^{\mu_2}}d_{D'}(w'_2)\\ \nonumber&>&
\frac{1}{3a(2^6\rho_1\rho_2\rho_6\mu_2)^{\mu_2}}d_{D'_{1,0}}(w'_2),
\end{eqnarray*}
which, together with the assumption in the lemma, yields that
\beq\label{mmeir-5-3''}
k_{D'_{1,0}}(x'_{3,t_2},w'_1)&\geq& \log
\frac{d_{D'_{1,0}}(x'_{3,t_2})}{d_{D'_{1,0}}(w'_1)}\\ \nonumber
&\geq&
\log\frac{d_{D'_{1,0}}(w'_2)}{d_{D'_{1,0}}(w'_1)}-\mu_2\log(2^8a\rho_1\rho_2\rho_6\mu_2)
\\ \nonumber  &>&
\frac{1}{2}\log \Gamma(t),
\eeq
which, together with Theorem \Ref{ThmF}, shows that
\beq \label{ind-s-4}
k_{D_{1,0}}(x_{3,t_2},w_1)>1.
\eeq

Let us leave the proof of the lemma for a moment and prove the following inequality:
\beq \label{ind-s-5}
\diam(\gamma_{2,0}[z_{0,0},w_1])\leq 2^7\mu_8^2(3\rho_1+1)\rho_2\rho_5d_{D_{1,0}}(w_1).
\eeq

For a proof of this inequality, it follows from Corollary
\ref{Cor-1} that if
$$\min\{\diam(\gamma_{2,0}[z_{0,0},w_1]),\diam(\gamma_{2,0}[w_1,z_{1,0}])\}
=\diam(\gamma_{2,0}[z_{0,0},w_1]),
$$
then
$$\diam(\gamma_{2,0}[z_{0,0},w_1])\leq \mu_8d_{D_{1,0}}(w_1).
$$

For the remaining case
$$\min\{\diam(\gamma_{2,0}[z_{0,0},w_1]),\diam(\gamma_{2,0}[w_1,z_{1,0}])\}=\diam(\gamma_{2,0}[w_1,z_{1,0}]),
$$
Corollary \ref{Cor-1} shows that
$$\diam(\gamma_{2,0}[w_1,z_{1,0}])\leq \mu_8d_{D_{1,0}}(w_1),
$$
and by Lemma \ref{lem-4-1''},
$$d_{D_{1,0}}(z_{1,0})\leq 2\mu_8d_{D_{1,0}}(w_1).
$$
It follows from Corollary
\ref{Cor-1} and (\ref{mei-2}) that $\diam(\gamma_{2,0})\leq \mu_8|z_{0,0}-z_{1,0}|$, and thus
by Lemmas \ref{lem4-h-1} and \ref{lem4-h-0'}, we have
\beq \nonumber \label{mon-1}
\diam(\gamma_{2,0}[z_{0,0},w_1]) &\leq &
\diam(\gamma_{2,0})\leq \mu_8|z_{0,0}-z_{1,0}|
\\ \nonumber &\leq & \mu_8\rho_5d_{D_1}(z_{1,0})\leq 2^6\mu_8(3\rho_1+1)\rho_2\rho_5d_{D_{1,0}}(z_{1,0})
\\ \nonumber &\leq &2^7\mu_8^2(3\rho_1+1)\rho_2\rho_5d_{D_{1,0}}(w_1).
\eeq
Hence the inequality \eqref{ind-s-5} has been proved.
\medskip

Let us continue the discussion with the aid of \eqref{ind-s-5}.
Obviously, it follows from \eqref{ind-s-5} that
\beq\label{mmeir-5-4'}
|x_{3,t_2}-w_1|\leq \diam(\gamma_{2,0}[z_{0,0},w_1])\leq
2^7\mu_8^2(3\rho_1+1)\rho_2\rho_5d_{D_{1,0}}(w_1).
\eeq
Also we see from (\ref{mmeir-5-3''}), (\ref{ind-s-4}), Theorem \Ref{ThmF} and Corollary
\ref{Cor-1} that
\beq  \nonumber \label{mon-1'}
\log \Gamma(t) &\leq & 2k_{D'_{1,0}}(x'_{3,t_2},w'_1)\leq
2\mu_2k_{D_{1,0}}(x_{3,t_2},w_1)\\ \nonumber &\leq & 2\mu_2\mu_8
\log\Big(1+\frac{|x_{3,t_2}-w_1|}{\min\{d_{D_{1,0}}(x_{3,t_2}),d_{D_{1,0}}(w_1)\}}\Big).
\eeq
Necessarily, (\ref{mmeir-5-4'}) leads to
$$\min\{d_{D_{1,0}}(x_{3,t_2}),d_{D_{1,0}}(w_1)\}=d_{D_{1,0}}(x_{3,t_2}),
$$
and so, we have %\eqref{mon-1} implies
\be\label{z-4'}
|x_{3,t_2}-w_1|\geq \Big(\Gamma(t)^{\frac{1}{2\mu_2\mu_8}}-1\Big)d_{D_{1,0}}(x_{3,t_2}).
\ee

 It follows from
(\ref{india-s-1}) and Proposition \ref{00-11} that
\begin{eqnarray*}
\delta_{D'_{3,0}}(x'_{3,t_3},w'_1)&\leq&
t\min\{d_{D'_{1,0}}(x'_{3,t_3}),d_{D'_{1,0}}(x'_{3,t_2})\}\\
&\leq&
3\rho_7t\min\{d_{D'_{3,0}}(x'_{3,t_3}),d_{D'_{3,0}}(x'_{3,t_2})\}
,\end{eqnarray*}
and thus we obtain from Theorem \Ref{ThmF'} and (\ref{mmeir-5-5}) that
$$k_{D'_{3,0}}(x'_{3,t_2},x'_{3,t_3})\leq
\mu(t)\log\Big(1+\frac{|x'_{3,t_2}-x'_{3,t_3}|}{\min\{d_{D'_{3,0}}
(x'_{3,t_2}),d_{D'_{3,0}}(x'_{3,t_3})\}}\Big)<\mu(t)\log(1+147\rho_7t),
$$
since
$$|x'_{3,t_2}-x'_{3,t_3}|\leq |x'_{3,t_2}-w'_1|+|w'_1-x'_{3,t_3}|\leq 49\delta_{D'_{3,0}}(x'_{3,t_3},w'_1),
$$
where $\mu(t)=\mu_1(2^{47+(3a\rho_7t)^2}a^2\rho_7t^2)
$ and $\mu_1$ is from Theorem \Ref{ThmF'}. Then, Theorem \Ref{ThmF} shows that
$$\log\Big(1+\frac{|x_{3,t_2}-x_{3,t_3}|}{d_{D_{3,0}}(x_{3,t_2})}\Big)\leq k_{D_{3,0}}(x_{3,t_2},x_{3,t_3})
\leq \mu_2\mu(t)\log(1+147\rho_7t),
$$
and obviously,
$$|x_{3,t_2}-x_{3,t_3}|< (1+147\rho_7t)^{\mu_2\mu(t)}d_{D_{3,0}}(x_{3,t_2}).
$$
Hence it follows from the fact `` $D_{3,1}$ being a $\rho(t)$-uniform domain" that
\beq\label{0-2-1'}
\nonumber\delta_{D_{3,0}}(x_{3,t_2},x_{3,t_3})&\leq  & \delta_{D_{3,1}}(x_{3,t_2},x_{3,t_3})\leq
\rho(t)|x_{3,t_2}-x_{3,t_3}|\\ \nonumber
&\leq & \rho(t)(1+147\rho_7t)^{\mu_2\mu(t)}d_{D_{3,0}}(x_{3,t_2}),
\eeq
which, together with (\ref{z-4'}), shows that
\beq\label{meiren0-2-1'}\nonumber
\delta_{D_{3,0}}(x_{3,t_3},w_1)&\geq&
|x_{3,t_2}-w_1|-\delta_{D_{3,0}}(x_{3,t_3},x_{3,t_2})\\
\nonumber&\geq& \Big(\Gamma(t)^{\frac{1}{2\mu_2\mu_8}}
-\rho(t)(1+147\rho_7t){^{\mu_2\mu(t)}} -1\Big)d_{D_{3,0}}
(x_{3, t_2}),
\eeq
whence
$$\frac{\delta_{D_{3,0}}(x_{3,t_3},w_1)}{\delta_{D_{3,0}}(x_{3,t_2},x_{3,t_3})}\geq \frac{\Gamma(t)^{\frac{1}{2\mu_2\mu_8}}
-\rho(t)(1+147\rho_7t){^{\mu_2\mu(t)}} -1}{\rho(t)(1+147\rho_7t)^{\mu_2\mu(t)}}>H(t).
$$
But the weak quasisymmetry of $f^{-1}$ together with (\ref{d-man-1}) shows
that
$$\frac{\delta_{D_{3,0}}(x_{3,t_3},w_1)}{\delta_{D_{3,0}}(x_{3,t_2},x_{3,t_3})}\leq H(t).
$$
This obvious contradiction completes the proof of Lemma \ref{India-samy}.\medskip

As an application of Lemma \ref{India-samy}, we get the following result.

\begin{lem}\label{india-wang-0} There must exist some $i\in\{m_1,\ldots, m\}$ such that
$$k_{D'_{1,0}}(y'_{i+1},y'_i)\geq \frac{1}{2}\log \Lambda_{T_1+3},
$$
where $m_1=m-\big[\frac{\Gamma(\Lambda_{T_1+3})}{2}\big]$, $\Gamma(\Lambda_{T_1+3})$
and $\Lambda_{T_1+3}$ are defined as in Section \ref{sec-6}. Here $m$ is
defined by \eqref{india-mw-1} and we refer to Figure \ref{fig12} for the points $y'_i$ in $\gamma'_{2,0}$.
\end{lem}

Suppose on the contrary that
$$k_{D'_{1,0}}(y'_{i+1},y'_i)< \frac{1}{2}\log \Lambda_{T_1+3}
$$
for each $i\in\{m_1,\ldots,m\}$.

With the aid of Lemma \ref{India-samy}, we shall gain a contradiction by considering the arclength of the part
$\gamma'_{2,0}[y'_{m+1},w']$ for any $w'$ in $\gamma'_{2,0}[y'_{m+1},y'_{m_1}]$.
It follows from the obvious fact
$$\log\Big(1+\frac{\ell(\gamma'_{2,0}[y'_{i+1},y'_i])}{d_{D'_{1,0}(y'_{i+1})}}\Big)\leq k_{D'_{1,0}}(y'_{i+1},y'_i)
$$
that
$$\ell(\gamma'_{2,0}[y'_{i+1},y'_i])< \Lambda_{T_1+3}^{\frac{1}{2}}d_{D'_{1,0}}(y'_{i+1}),
$$
and further, for each $y'\in\gamma'_{2,0}[y'_{i+1},y'_i]$,
\be\label{india-wang-1}
\ell(\gamma'_{2,0}[y'_{i+1},y'_i])< \Lambda_{T_1+3}^{\frac{1}{2}}d_{D'_{1,0}}(y').
\ee

Since for every $w'\in\gamma'_{2,0}[y'_{m+1},y'_{m_1}]$, there exists some $i\in\{m_1,\ldots,m\}$
such that $w'\in\gamma'_{2,0}[y'_i,y'_{i+1}]$, we see from (\ref{india-wang-1}) that
\beq\label{in-man-1}
\ell(\gamma'_{2,0}[y'_{m+1},w'])&=  & \ell(\gamma'_{2,0}[y'_{m+1},y'_m])+
\cdots+\ell(\gamma'_{2,0}[y'_{i+2},y'_{i+1}])\\ \nonumber &&
+\ell(\gamma'_{2,0}[y'_{i+1},w'])\\ \nonumber
&\leq &\Lambda_{T_1+3}^{\frac{1}{2}}\big(d_{D'_{1,0}}(y'_{m+1})+\cdots+d_{D'_{1,0}}(y'_{i+2})
\\ \nonumber &&
+d_{D'_{1,0}}(y'_{i+1})\big)\\ \nonumber
&< &
2\Lambda_{T_1+3}^{\frac{1}{2}}d_{D'_{1,0}}(y'_{i+1})\\ \nonumber
&\leq & 2\Lambda_{T_1+3}^{\frac{1}{2}}d_{D'_{1,0}}(w'),
\eeq
and since
$$d_{D'_{1,0}}(y'_{m_1})=\frac{d_{D'_{1,0}}(y'_1)}{2^{m_1-1}}
=2^{\big[\frac{\Gamma(\Lambda_{T_1+3})}{2}\big]+1}d_{D'_{1,0}}(y'_{m+1})>\Gamma(\Lambda_{T_1+3})d_{D'_{1,0}}(y'_{m+1}),
$$
by taking  $w'_1=y'_{m+1}$, $w'_2=y'_{m_1}$ and $t=\Lambda_{T_1+3}$, we see from Lemma \ref{India-samy}
that there exists $y'\in \gamma'_{2,0}[y'_{m+1}, y'_{m_1}]$ such that
$$\ell(\gamma'_{2,0}[y'_{m+1},y'])\geq td_{D'_{1,0}}(y').
$$
A contradiction can be reached by taking $w'=y'$ in (\ref{in-man-1}) since, obviously,
$\Lambda_{T_1+3}> 2\Lambda_{T_1+3}^{\frac{1}{2}}$. The proof of the lemma is complete.
\medskip

It follows from Lemma \ref{india-wang-0} that there exists
$i_1\in\Big\{m-\big[\frac{\Gamma(\Lambda_{T_1+3})}{2}\big], \ldots, m\Big\}$
such that
\be\label{indian-1}
k_{D'_{1,0}}(y'_{i_1}, y'_{i_1+1})\geq \frac{1}{2}\log \Lambda_{T_1+3}.
\ee

Now, we prove a lower bound for the number $m$ below. Since Lemma \ref{14-2} implies
$$  \frac{\ell(\gamma'_{2,0})}{d_{D'_{1,0}}(\omega'_0)}\leq \int_{\gamma'_{2,0}}\frac{|dw'|}{d_{D'_{1,0}}(w')}
=k_{D'_{1,0}}(z'_{1,0}, z'_{0,0})\leq \rho_8\rho_9,$$ we see from (\ref{eq(h-3-3)}) and (\ref{24-1-mei}) that
\beq
\nonumber d_{D'_{1,0}}(\omega'_0)&\geq&
\frac{1}{\rho_8\rho_9}\ell(\gamma'_{2,0})\geq
\frac{1}{\rho_8\rho_9}|z'_{1,0}-z'_{0,0}|
\\ \nonumber
&>& 2^{\frac{1}{2}\rho_9\log\rho_8}\,d_{D'_{1,0}}(z'_{1,0}),
\eeq
and so \eqref{india-mw-1} leads to
\be\label{24-1-mei-1}
m>\frac{1}{2}\rho_9\log\rho_8-1.
\ee

Since $d_{D'_{1,0}}(y'_{i_1})=\frac{1}{2^{i_1-1}}d_{D'_{1,0}}(y'_{1})$, we know from \eqref{24-1-mei-1} that
\be\label{dear-wang-2-0}
d_{D'_{1,0}}(y'_{i_1})\leq \frac{1}{2^{m-\big[\frac{\Gamma(\Lambda_{T_1+3})}{2}\big]-1}}d_{D'_{1,0}}(y'_{1})
<\frac{1}{2^{5\Gamma(\Lambda_{T_1+3})}}d_{D'_{1,0}}(y'_{1}),
\ee
which implies $|y'_1-y'_{i_1}|>\frac{1}{2}d_{D'_{1,0}}(y'_1)$, and so
$$|y'_1-y'_{i_1}|> 2^{5\Gamma(\Lambda_{T_1+3})-1}d_{D'_{1,0}}(y'_{i_1}).
$$

\begin{figure}[!ht]
\begin{center}
\includegraphics%[width=10cm]
{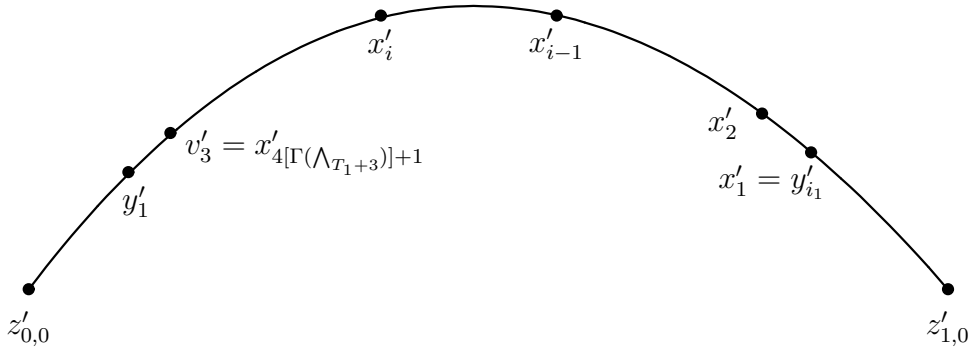}
\end{center}
\caption{A partition to $\gamma'_{2,0}[y'_{i_1}, v'_3]$ \label{figurep46}}
\end{figure}

Let $v'_3$ be the first point in $\gamma'_{2,0}[y'_{i_1},y'_{1}]$
from $y'_{i_1}$ to $y'_{1}$ such that (see Figure \ref{figurep46})
\beq\label{h-m-1}
\ell(\gamma'_{2,0}[y'_{i_1}, v'_3])=2^{4[\Gamma(\Lambda_{T_1+3})]}d_{D'_{1,0}}(y'_{i_1}).
%\;\;\mbox{(see Figure \ref{figurep46})}.
\eeq

Now, we partition the part $\gamma'_{2,0}[y'_{i_1},v'_3]$ of $\gamma'_{2,0}$ as follows.
We set $x_1'=y'_{i_1}$ and let $x'_2,\ldots ,x'_{4[\Gamma(\Lambda_{T_1+2})]+1}\in
\gamma'_{2,0}[y'_{i_1},v'_3]$ be the points such that for each $i\in
\{2,\ldots,4[\Gamma(\Lambda_{T_1+3})]+1\}$, $x'_i$ denotes the first point from
$x'_1$ to $v'_3$ with (see Figure \ref{figurep46})
%\be\label{hws-eq(4.3)}
\beq\label{m-w-l-s-3}
\ell(\gamma'_{2,0}[x'_1, x'_i])=2^{i-1}\, d_{D'_{1,0}}(x'_1)
%\;\;\mbox{(see Figure \ref{figurep46})}
.
\eeq
Obviously, $x'_{4[\Gamma(\Lambda_{T_1+3})]+1}=v'_3$. Then we have

%{\color{red} A figure for A partition to $\gamma'_{2,0}[y'_{i_1}, v'_3]$}

\begin{lem}\label{india-w-h-1}
For each pair $h$, $j\in\{1,\ldots,[\Gamma(\Lambda_{T_1+3})]+1\}$, if $j\geq h+[\Gamma(\Lambda_{i})]$ for some
$i\in\{1, \ldots, T_1+3\}$, then there exists some $s\in\{h+1,\ldots, j\}$ such that
$k_{D'_{1,0}}(x'_{s-1},x'_s)\geq \frac{1}{2}\log \Lambda_{i}$, where $\Lambda_{i}$ and $\Gamma(\Lambda_{i})$ are defined as in Section \ref{sec-6}.\end{lem}

Suppose on the contrary that
$$k_{D'_{1,0}}(x'_{s-1},x'_s)< \frac{1}{2}\log \Lambda_{i}
$$
for each $s\in\{h+1,\ldots, j\}$.

By considering the arclength of the part $\gamma'_{2,0}[x'_{h},w']$ for any $w'$ in
$\gamma'_{2,0}[x'_{h},x'_{j}]$ and by using Lemma \ref{India-samy}, we shall get a contradiction.
Since $\gamma'_{2,0}$ is a quasihyperbolic geodesic, we see from \eqref{h-w-01} that
for $x'\in \gamma'_{2,0}[x'_{s-1},x'_s]$,
$$\log\Big(1+\frac{\ell(\gamma'_{2,0}[x'_{s-1},x'_s])}{d_{D'_{1,0}}(x')}\Big)\leq k_{D'_{1,0}}(x'_{s-1},x'_s)<\frac{1}{2}\log \Lambda_{i},
$$
and so
\be\label{india-wang-2}
\ell(\gamma'_{2,0}[x'_{s-1},x'_s])< \Lambda_{i}^{\frac{1}{2}}d_{D'_{1,0}}(x').\
\ee

For $w'\in\gamma'_{2,0}[x'_{h},x'_{j}]$, obviously, there exists some $s_0\in\{h+1,\ldots,j\}$ such that
$w'\in\gamma'_{2,0}[x'_{s_0-1},x'_{s_0}]$. Then we see from (\ref{india-wang-2}) that
\beq\label{de-wang-0}
\ell(\gamma'_{2,0}[x'_{i},w'])&=  & \ell(\gamma'_{2,0}[x'_h,x'_{h+1}])+
\cdots+\ell(\gamma'_{2,0}[x'_{s_0-2},x'_{s_0-1}])\\ \nonumber &&+\ell(\gamma'_{2,0}[x'_{s_0-1},w'])
\\ \nonumber &\leq& \frac{1}{2}(\ell(\gamma'_{2,0}[x'_1,x'_{h+1}])+
\cdots+\ell(\gamma'_{2,0}[x'_1,x'_{s_0-1}])\\ \nonumber &&+\ell(\gamma'_{2,0}[x'_1,x'_{s_0}]))\\ \nonumber
&=&2\ell(\gamma'_{2,0}[x'_{s_0-1},x'_{s_0}])\\ \nonumber &<&
2\Lambda_{i}^{\frac{1}{2}}\,d_{D'_{1,0}}(w'),
\eeq
and consequently,
\be\label{de-wang-1}
\ell(\gamma'_{2,0}[x'_h,x'_j])<2\Lambda_{i}^{\frac{1}{2}}d_{D'_{1,0}}(x'_j).
\ee
Since
\beq\nonumber \ell(\gamma'_{2,0}[x'_{h},x'_j])
&=  & \ell(\gamma'_{2,0}[x'_1,x'_j])-\ell(\gamma'_{2,0}[x'_1,x'_h])= (2^{j-1}-2^{h-1})d_{D'_{1,0}}(x'_1)\\ \nonumber
&=& 2^{h-1}(2^{j-h}-1)\,d_{D'_{1,0}}(x'_1)
\eeq
and
$$d_{D'_{1,0}}(x'_h)\leq \ell(\gamma'_{2,0}[x'_1,x'_h])+d_{D'_{1,0}}(x'_1)\leq (2^{h-1}+1)d_{D'_{1,0}}(x'_1),
$$
we know from (\ref{de-wang-1}) that
$$\frac{d_{D'_{1,0}}(x'_j)}{d_{D'_{1,0}}(x'_h)}
\geq \frac{\ell(\gamma'_{2,0}[x'_h,x'_j])}{2\Lambda_{i}^{\frac{1}{2}}d_{D'_{1,0}}(x'_h)}
\geq \frac{2^{h-1}(2^{j-h}-1)\,d_{D'_{1,0}}(x'_1)}{2\Lambda_{i}^{\frac{1}{2}}(2^{h-1}+1)d_{D'_{1,0}}(x'_1)}
> 2^{\frac{\Gamma(\Lambda_{i})}{2}}>\Gamma(\Lambda_{i}).
$$
By taking  $w'_1=x'_h$, $w'_2=x'_j$ and $t=\Lambda_{i}$, Lemma \ref{India-samy} implies that
there exists $y'\in \gamma'_{2,0}[x'_h, x'_j]$ such that
$$\ell(\gamma'_{2,0}[x'_h,y'])\geq td_{D'_{1,0}}(y').
$$
A contradiction can be reached by taking $w'=y'$ in (\ref{de-wang-0}) since, obviously,
$\Lambda_{i}> 2\Lambda_{i}^{\frac{1}{2}}$.
 We complete the proof of Lemma \ref{india-w-h-1}.
\medskip

By using Lemma \ref{india-w-h-1}, we first choose some points
from $\{x'_{[\Gamma(\Lambda_{T_1+1})]+1}, \ldots, x'_{4[\Gamma(\Lambda_{T_1+3})]+1}\}$.

It follows from Lemma \ref{india-w-h-1} that
there exists $j\in\{[\Gamma(\Lambda_{T_1+1})]+1,\cdots,
2[\Gamma(\Lambda_{T_1+1})]\}$  such that
\be\label{india-5}
k_{D'_{1,0}}(x'_{j-1},x'_{j})\geq \frac{1}{2}\log\Lambda_{T_1+1}.
\ee
Let $x'_{j_1}$ denote the first point along the direction from $x'_{[\Gamma(\Lambda_{T_1+1})]+1}$
to $x'_{2[\Gamma(\Lambda_{T_1+1})]}$, which satisfies \eqref{india-5}. Easily, we have
$$\ell(\gamma'_{2,0}[x'_1,x'_{j_1}])\leq \ell(\gamma'_{2,0}[x'_1,x'_{2[\Gamma(\Lambda_{T_1+1})]}])=
\frac{2^{2[\Gamma(\Lambda_{T_1+1})]+1}}{4}d_{D'_{1,0}}(x'_1).
$$

Again, it follows from Lemma \ref{india-w-h-1} that there exists $j_2\in\{j_1+2,\ldots,
j_1+2\Gamma(\Lambda_{T_1-1})+1\}$  such that $x'_{j_2}$ is the first point along the
direction from $x'_{j_1+2}$ to $x'_{j_1+2\Gamma(\Lambda_{T_1-1})+1}$ satisfying
$$k_{D'_{1,0}}(x'_{j_2-1},x'_{j_2})\geq \frac{1}{2}\log\Lambda_{T_1-1}.
$$

Also, we get
$$\ell(\gamma'_{2,0}[x'_1,x'_{j_2}])\leq \ell(\gamma'_{2,0}[x'_1,x'_{j_1+2\Gamma(\Lambda_{T_1-1})+1}])\leq
\frac{2^{2\big([\Gamma(\Lambda_{T_1+1})]+\Gamma(\Lambda_{T_1-1})\big)+2}}{4}d_{D'_{1,0}}(x'_1).
$$

Once again, we see that there exists some $j_3\in\{j_2+2,\ldots,
j_2+2[\Gamma(\Lambda_{T_1-3})]+1\}$  such that $x'_{j_3}$ is the first point along the direction from $x'_{j_2+2}$ to $x'_{j_2+2[\Gamma(\Lambda_{T_1-3})]+1}$ satisfying
$$k_{D'_{1,0}}(x'_{j_3-1},x'_{j_3})\geq \frac{1}{2}\log\Lambda_{T_1-3}.
$$
Moreover, we obtain
\beqq
\ell(\gamma'_{2,0}[x'_1,x'_{j_3}])& \leq & \ell(\gamma'_{2,0}[x'_1,x'_{j_2+2[\Gamma(\Lambda_{T_1-3})]+1}])\\
&\leq &\frac{2^{2([\Gamma(\Lambda_{T_1+1})]+\Gamma(\Lambda_{T_1-1})+[\Gamma(\Lambda_{T_1-3})])+3}}{4}d_{D'_{1,0}}(x'_1).
\eeqq

\begin{figure}[!ht]
\begin{center}
\includegraphics%[width=10cm]
{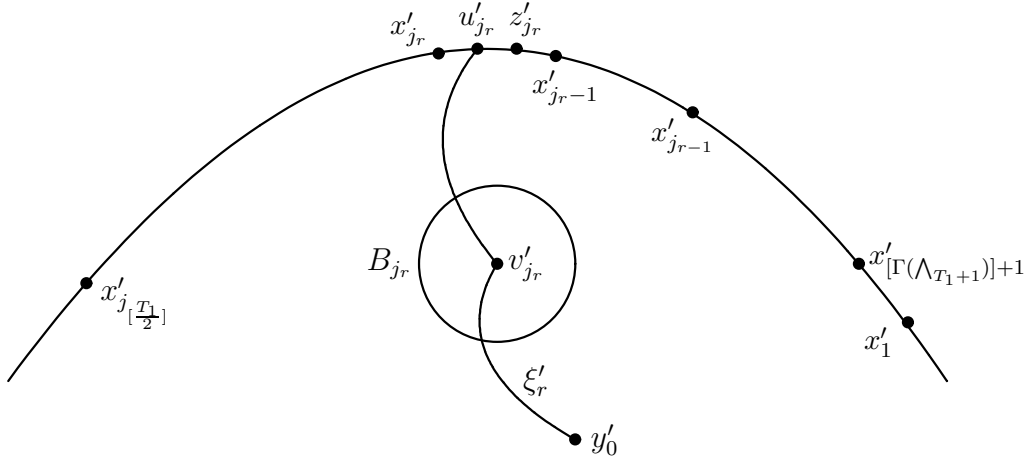}
\caption{Points $x'_{j_r}$ $(r\in \{1, \ldots, [\frac{T_1}{2}]\})$ in $\gamma'_{2,0}[y'_{i_1}, v'_3]$, $z'_{j_r}$
and $u'_{j_r}$ in $\gamma'_{2,0}[x'_{j_r-1}, x'_{j_r}]$, the carrot arcs $\xi'_{r}$ and the balls $B_{j_r}$ with center
$v'_{j_r}$ in $\xi'_{r}$ \label{figurep48}}
\end{center}
\end{figure}

By repeating this procedure as above $m$ $(m\geq 3)$ times, we shall determine a point $x'_{j_m}$ in $\gamma'_{2,0}$ such that
\beq\label{dear-man-1}
\ell(\gamma'_{2,0}[x'_1,x'_{j_m}])&\leq& \ell(\gamma'_{2,0}[x'_1,x'_{j_{m-1}+2[\Gamma(\Lambda_{T_1+3-2m})]+1}])\\
\nonumber&\leq& \frac{2^{\sum_{r=1}^{m}(2\Gamma(\Lambda_{T_1+3-2r})+1)}}{4} d_{D'_{1,0}}(x'_1),
\eeq
and \eqref{h-m-1} shows
\beq\label{dear-wang}
\ell(\gamma'_{2,0}[x'_1,x'_{j_{[\frac{T_1}{2}]}}])&\leq&
\frac{2^{\sum_{r=1}^{[\frac{T_1}{2}]}{(2\Gamma(\Lambda_{T_1+3-2r})+1)}}}{4}d_{D'_{1,0}}(x'_1)\\
\nonumber&<& 2^{T_1\Gamma(\Lambda_{T_1+1})+\frac{T_1}{2}-2}d_{D'_{1,0}}(x'_1)\\
\nonumber&<& 2^{\Gamma(\Lambda_{T_1+2})-2}d_{D'_{1,0}}(x'_1)\\
\nonumber&<& \ell(\gamma'_{2,0}[y'_{i_1}, v'_3]),
\eeq
whence we can find a sequence $\{x'_{j_r}\}_{r=1}^{[\frac{T_1}{2}]}$  such
that for each $r\in \{2, \ldots, [\frac{T_1}{2}]\}$, $x'_{j_r}$ denotes the first point along the
direction from $x'_{j_{r-1}+2}$ to $x'_{j_{r-1}+2[\Gamma(\Lambda_{T_1+3-2r})]+1}$ satisfying
\be\label{india-1-0}
k_{D'_{1,0}}(x'_{j_{r}-1},x'_{j_r})\geq \frac{1}{2}\log\Lambda_{T_1+3-2r} \;\;\mbox{(see Figure \ref{figurep48})}.
\ee

%{\color{red} A figure for points $x'_{j_r}$ $(r\in \{1, \cdots, [\frac{T_1}{2}]\})$ in $\gamma'_{2,0}[y'_{i_1}, v'_3]$, $z'_{j_r}$ and $u'_{j_r}$ in $\gamma_{2,0}[x'_{j_r-1}, x'_{j_r}]$, the carrot arcs $\xi'_{r}$ and the balls $B_{j_r}$ with center $v'_{j_r}$ in $\xi'_{r}$}

In order to determine some finite sequences of points in $\gamma'_{2,0}$,
we fix $r\in \{1, \ldots, [\frac{T_1}{2}]\}$ and emphasize our discussions on the part
$\gamma'_{2,0}[x'_{j_r-1},x'_{j_r}]$. It follows from (\ref{dear-man-1}) that
\beq\label{dear-man-2}
\ell(\gamma'_{2,0}[x'_1,x'_{j_{[\frac{T_1}{2}]}}])
&<& 2^{j_r+\sum_{i=r+1}^{[\frac{T_1}{2}]}\big(2\Gamma(\Lambda_{T_1+3-2i})+1\big)}d_{D'_{1,0}}(x'_1)\\
\nonumber &<& 2^{j_r+T_1\Gamma(\Lambda_{T_1+1-r})}d_{D'_{1,0}}(x'_1)\\
\nonumber &=&2^{1+T_1\Gamma(\Lambda_{T_1+1-r})}\ell(\gamma'_{2,0}[x'_1,x'_{j_r}]).
\eeq
By letting (see Figure \ref{figurep48})
$$d_{D'_{1,0}}(z'_{j_r})=\min_{z'\in\gamma'_{2,0}[x'_{j_r-1},x'_{j_r}]}d_{D'_{1,0}}(z'),
%\;\;\mbox{(see Figure \ref{figurep48})},
$$
we see that  $z'_{j_r}$ divides $\gamma'_{2,0}[x'_{j_r-1},x'_{j_r}]$ into two parts.
Without loss of generality, we may assume that
$$\max\{\ell(\gamma'_{2,0}[x'_{j_r-1},z'_{j_r}]),\ell(\gamma'_{2,0}[z'_{j_r},x'_{j_r}])\}=\ell(\gamma'_{2,0}[z'_{j_r},x'_{j_r}]).
$$
Let $u'_{j_r}\in\gamma'_{2,0}[z'_{j_r}, x'_{j_r}]$ be
the last point along the direction from $z'_{j_r}$ to $x'_{j_r}$ such that (see Figure \ref{figurep48})
\be\label{hws-z-1}
\ell(\gamma'_{2,0}[z'_{j_r},u'_{j_r}])\leq \Lambda_{T_1-2r}\;d_{D'_{1,0}}(u'_{j_r}).
%\;\;\mbox{(see Figure \ref{figurep48})}.
\ee
Then for each $z'\in \gamma'_{2,0}[u'_{j_r}, x'_{j_r}]$, we have
$$\ell(\gamma'_{2,0}[z'_{j_r},z'])\geq \Lambda_{T_1-2r}\,d_{D'_{1,0}}(z'),
$$
which, together with \eqref{hws-z-1}, implies
\beq
\nonumber d_{D'_{1,0}}(z')&\leq & \frac{1}{\Lambda_{T_1-2r}}\ell(\gamma'_{2,0}[z'_{j_r},z'])=
\frac{1}{\Lambda_{T_1-2r}}\big(   \ell(\gamma'_{2,0}[z'_{j_r},u'_{j_r}])+\ell(\gamma'_{2,0}[u'_{j_r}, z'])\big)
\\ \nonumber
 &\leq & d_{D'_{1,0}}(u'_{j_r})+\frac{1}{\Lambda_{T_1-2r}}\ell(\gamma'_{2,0}[u'_{j_r}, z']),
\eeq
whence
\beq\label{hws-z-4}
\ell_{k_{D'_{1,0}}}(\gamma'_{2,0}[u'_{j_r}, x'_{j_r}])&=&
\int_{\gamma'_{2,0}[u'_{j_r},
 x'_{j_r}]}\frac{|dz'|}{d_{D'_{1,0}}(z')}\\ \nonumber
 &\geq& \Lambda_{T_1-2r}\,
 \log\left(1+\frac{\ell(\gamma'_{2,0}[u'_{j_r}, x'_{j_r}])}{\Lambda_{T_1-2r}d_{D'_{1,0}}(u'_{j_r})}\right ).
\eeq

We take $\xi'_r$ to be an $a$-carrot arc joining $u'_{j_r}$ and $y'_0$ in $D'$ with center
$y'_0$ (see Figure \ref{figurep48}). The existence of $\xi'_r$ comes from the assumption
``$D'$ being an $a$-John domain with center $y'_0$".

Let $\lambda=j_{[\frac{T_1}{2}]}$. The continuation of the discussions needs a lower
bound for arclength of $\xi'_r$ which is as follows:
\beq\label{d-man-2}
\ell(\xi'_r)\geq 2^{\lambda+1}d_{D'_{1,0}}(x'_1).
\eeq
Since $u'_{j_r}\in \mathbb{\overline{B}}(x'_1, 2^{\lambda-1} d_{D'_{1,0}}(x'_1))$,
to get this lower bound, it suffices to show
$$y'_0 \notin\mathbb{\overline{B}}(x'_1, 2^{\lambda+1} d_{D'_{1,0}}(x'_1)).
$$
Suppose $|y'_0-x'_1|\leq 2^{\lambda+1} d_{D'_{1,0}}(x'_1)$.  Then Lemma \ref{lem5-h-0'} and (\ref{dear-wang}) imply
$$d_{D'}(y'_0)\leq |y'_0-x'_1|+d_{D'}(x'_1)\leq
(2^{\lambda+1}+(2^6\rho_1\rho_2\rho_6\mu_2)^{\mu_2})d_{D'_{1,0}}(x'_1)< 2^{\Gamma(\Lambda_{T_1+2})}d_{D'_{1,0}}(x'_1).
$$
But
it follows from the assumption ``$D'$ being an $a$-John domain with center $y'_0$" and (\ref{dear-wang-2-0})
that
$$ d_{D'}(y'_0)\geq \frac{1}{3a}\diam(D')
\geq \frac{2}{3a}d_{D'_{1,0}}(y'_1)
\geq \frac{2^{5\Gamma(\Lambda_{T_1+3})+1}}{3a}d_{D'_{1,0}}(y'_{i_1})> 2^{\Gamma(\Lambda_{T_1+2})}d_{D'_{1,0}}(x'_1).
$$
This obvious contradiction shows that \eqref{d-man-2} holds.

It follows from \eqref{d-man-2} that there is a point, denoted by $v'_{j_r}$ (see Figure \ref{figurep48}),
in the intersection $\xi'_r\cap \mathbb{B}(x'_1, 2^{\lambda+1} d_{D'_{1,0}}(x'_1))$ such that
\be\label{india-8-man}
\ell(\xi'_r[u'_{j_r},v'_{j_r}])=2^\lambda d_{D'_{1,0}}(x'_1).
\ee
Then
\be\label{india-8}
d_{D'}(v'_{j_r})\geq \frac{2^\lambda}{a}d_{D'_{1,0}}(x'_1).
\ee

Let $B_0=\mathbb{B}(x'_1, 2^{\lambda +2}d_{D'_{1,0}}(x'_1)), $
and for each $ r\in\{1,\ldots,[\frac{T_1}{2}]\}$, we let (see Figure \ref{figurep48})
$$B_{j_r}=\mathbb{B}\big(v'_{j_r}, \frac{2^{\lambda-1} }{a}d_{D'_{1,0}}(x'_1)\big).
%\;\;\mbox{(see Figure \ref{figurep48})}.
$$
We easily know that
$$B_{j_r}\subset B_0.
$$
Hence we have determined finite sequences of points $\{z'_{j_r}\}_{r=1}^{[\frac{T_1}{2}]}$ and
$\{u'_{j_r}\}_{r=1}^{[\frac{T_1}{2}]}$ in $\gamma_{2,0}$ together with a finite sequence of balls
$\{B_{j_r}\}_{r=1}^{[\frac{T_1}{2}]}$ in $B_0$. Similarly, we shall determine another finite sequence of points
$\{\zeta'_{j_{t_p}}\}_{p=1}^{p_0}$
in $\gamma'_{2,0}$ and another finite sequence of balls $\{B_{1p}\}_{p=1}^{p_0}$ in $B_0$.
The procedure is as follows:

Let
$$J=\Big\{r:\; r\in \{2,\ldots,[\frac{T_1}{2}]\}\;\;\mbox{and}\;\; j_r>j_{r-1}+2\Big\}.
$$
It is possible that $J=\emptyset$. If $J\not=\emptyset$, we let $J=\{t_1, \ldots, t_{p_0}\}$.
In particular, in what follows, we always assume that $p_0=0$ provided $J=\emptyset$. Then, obviously, $ p_0\leq [\frac{T_1}{2}]$.
If ${\rm card}\,(J)>1$, we always assume that $t_p<t_q$ if $p<q$.

\begin{figure}[!ht]
\begin{center}
\includegraphics%[width=10cm]
{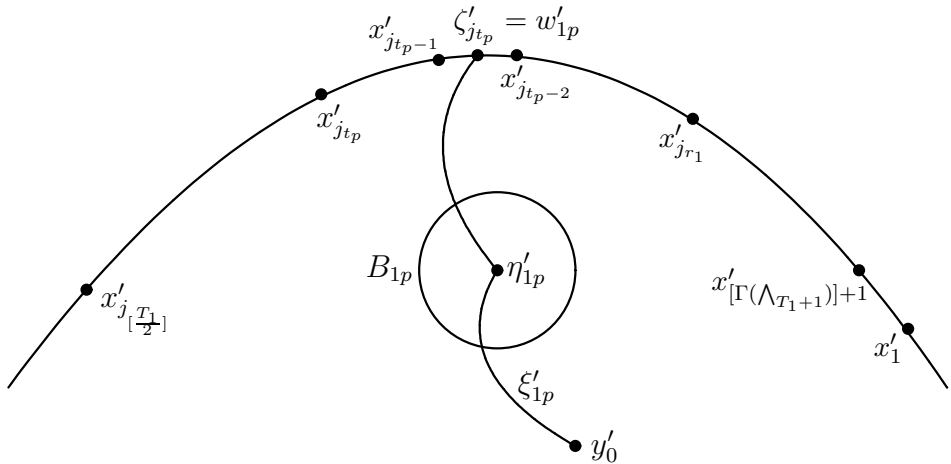}
\end{center}
\caption{Point $w'_{1p}$ in $\gamma'_{2,0}[x'_{j_{t_p}-2}, x'_{j_{t_p}-1}]$,
the carrot arcs $\xi'_{1p}$ and the balls $B_{1p}$ with center $\eta'_{1p}$ in $\xi'_{1p}$
\label{figurep50}}
\end{figure}

 For each $t_p\in J$, we let
\beq \label{manzi-1}
d_{D'_{1,0}}(\zeta'_{j_{t_p}})=\max_{z'\in\gamma'_{2,0}[x'_{j_{t_p}-2}, x'_{j_{t_p}-1}]}d_{D'_{1,0}}(z').
\eeq
In this way, we get a finite sequence of points
$$\{\zeta'_{j_{t_1}}, \ldots, \zeta'_{j_{t_{p_0}}}\}
$$
in $\gamma'_{2,0}$ (see Figure \ref{figurep50}).

%{\color{red} A figure for points $w'_{1p}$, $u'_{j_r}$ in $\gamma_{2,0}[x'_{j_r-1}, x'_{j_r}]$,
%the carrot arcs $\xi'_{1p}$ and the balls $B_{1p}$ with center
%$\eta'_{1p}$ in $\xi'_{1p}$}

For convenience, for each $p\in \{1, \ldots, p_0\}$, we denote the point $\zeta'_{j_{t_p}}$ by $w'_{1p}$,
and set
\beq \label{manzi-2}
k_{D'_{1,0}}(x'_{j_{t_p}-2},x'_{j_{t_p}-1})=\vartheta_{1{t_p}}.
\eeq
Obviously, the choice of $x'_{j_r}$ implies
\be\label{india-pwh-1}
\vartheta_{1{t_p}}<\frac{1}{2}\log\Lambda_{T_1+3-2{t_p}}.
\ee

For each $ p\in\{1,\ldots, p_0\}$,  let $\xi'_{1p}$ be an $a$-carrot arc joining $w'_{1p}$ and $y'_0$ in $D'$ with center
$y'_0$ (see Figure \ref{figurep50}). Denote by $\eta'_{1p}$ the point in the intersection
$\xi'_{1p}\cap \mathbb{B}(x'_1, 2^{\lambda+1} d_{D'_{1,0}}(x'_1))$  such that (see Figure \ref{figurep50})
$$\ell(\xi'_{1p}[w'_{1p},\eta'_{1p}])=2^\lambda d_{D'_{1,0}}(x'_1),
%\;\;\mbox{(see Figure \ref{figurep50})},
$$
and then
\be\label{manhuang-1}
d_{D'}(\eta'_{1p})\geq \frac{2^\lambda}{a}d_{D'_{1,0}}(x'_1).
\ee
Let
$$B_{1p}=\mathbb{B}\big(\eta'_{1p}, \frac{2^{\lambda-1}}{a}d_{D'_{1,0}}(x'_1)\big)\;\;\mbox{ (see Figure \ref{figurep50})}
$$
for each $ p\in\{1,\ldots, p_0\}$. Obviously,
\be\label{dear-india-1}B_{1p}\subset B_0.
\ee

\begin{figure}[!ht]
\begin{center}
\includegraphics%[width=10cm]
{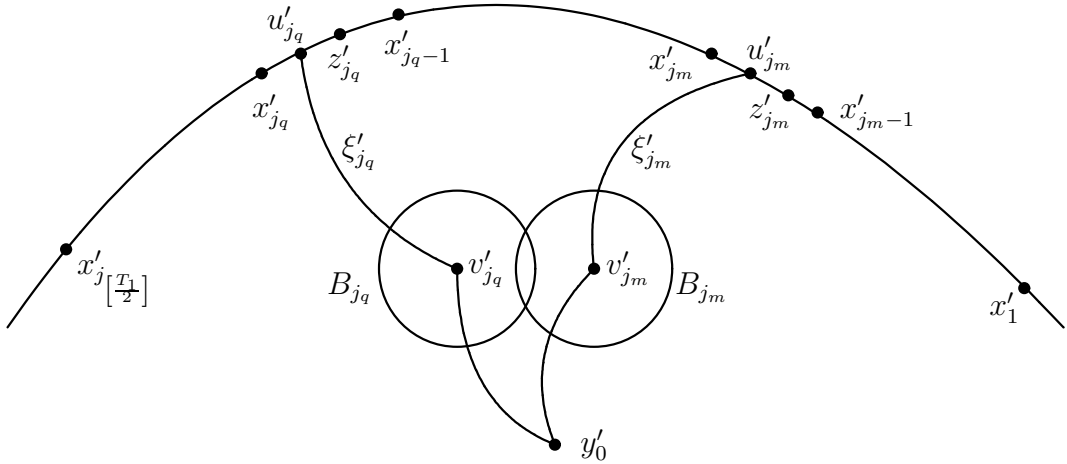}
\end{center}
\caption{The balls $B_{j_m}$ and $B_{j_q}$ \label{figurep51}}
\end{figure}

\noindent Then we have

\begin{lem}\label{india-pmw-1}
\bee
\item[$(1)$] Suppose that there exist integers $m<q\in \{1,\ldots, [\frac{T_1}{2}]\}$ such that
$B_{j_m}\cap B_{j_q}\not=\emptyset$ (see Figure \ref{figurep51}). Then
$$k_{D'}(u'_{j_{m}},u'_{j_q})
\leq 2a\bigg(\log\left (1+\frac{2^{\lambda+1}d_{D'_{1,0}}(x'_1)}{d_{D'}(u'_{j_m})}\right )
+\log\left (1+\frac{2^{\lambda+1}d_{D'_{1,0}}(x'_1)}{d_{D'}(u'_{j_q})}\right )+1\bigg);
$$

\item[$(2)$]  Suppose that there exist integers $h<s\in \{1,\ldots, p_0\}$ such that
$B_{1h}\cap B_{1s}\not=\emptyset$. Then
$$k_{D'}(w'_{1h},w'_{1s})
\leq 2a\bigg(\log\left (1+\frac{2^{\lambda+1}d_{D'_{1,0}}(x'_1)}{d_{D'}(w'_{1h})}\right )
+\log\left (1+\frac{2^{\lambda+1}d_{D'_{1,0}}(x'_1)}{d_{D'}(w'_{1s})}\right )+1\bigg).
$$
\eee
\end{lem}

%{\color{red} A figure for the balls $B_{j_m}$ and $B_{j_q}$}

To prove Lemma \ref{india-pmw-1}, we only need to prove $(1)$ since the proof of $(2)$ is similar.
%We will prove it by using the triangle inequalities: $k_{D'}(u'_{j_m},u'_{j_q})
%\leq k_{D'}(u'_{j_m},v'_{j_m})+k_{D'}(v'_{j_m}, v'_{j_q})+k_{D'}(u'_{j_q},v'_{j_q}).$

Since Lemma \ref{lem-4-1} and (\ref{india-8-man}) imply
\begin{eqnarray*}
k_{D'}(u'_{j_m},v'_{j_m}) & \leq & \ell_{k_{D'}}(\xi'_{j_m}[u'_{j_m},v'_{j_m}])
 \leq    2a\log\left (1+ \frac{2\ell(\xi'_{j_m}[u'_{j_m},v'_{j_m}])}{d_{D'}(u'_{j_m})}\right )\\
 & = &
2a\log\left (1+\frac{2^{\lambda+1}d_{D'_{1,0}}(x'_1)}{d_{D'}(u'_{j_m})}\right ),
\end{eqnarray*}
$$ k_{D'}(u'_{j_q},v'_{j_q})\leq 2a\log\left (1+\frac{2^{\lambda+1}d_{D'_{1,0}}(x'_1)}{d_{D'}(u'_{j_q})}\right ),
$$
and by \eqref{india-8}, obviously,
$$k_{D'}(v'_{j_m}, v'_{j_q})<2,
$$
we see that
\begin{eqnarray*}
k_{D'}(u'_{j_m},u'_{j_q})
& \leq & k_{D'}(u'_{j_m},v'_{j_m})+k_{D'}(v'_{j_m}, v'_{j_q})+k_{D'}(u'_{j_q},v'_{j_q})\\
& \leq &2a\bigg(\log\left (1+\frac{2^{\lambda+1}d_{D'_{1,0}}(x'_1)}{d_{D'}(u'_{j_m})}\right )+\log\left (1+\frac{2^{\lambda+1}d_{D'_{1,0}}(x'_1)}{d_{D'}(u'_{j_q})}\right )+1\bigg)
\end{eqnarray*}
as required.\medskip

Our next Lemma below is on the upper bound for $p_0$. Its proof needs
Propositions \ref{cl8-1-0} and \ref{scl-1-m-1} below.

\begin{prop} \label{cl8-1-0}
For all $x', y'\in\gamma'_{2,0}$, if $k_{D'_{1,0}}(x',y')\leq \vartheta$,
then there exists some point $z'\in\gamma'_{2,0}[x',y']$ such that
$\ell(\gamma'_{2,0}[x',y'])\leq \vartheta\;d_{D'_{1,0}}(z')$.
\end{prop}

For a proof of this result, we let $z'\in \gamma'_{2,0}[x',y']$ be such that
$$d_{D'_{1,0}}(z')=\sup_{w'\in D'_{1,0}}d_{D'_{1,0}}(w').
$$
Then
$$\frac{\ell(\gamma'_{2,0}[x',y'])}{d_{D'_{1,0}}(z')}
\leq\int_{\gamma'_{2,0}[x',y']}\frac{|dw'|}{d(w')}= k_{D'_{1,0}}(x',y')\leq \vartheta,
$$
which implies that
$$\ell(\gamma'_{2,0}[x',y'])\leq \vartheta\;d_{D'_{1,0}}(z')
$$
as required.
\medskip

\begin{prop}\label{scl-1-m-1}
For positive constants $r$, $\lambda_1$ and $\lambda_3$ with $\lambda_1>\lambda_3$, there
exists a positive integer $\lambda_2<\big(\frac{\lambda_1}{\lambda_3}\big)^n$ such that
the ball $\mathbb{B}(x, \lambda_1r)$ contains at most $\lambda_2$ disjoint
balls $B_i$, each of which has radius $\lambda_3r$.
\end{prop}

The proof easily follows from the following estimate:
$$\frac{\sum_{i=1}^{\lambda_2} {\rm Vol\,}(B_i)}{{\rm Vol\,}(\mathbb{B}(x, \lambda_1r))}
=\lambda_2\bigg(\frac{\lambda_3}{\lambda_1}\bigg)^n<1,
$$
where ``${\rm Vol\,}$" denotes the volume.

Now, we are in a position to state and prove our next Lemma.

\begin{lem} \label{wang-huang-1} $p_0\leq [8a]^n+2$.\end{lem}

Suppose on the contrary that $p_0 > [8a]^n+2$.
%there exist at least $[8a]^n$ points $w'_{1r}$ $(r\in\{1,\ldots, p_0\})$.
Under this assumption, (\ref{dear-india-1}) and Proposition \ref{scl-1-m-1} show that there exist
$p_1<p_2\in\{1,\ldots,p_0\}$ such that
$$B_{1p_1}\cap B_{1p_2}\not=\emptyset.
$$
By the construction of $w_{1p}$, we know that there  exist $t_{p_1}<t_{p_2}\in\{1,\ldots, [\frac{T_1}{2}]\}$ such that
$$\zeta'_{j_{t_{p_1}}}=w'_{1p_1}\;\;\mbox{and}\;\;\zeta'_{j_{t_{p_2}}}=w'_{1p_2}.
$$

%Obviously, by (\ref{india-1-0}),
%\be\label{manhuang-2}k_{D'_{1,0}}(\zeta'_{j_{t_{p_1}}},\zeta'_{j_{t_{p_2}}})>k_{D'_{1,0}}(x'_{j_{t_{p_1}}-1},x'_{j_{t_{p_1}}})\geq \frac{1}{2}\log \Lambda_{T_1+3-2t_{p_1}}.\ee
By estimating $k_{D'_{1,0}}(\zeta'_{j_{t_{p_1}}},\zeta'_{j_{t_{p_2}}})$, we shall get a contradiction.
First, we see from Proposition \ref{cl8-1-0} and (\ref{manzi-2}) that
$$\ell(\gamma'_{2,0}[x'_{j_{t_{p_1}}-2}, x'_{j_{t_{p_1}}-1}])\leq \vartheta_{1t_{p_1}}\,d_{D'_{1,0}}(\zeta'_{j_{t_{p_1}}})
$$
and
$$\ell(\gamma'_{2,0}[x'_{j_{t_{p_2}}-2}, x'_{j_{t_{p_2}}-1}])\leq \vartheta_{1t_{p_2}}\,d_{D'_{1,0}}(\zeta'_{j_{t_{p_2}}}).
$$

Further, by (\ref{dear-man-2}), we get
\begin{align*}
\ell(\gamma'_{2,0}[x'_{j_{t_{p_1}}-2}, x'_{j_{t_{p_1}}-1}])=&
\frac{1}{2}\ell(\gamma'_{2,0}[x'_1, x'_{j_{t_{p_1}}-1}])=\frac{1}{4}\ell(\gamma'_{2,0}[x'_1, x'_{j_{t_{p_1}}}])\\ \nonumber
\geq& \frac{1}{2^{3+T_1\Gamma(\Lambda_{T_1+1-2t_{p_1}})}}\ell(\gamma'_{2,0}[x'_1, x'_{\lambda}])\\ \nonumber
=&2^{\lambda-4-T_1\Gamma(\Lambda_{T_1+1-2t_{p_1}})}d_{D'_{1,0}}(x'_1),
\end{align*}
and similarly, we gain
$$\ell(\gamma'_{2,0}[x'_{j_{t_{p_2}}-2}, x'_{j_{t_{p_2}}-1}])\geq 2^{\lambda-4-T_1\Gamma(\Lambda_{T_1+1-2t_{p_2}})}d_{D'_{1,0}}(x'_1).
$$
Then (\ref{india-1-0}),  Lemma \ref{cl5-h-w-1} and Lemma \ref{india-pmw-1}(2) lead to
\begin{eqnarray*}
\mu_2^2\rho_7&<&\frac{1}{2}\log \Lambda_{T_1+3-2t_{p_1}}\\ \nonumber & \leq&
 k_{D'_{1,0}}(x'_{j_{t_{p_1}}-1},x'_{j_{t_{p_1}}})
\\ \nonumber&<& k_{D'_{1,0}}(\zeta'_{j_{t_{p_1}}},\zeta'_{j_{t_{p_2}}}) \leq  \mu_2^2\rho_7k_{D'}(\zeta'_{j_{t_{p_1}}},\zeta'_{j_{t_{p_2}}})\\
\nonumber
&\leq& 2a\mu_2^2\rho_7\bigg(\log\Big (1+\frac{2^{\lambda+1}d_{D'_{1,0}}(x'_1)}{d_{D'}(\zeta'_{j_{t_{p_1}}})}\Big )
+\log\Big (1+\frac{2^{\lambda+1}d_{D'_{1,0}}(x'_1)}{d_{D'}(\zeta'_{j_{t_{p_2}}})}\Big )+1\bigg)
%\\
\end{eqnarray*}
\begin{eqnarray*}
\nonumber & \leq& 2a\mu_2^2\rho_7\bigg(\log\Big (1+\frac{2^{\lambda+1}\vartheta_{1t_{p_1}}d_{D'_{1,0}}(x'_1)}{\ell(\gamma'_{2,0}[x'_{j_{t_{p_1}}-2}, x'_{j_{t_{p_1}}-1}])}\Big )\\ \nonumber
&&+\log\Big (1+\frac{2^{\lambda+1}\vartheta_{1t_{p_2}}d_{D'_{1,0}}(x'_1)}{\ell(\gamma'_{2,0}[x'_{j_{t_{p_2}}-2}, x'_{j_{t_{p_2}}-1}])}\Big )+1 \bigg )\\
\nonumber &\leq& 2a\mu_2^2\rho_7\bigg(\log\Big (1+\frac{2^{\lambda+1}\vartheta_{1t_{p_1}}d_{D'_{1,0}}(x'_1)}{2^{\lambda-4-T_1
\Gamma(\Lambda_{T_1+1-2t_{p_1}})}d_{D'_{1,0}}(x'_1)}\Big )\\ \nonumber
&&+\log\Big (1+\frac{2^{\lambda+1}\vartheta_{1t_{p_2}}d_{D'_{1,0}}(x'_1)}{2^{\lambda-4
-T_1\Gamma(\Lambda_{T_1+1-2t_{p_2}})}d_{D'_{1,0}}(x'_1)}\Big )+1\bigg)\\
\nonumber &<& 4a\mu_2^2\rho_7\Big(5+ T_1 \Gamma(\Lambda_{T_1+1-2t_{p_1}})+\log (1+\max\{\vartheta_{1t_{p_1}},\vartheta_{1t_{p_2}}\})\Big)+2a\mu_2^2\rho_7.
\end{eqnarray*}
But it follows from (\ref{india-pwh-1}) that
$$\max\{\vartheta_{1t_{p_2}},\vartheta_{1t_{p_1}}\}< \frac{1}{2}\log \Lambda_{T_1+3-2t_{p_1}},
$$
and so
$$4a\mu_2^2\rho_7\Big(5+ T_1 \Gamma(\Lambda_{T_1+1-2t_{p_1}})
+\log (1+\max\{\vartheta_{1t_{p_1}},\vartheta_{1t_{p_2}}\})\Big)+2a\mu_2^2\rho_7 < \frac{1}{2}\log \Lambda_{T_1+3-2t_{p_1}}.
$$
This obvious contradiction completes the proof of Lemma \ref{wang-huang-1}.
\medskip

\begin{lem}\label{india-samy-2-1} $k_{D'_{1,0}}(u'_{j_r},x'_{j_r})\leq 2\rho_8 k_{D'_{1,0}}(z'_{j_r}, u'_{j_r})$.
\end{lem}

Proof of  this Lemma easily follows from the upper bounds for
$$k_{D'_{1,0}}(x'_{j_r-1},x'_{j_r}) ~\mbox{ and }~\frac{d_{D'_{1,0}}(z'_{j_r})}{d_{D'_{1,0}}(u'_{j_r})}.
$$
First, we get an upper bound for $k_{D'_{1,0}}(x'_{j_r-1},x'_{j_r})$.
Since
$$k_{D'_{1,0}}(x'_{j_r-1},x'_{j_r})\geq \frac{1}{2}\log \Lambda_{T_1+3-2r}>8\mu_2^6\rho_7^4,
$$
we see from  Lemma \ref{cl5-h-w-1} that
$$k_{D'}(x'_{j_r-1},x'_{j_r})\geq \frac{1}{\mu_2^2\rho_7}k_{D'_{1,0}}(x'_{j_r-1},x'_{j_r})> 8\mu_2^4\rho_7^3.
$$
Then by Lemma \ref{cl8-0''}, we have
\be\label{dear-wmp-1}
\frac{1}{2}\log \Lambda_{T_1+3-2r}\leq k_{D'_{1,0}}(x'_{j_r-1},x'_{j_r})
\leq 48a^2\mu_2^6\rho_7^4\log \Big(1+\frac{\ell(\gamma'_{2,0}[x'_{j_r-1},x'_{j_r}])}{d_{D'_{1,0}}(z'_{j_r})}\Big),
\ee
whence
\be\label{dear-wm-1}
\ell(\gamma'_{2,0}[x'_{j_r-1},x'_{j_r}])>\Lambda_{T_1+2-2r}d_{D'_{1,0}}(z'_{j_r}),
\ee
which shows that there exists an integer $m_2$ such that
\be\label{dear-wm-3}
2^{m_2}\;d_{D'_{1,0}}(z'_{j_r})\leq \ell(\gamma'_{2,0}[x'_{j_r-1},x'_{j_r}])< 2^{m_2+1}\;d_{D'_{1,0}}(z'_{j_r}).
\ee
Hence it follows from (\ref{dear-wm-1}) that
\be\label{dear-wm-2}
m_2> \log_2 \Lambda_{T_1+2-2r}-1>\max\{\Lambda_{T_1+1-2r},3\log_2\Lambda_{T_1-2r},\Lambda_1,2\rho_8\}.
\ee

By \eqref{dear-wmp-1} and (\ref{dear-wm-3}), we have
\beq\label{india-yang-2}
k_{D'_{1,0}}(x'_{j_r-1},x'_{j_r})
&\leq& 48a^2\mu_2^6\rho_7^4\log\Big(1+\frac{\ell(\gamma'_{2,0}[x'_{j_r-1},x'_{j_r}])}{d_{D'_{1,0}}(z'_{j_r})}\Big)\\
\nonumber&\leq & 48a^2\mu_2^6\rho_7^4\log\Big(1+\frac{2^{m_2+1}\;d_{D'_{1,0}}(z'_{j_r})}{d_{D'_{1,0}}(z'_{j_r})}\Big)\\
\nonumber&=& 48a^2\mu_2^6\rho_7^4(1+2^{m_2+1})\\
\nonumber &<& 48a^2\mu_2^6\rho_7^4(m_2+1)\\ \nonumber&<&\frac{2\rho_8m_2}{3}.
\eeq

An upper bound for $\frac{d_{D'_{1,0}}(z'_{j_r})}{d_{D'_{1,0}}(u'_{j_r})}$ is as follows:
\be\label{india-y-1}
d_{D'_{1,0}}(u'_{j_r})\geq  2^{\frac{m_2}{2}+1}\;d_{D'_{1,0}}(z'_{j_r}).
\ee
Otherwise, (\ref{hws-z-1}) and (\ref{dear-wm-3}) imply that
\beq \nonumber
\ell(\gamma'_{2,0}[u'_{j_r}, x'_{j_r}])
&=& \ell(\gamma'_{2,0}[z'_{j_r}, x'_{j_r}])-\ell(\gamma'_{2,0}[z'_{j_r}, u'_{j_r}])  \\
\nonumber &\geq&\frac{1}{2}\ell(\gamma'_{2,0}[x'_{j_r-1},x'_{j_r}])- \Lambda_{T_1-2r}d_{D'_{1,0}}(u'_{j_r})\\
\nonumber &\geq& \frac{1}{2}\ell(\gamma'_{2,0}[x'_{j_r-1},x'_{j_r}])-2^{\frac{m_2}{2}+1}\Lambda_{T_1-2r}\;d_{D'_{1,0}}(z'_{j_r}) \\
\nonumber &\geq& (\frac{1}{2}-\frac{\Lambda_{T_1-2r}}{2^{\frac{m_2}{2}-1}})\ell(\gamma'_{2,0}[x'_{j_r-1},x'_{j_r}])\\
\nonumber &>& \frac{1}{4}\ell(\gamma'_{2,0}[x'_{j_r-1},x'_{j_r}]),
\eeq
since (\ref{dear-wm-2}) implies that $\frac{\Lambda_{T_1-2r}}{2^{\frac{m_2}{2}-1}}<\frac{1}{4}$. Then it
follows from (\ref{hws-z-4}) and (\ref{dear-wm-3}) that
\beq \nonumber
k_{D'_{1,0}}(u'_{j_r},x'_{j_r})&\geq & \Lambda_{T_1-2r}\log\Big(1+\frac{\ell(\gamma'_{2,0}[u'_{j_r},x'_{j_r}])}{\Lambda_{T_1-2r}d_{D'_{1,0}}(u'_{j_r})}\Big)\\
\nonumber&\geq& \Lambda_{T_1-2r}\log\Big(1+\frac{\ell(\gamma'_{2,0}[x'_{j_r-1},x'_{j_r}])}{4\Lambda_{T_1-2r}d_{D'_{1,0}}(u'_{j_r})}\Big) \\
\nonumber &\geq & \Lambda_{T_1-2r}\log\Big(1+\frac{2^{m_2}d_{D'_{1,0}}(z'_{j_r})}{2^{\frac{m_2}{2}+3}\Lambda_{T_1-2r}d_{D'_{1,0}}(z'_{j_r})}\Big) \\
\nonumber &\geq & \Lambda_{T_1-2r}\log\Big(1+\frac{2^{\frac{m_2}{2}-1}}{4\Lambda_{T_1-2r}}\Big),
\eeq
whence (\ref{dear-wm-2}) implies that
\beq \nonumber
\Lambda_{T_1-2r}\log\Big(1+\frac{2^{\frac{m_2}{2}-1}}{4\Lambda_{T_1-2r}}\Big) &\geq & 2\Lambda_{T_1-2r}^{\frac{1}{2}}\log(1+2^{\frac{m_2}{2}-1})\\
\nonumber &>& 48a^2\mu_2^6\rho_7^4(m_2+1),
\eeq
which, apparently, contradicts (\ref{india-yang-2}) since $k_{D'_{1,0}}(u'_{j_r},x'_{j_r})\leq k_{D'_{1,0}}(x'_{j_r-1},x'_{j_r})$.
Thus \eqref{india-y-1} holds.

Now, we obtain from (\ref{india-yang-2}) and (\ref{india-y-1}) that
$$k_{D'_{1,0}}(z'_{j_r}, u'_{j_r})
\geq \log\frac{d_{D'_{1,0}}(u'_{j_r})}{d_{D'_{1,0}}(z'_{j_r})}
> \frac{m_2}{3}> \frac{1}{2\rho_8}k_{D'_{1,0}}(x'_{j_r-1},x'_{j_r})
$$
as required.

%%%%%%%%%%%%%%%%%%%%%%%%%%%%%%%%%%%
%%%%%%%%%%%%%%%%%%%%%%%%%%%%%%%%%%%
\subsubsection{\bf Completion of the proof of Theorem \ref{thm6.1}}
%%%%%%%%%%%%%%%%%%%%%%%%%%%%%%%%%%%
%%%%%%%%%%%%%%%%%%%%%%%%%%%%%%%%%%%
To complete the proof, we first do some preparation. Since
$$\frac{[\frac{T_1}{2}]}{2[16a]^n}\geq 2[16a]^n>2([8a]^n+2),
$$
Lemma \ref{wang-huang-1} implies that there are $s_0$, $s_1\in \{1,\ldots, [\frac{T_1}{2}]\}$ with $s_1-s_0\geq 2[16a]^n$
such that $\gamma'_{2,0}[x'_{j_{s_0}}, x'_{j_{s_1}}]$ doesn't contain any $w'_{1p}$.

Since $s_1-s_0 > 2([8a]^n+2)$, by Proposition \ref{scl-1-m-1}, we know that  there exist integers
$p<q\in \{s_0,\ldots, s_1\}$ such that
$B_{j_p}\cap B_{j_q}\not=\emptyset$, $q-p\geq 2$ and $q\leq s_0+[16a]^n$.
It follows from (\ref{hws-z-4}) and Lemma \ref{india-samy-2-1} that
\beq\label{dear-man-5}
k_{D'_{1,0}}(u'_{j_p}, u'_{j_q})&>& k_{D'_{1,0}}(u'_{j_p}, x'_{j_p})+k_{D'_{1,0}}(z'_{j_q}, u'_{j_q})\\
\nonumber&>& \frac{1}{2\rho_8}\big(k_{D'_{1,0}}(u'_{j_p}, x'_{j_p})+k_{D'_{1,0}}(u'_{j_q},x'_{j_q})\big) \\
\nonumber&\geq& \frac{\Lambda_{T_1-2p}}{2\rho_8}\,\log\left(1+\frac{\ell(\gamma'_{2,0}[u'_{j_p}, x'_{j_{p}}])}{\Lambda_{T_1-2p}d_{D'_{1,0}}(u'_{j_p})}\right )\\
\nonumber&&+\frac{\Lambda_{T_1-2q}}{2\rho_8}\,\log\left(1+\frac{\ell(\gamma'_{2,0}[u'_{j_q}, x'_{j_{q}}])}{\Lambda_{T_1-2q}d_{D'_{1,0}}(u'_{j_q})}\right ).
\eeq

By (\ref{india-1-0}), we get
\beq\label{dear-man-6}
k_{D'_{1,0}}(u'_{j_p}, u'_{j_q})&>&k_{D'_{1,0}}(x'_{j_{p+1}-1}, x'_{j_{p+1}})\geq \frac{1}{2}\log \Lambda_{T_1+1-2p}\\
\nonumber&>&\Lambda_{T_1-2p}>\mu_2^2\rho_7,
\eeq
and then it follows from Lemmas \ref{cl5-h-w-1} and \ref{india-pmw-1} that
\beq \nonumber
k_{D'_{1,0}}(u'_{j_p},u'_{j_q})&\leq&\mu_2^2\rho_7k_{D'}(u'_{j_p},u'_{j_q})\\
\nonumber&\leq& 2a\mu_2^2\rho_7\bigg(\log\Big (1+\frac{2^{\lambda+1}d_{D'_{1,0}}(x'_1)}{d_{D'}(u'_{j_p})}\Big )
%\\\nonumber&&
+\log\Big (1+\frac{2^{\lambda+1}d_{D'_{1,0}}(x'_1)}{d_{D'}(u'_{j_q})}\Big )+1\bigg).
\eeq

Without loss of generality, we may assume that
$d_{D'}(u'_{j_p})\leq d_{D'}(u'_{j_q}),
$
and thus
\beq\label{dear-man-7}
k_{D'_{1,0}}(u'_{j_p},u'_{j_q}) \leq 4a\mu_2^2\rho_7\log\left (1+\frac{2^{\lambda+1}d_{D'_{1,0}}(x'_1)}{d_{D'}(u'_{j_p})}\right )+2a\mu_2^2\rho_7.
\eeq
Then we have
\be\label{dear-man-8}
2^{\lambda+1}d_{D'_{1,0}}(x'_1)\geq \Lambda_{T_1-2p}^2d_{D'}(u'_{j_p}),
\ee
because otherwise, (\ref{dear-man-7}) implies that
$$k_{D'_{1,0}}(u'_{j_p},u'_{j_q})< 4a\mu_2^2\rho_7\log(1+\Lambda_{T_1-2p}^2)+2a\mu_2^2\rho_7<\frac{1}{2}\log \Lambda_{T_1+1-2p},
$$
which contradicts (\ref{dear-man-6}). Hence \eqref{dear-man-8} holds.

Next, since
$$j_{p+[16a]^n}-j_p=2[16a]^n
~\mbox{ and }~\lambda-j_{p+[16a]^n}\leq \sum_{i=p+[16a]^n}^{[\frac{T_1}{2}]} (2\Gamma(\Lambda_{T_1+3-2i})+1),
$$
we have
$$\lambda\leq j_p+2[16a]^n+\sum_{i=p+[16a]^n}^{[\frac{T_1}{2}]} (2\Gamma(\Lambda_{T_1+3-2i})+1).
$$
Then elementary computations show that
\beq\label{dear-wang-xian-1}
\ell(\gamma'_{2,0}[x'_1,x'_{\lambda}])&=& 2^{\lambda-1}d_{D'_{1,0}}(x'_1)\\
\nonumber &\leq& 2^{j_p-1+2[16a]^n+\sum_{i=p+[16a]^n}^{[\frac{T_1}{2}]}(2\Gamma(\Lambda_{T_1+3-2i})+1)}d_{D'_{1,0}}(x'_1)\\
\nonumber &=& 2^{[\frac{T_1}{2}]+2[16a]^n+2\sum_{i=p+[16a]^n}^{[\frac{T_1}{2}]} \Gamma(\Lambda_{T_1+3-2i})}\ell(\gamma'_{2,0}[x'_1,x'_{j_p}])\\
\nonumber &<&2^{1+[\frac{T_1}{2}]+2[16a]^n+T_1 \,\Gamma(\Lambda_{T_1+3-2[16a]^n-2p})}\ell(\gamma'_{2,0}[x'_{j_p-1},x'_{j_p}])\\
\nonumber &<& \Lambda_{T_1-2p}\ell(\gamma'_{2,0}[x'_{j_p-1},x'_{j_p}]).
\eeq
%\medskip

%Now, we are ready to get the desired contradiction which completes the proof of Theorem \ref{thm6.1}.
%\medskip
\medskip

Now, we are in a position to finish the proof. To this end, we separate the discussions
into two cases. For the first case where
$$\ell(\gamma'_{2,0}[u'_{j_p},x'_{j_p}])\leq \frac{1}{4}\ell(\gamma'_{2,0}[x'_{j_p-1},x'_{j_p}]),
$$
we have
\begin{eqnarray*}
\ell(\gamma'_{2,0}[z'_{j_p},u'_{j_p}])&=&\ell(\gamma'_{2,0}[z'_{j_p},x'_{j_p}])-
\ell(\gamma'_{2,0}[u'_{j_p},x'_{j_p}])\\ \nonumber&\geq&\frac{1}{2}
\ell(\gamma'_{2,0}[x'_{j_p-1},x'_{j_p}])-\frac{1}{4}\ell(\gamma'_{2,0}[x'_{j_p-1},x'_{j_p}])
\\ \nonumber&\geq& \frac{1}{4}\ell(\gamma'_{2,0}[x'_{j_p-1},x'_{j_p}]),
\end{eqnarray*}
which, together with (\ref{hws-z-1}) and (\ref{dear-wang-xian-1}), shows that
\beq
\nonumber d_{D'_{1,0}}(u'_{j_p}) &\geq & \frac{\ell(\gamma'_{2,0}[z'_{j_p},u'_{j_p}])}{\Lambda_{T_1-2p}}
\geq \frac{\ell(\gamma'_{2,0}[x'_{j_p-1},x'_{j_p}])}{4\Lambda_{T_1-2p}}
\geq \frac{\ell(\gamma'_{2,0}[x'_1,x'_{\lambda}])}{4\Lambda_{T_1-2p}^2}\\
\nonumber &=& \frac{2^{\lambda-1}d_{D'_{1,0}}(x'_1)}{4\Lambda_{T_1-2p}^2},
\eeq
whence (\ref{dear-man-7}) implies
\be\label{dear-wang-2-4} \nonumber
k_{D'_{1,0}}(u'_{j_p}, u'_{j_q})< 4a\mu_2^2\rho_7\log(1+16\Lambda_{T_1-2p}^2)+2a\mu_2^2\rho_7<\frac{1}{2}\log \Lambda_{T_1+1-2p},
\ee
which contradicts (\ref{dear-man-6}).

 For the other case, namely,
$$\ell(\gamma'_{2,0}[u'_{j_p},x'_{j_p}])>\frac{1}{4}\ell(\gamma'_{2,0}[x'_{j_p-1},x'_{j_p}]),
$$
we infer from \eqref{dear-man-5} and \eqref{dear-wang-xian-1} that
\beq \label{In-1} \nonumber
k_{D'_{1,0}}(u'_{j_p}, u'_{j_q})
&>& \frac{\Lambda_{T_1-2p}}{2\rho_8}\,\log\Big(1+\frac{\ell(\gamma'_{2,0}[u'_{j_p}, x'_{j_{p}}])}{\Lambda_{T_1-2p}d_{D'_{1,0}}(u'_{j_p})}\Big )\\
\nonumber &\geq &
\frac{\Lambda_{T_1-2p}}{2\rho_8}\,\log\Big(1+
\frac{\ell(\gamma'_{2,0}[x'_{j_p-1},x'_{j_p}])}{4\Lambda_{T_1-2p}d_{D'_{1,0}}(u'_{j_p})}\Big )\\ \nonumber
&\geq &
\frac{\Lambda_{T_1-2p}}{2\rho_8}\,\log\Big(1+
\frac{\ell(\gamma'_{2,0}[x'_1,x'_{\lambda}])}{4\nu^2_{T_1-2p}d_{D'_{1,0}}(u'_{j_p})}\Big )\\ \nonumber
&= &
\frac{\Lambda_{T_1-2p}}{2\rho_8}\log\Big(1+\frac{2^{\lambda-1}d_{D'_{1,0}}(x'_1)}{4\Lambda_{T_1-2p}^2d_{D'_{1,0}}
(u'_{j_p})}\Big )
\\ \nonumber
&\geq &
\frac{\Lambda_{T_1-2p}}{2\rho_8}\log\Big(1+\frac{2^{\lambda-1}d_{D'_{1,0}}(x'_1)}{4\Lambda_{T_1-2p}^2d_{D'}
(u'_{j_p})}\Big ),
\eeq
whence (\ref{dear-man-7}) implies
$$\frac{\Lambda_{T_1-2p}}{2\rho_8}\log\left (1+\frac{2^{\lambda-1}d_{D'_{1,0}}(x'_1)}{4\Lambda_{T_1-2p}^2d_{D'}
(u'_{j_p})}\right )
 \leq 4a\mu_2^2\rho_7\log\left (1+\frac{2^{\lambda+1}d_{D'_{1,0}}(x'_1)}{d_{D'}(u'_{j_p})}\right )+2a\mu_2^2\rho_7.
 $$
Further, \eqref{dear-man-8} leads to
\beq \label{d-man-4}
\hspace{.6cm}\log\Big(1+\frac{t}{16\Lambda_{T_1-2p}^2} \Big )  < \frac{\rho_8^2}{\Lambda_{T_1-2p}}
\log (1+t ).
\eeq
where
$$t=\frac{2^{\lambda+1}d_{D'_{1,0}}(x'_1)}{d_{D'}(u'_{j_p})}.
$$
Then \eqref{dear-man-8} implies that $t\geq \Lambda_{T_1-2p}^2$. Finally, set
$$f(t)=\log\Big(1+\frac{t}{16\Lambda_{T_1-2p}^2}\Big )-
\frac{\rho_8^2}{\Lambda_{T_1-2p}} \log(1+t ).
$$
Elementary computations show that $f$ is increasing in $(\Lambda_{T_1-2p}^2,+\infty)$ and $f(\Lambda_{T_1-2p}^2)>0$,
which implies that \eqref{d-man-4} is impossible. This completes the proof of Theorem \ref{thm6.1}.
\qed

%%%%%%%%%%%%%%%%%%%%%%%%%%%%%%%%%%
\subsection{The proof of Theorem \ref{thm1.1}}
%%%%%%%%%%%%%%%%%%%%%%%%%%%%%%%%%%

It follows from Basic assumptions in Section \ref{sec-6},
Theorem \ref{thm6.1} and the arbitrariness of $z_1$ in $D_1$ that $D'_1$ is a
$2\rho_{10}$-John domain with center $z'_0\in D'_1$ in
the diameter metric. By Theorem \Ref{ThmF-1}, we see that Theorem
\ref{thm1.1} is true.
\qed

%\br
%From the proof, we see that, in several occasions, we get the desired contradictions by comparing the coefficients. %Hence the exact expression for $\rho_{10}$ is necessary.
%\er

%%%%%%%%%%%%%%%%%%%%%%%%%%%%%%%%%%
%%%%%%%%%%%%%%%%%%%%%%%%%%%%%%%%%%
\section{Quasisymmetries}\label{sec-7}
%%%%%%%%%%%%%%%%%%%%%%%%%%%%%%%%%%
%%%%%%%%%%%%%%%%%%%%%%%%%%%%%%%%%%

The aim of this section is to prove Theorem \ref{thm1.2}. The proof
easily follows from Lemma \ref{lem5-A-1} below. But the proof of
Lemma \ref{lem5-A-1} depends on Theorem \ref{thm1.1} together with a relationship between inner uniformity and
carrot property of the related domains whose precise statement is
as follows.

\begin{lem}\label{lem5-A-0'}
Suppose that  $D$ is a bounded $b$-inner uniform
domain. Then there exists some $x_0\in D$ such that $D$
has the $4b^2$-carrot property with center $x_0$.
\end{lem}

\bpf
Let $x_0\in D$ be such that
$$d_D(x_0)=\sup_{x\in D} d_D(x).
$$

Since $D$ is $b$-inner uniform, we see that for $y\in D$, there exists a curve $\alpha$ in $D$ connecting $y$ and $x'_0$
such that for every $x\in \alpha$,
$$\min\{\ell(\alpha[y,x]),\ell(\alpha[x,x_0])\}\leq bd_D(x).
$$
We only need to prove that
$$\ell(\alpha[y,x])\leq 4b^2d_D(x).
$$

For this, we consider two cases. For the first case where
$\min\{\ell(\alpha[y,x]),\ell(\alpha[x,x_0])\}=\ell(\alpha[y,x])$,
apparently,
$$\ell(\alpha[y,x])\leq b\,d_D(x).
$$
For the remaining case where
$\min\{\ell(\alpha[y,x]),\ell(\alpha[x,x_0])\}=\ell(\alpha[x,x_0])$,
we get
$$|x-x_0|\leq \ell(\alpha[x, x_0])\leq bd_D(x),
$$
whence Lemma \ref{lem-4-1''} implies
$$d_D(x_0)\leq 2bd_D(x),
$$
and so
$$\ell(\alpha[y,x])\leq \ell(\alpha)\leq 2bd_D(x_0)\leq 4b^2d_D(x).
$$
Hence the proof of Lemma \ref{lem5-A-0'} is complete. \epf

Now, we are ready to state and prove the key lemma in this section.

\begin{lem}\label{lem5-A-1}
Suppose that $f: D\to D'$ is a $K$-quasiconformal mapping between bounded
domains $D$, $D'\subset \IR^n$, where $D$ is a $b$-inner uniform
domain. Suppose also that $A\subset D$ is a pathwise connected set
and that $A'=f(A)$ has the $c_1$-carrot property in $D'$ with center
$y'_0\in D'$. Then $\diam(A)\leq c_2\,d_D(y_0)$
with $c_2=c_2\Big(n, K, b,c_1,\frac{\diam(D')}{d_{D'}(x'_0)}\Big)$,
where $x_0\in D$ is determined in Lemma \ref{lem5-A-0'}.
\end{lem}
%$(2)$  If  $D$ is unbounded and $f$ extends to a homeomorphism
%$D\cup \{\infty\}\to D'\cup \{\infty\}$, then $\diam(A)\leq
%^c_2\,d_D(y_0)$ with $c_2=c_2(n, K, b,c_1)$.
\bpf Since $y'_0\in D'$, we see that $y_0\not=\infty$ and
$y'_0\not=\infty$. Let $y_1\in A$ be such that
\be\label{mxpl-5-1}
|y_1-y_0|>\frac{1}{3}\diam(A),
\ee
and let $\alpha'_0$ be a $c_1$-carrot arc joining $y'_1$ and $y'_0$ in $D'$ with center
$y'_0$, i.e.,
\be\label{huang-1}\ell(\alpha'_0[y'_1,y'])\leq
c_1d_{D'}(y')
\ee
for all $y'\in \alpha'_0$. Then from Lemmas \ref{lem-m-v-1}
and \ref{lem2.2-4-mv-3}, we obtain

\bcl\label{xx-11-1} There exists a simply connected $b_0$-uniform
domain $A'_1=\bigcup \limits_{i=1}^{t_0}B_i\subset D'$ such that
\bee
\item[$(i)$] \label{equl--1}
$y'_0$, $y'_1\in A'_1$;
\item[$(ii)$] \label{equl--2}
for each $i\in \{1,\ldots, t_0\}$,
$$\frac{1}{3c_1}\,d_{D'}(x'_i)\leq r_i\leq \frac{1}{c_1}d_{D'}(x'_i);
$$
\item[$(iii)$] \label{equl--3}
if $t_0\geq 3$, then for all $i,j\in\{1,\ldots,t_0\}$ with
$|i-j|\geq 2$,
$$\dist(B_i, B_j)\geq \frac{1}{2^{37+c_1^4}}\max\{r_i,r_j\};
$$
\item[$(iv)$] \label{equl--4}
if $t_0\geq 2$, then
$$\ds r_i+r_{i+1}-|x'_i-x'_{i+1}|\ge \frac{1}{2^{37+c_1^4}}\max\{r_i,r_{i+1}\}
$$
for each $i\in\{1,\ldots,k_0-1\}$;
\item[$(v)$] \label{equl--4-11}
$d_{D'}(y'_0)\leq 2^{40+c_1^4}c_1d_{A'_1}(y'_0);$
\eee
where $B_i=\mathbb{B}(x'_i, r_i)$, $x'_i\in \alpha'_0$, $x'_i\not\in B_{i-1}$
for each $i\in \{2, \ldots, t_0\}$, $x'_1=y'_1$ and
$b_0=2^{44+c_1^4}c_1^3$.
\ecl

First, we prove an estimate on the diameter of $A'_1$ in terms of
the distance from $y'_0$ to the boundary of $A'_1$. If $t_0>1$, then
we know from (\ref{huang-1}) and Claim \ref{xx-11-1}$(v)$ that
$$\diam(A'_1)\leq 2\sum_{i=1}^{t_0}r_i\leq 4\ell(\alpha'_0)\leq
4c_1d_{D'}(y'_0) \leq 2^{42+c_1^4}c_1^2d_{A'_1}(y'_0).
$$
If $t_0=1$, then by $(ii)$ and $(v)$ in Claim \ref{xx-11-1},
$$\diam(A'_1)=2r_1\leq \frac{2}{c_1}d_{D'}(y'_1)\leq
\frac{2}{c_1-1}d_{D'}(y'_0)< 2^{42+c_1^4}d_{A'_1}(y'_0),
$$
since $d_{D'}(y'_0)\geq (1-\frac{1}{c_1})d_{D'}(y'_1)$. Hence we have
proved
\be\label{mxpl-5-1'-0} \diam(A'_1)\leq
2^{42+c_1^4}c_1^2d_{A'_1}(y'_0).
\ee
Next, we prove that $A'_1$ has the carrot property with
center $y'_0$. For each $y'\in A'_1$, let $\alpha'$ be a curve in $A'_1$ joining $y'$ and $y'_0$ such that
$$\min\{\ell(\alpha'[y',x']),\alpha'[y'_0,x'])\}\leq b_0d_{A'_1}(x')
$$
for all $x'\in \alpha'$, since $A'_1$ is $b_0$-uniform.

For the case
$\min\{\ell(\alpha'[y',x']),\alpha'[y'_0,x'])\}=\ell(\alpha'[y',x'])$, we
obviously have
$$\ell(\alpha'[y',x'])\leq b_0d_{A'_1}(x').
$$
For the other case
$\min\{\ell(\alpha'[y',x']),\alpha'[y'_0,x'])\}=\ell(\alpha'[x',y'_0])$,
by Lemma \ref{lem-4-1''}, we have
$$d_{A'_1}(y'_0)\leq 2b_0 d_{A'_1}(x'),
$$
and so (\ref{mxpl-5-1'-0}) implies
$$\ell(\alpha'[y',x'])\leq \ell(\alpha')\leq b_0\diam(A'_1)
\leq 2^{42+c_1^4}b_0c_1^2d_{A'_1}(y'_0)\leq 2^{43+c_1^4}(b_0c_1)^2
d_{A'_1}(x').
$$
Hence $A'_1$ has the $2^{43+c_1^4}(b_0c_1)^2$-carrot property with center $y'_0$.

Now, we are ready to finish the proof of the lemma. Since both $D$ and
$D'$ are bounded, it follows from the assumption ``$D$ being $b$-inner uniform" and Lemma \ref{lem5-A-0'}
that $D$ is a $4b^2$-John domain with center $x_0$, where $x_0$ satisfies
$$d_D(x_0)=\sup_{x\in D} d_D(x).
$$
Then Theorem \ref{thm1.1} implies that $A_1$ is $c_0$-John with center $y_0$, where
$$c_0=\tau\Big(n,K,b,c_1,\frac{\diam(D')}{d_{D'}(x'_0)}\Big).
$$
By (\ref{mxpl-5-1}), we have
$$\diam(A)\leq 3|y_1-y_0|\leq 3c_0d_{A_1}(y_0)\leq 3c_0d_D(y_0).
$$
Hence the proof of Lemma \ref{lem5-A-1} is complete.\epf

%{\color{red} The place for Figure 7:.}
%\begin{figure}[!ht]
%\includegraphics[width=0.70\textwidth]{figure07} %,height=0.4\textwidth
%\caption{The domain $D'_1$ and the arc
%$\alpha'_1$\label{fig7}}
%\end{figure}

%%%%%%%%%%%%%%%%%%%%%%%%%%%%%%%%%%
\subsection{ The proof of Theorem \ref{thm1.2}.}
%%%%%%%%%%%%%%%%%%%%%%%%%%%%%%%%%%

The proof of Theorem \ref{thm1.2} easily follows from  Theorems
\Ref{ThmE}, \Ref{Lem4-1} and Lemma \ref{lem5-A-1}. \qed

\bigskip

\noindent
{\bf Acknowledgements:}
The research was partly supported by NSFs of China (No. 11071063 and No. 11101138).
The authors would like to thank Professor Matti Vuorinen
for his useful suggestions. In particular, the inequality
\eqref{xxt} is suggested by him.
The revision of this work was completed during the visit of
the first and the fourth authors to the Department of Mathematics,
IIT Madras, Chennai-600 036, India.
These authors thank IIT Madras for the hospitality and the
other supports.

\end{document}